\documentclass{amsart}
\usepackage{amssymb}
\usepackage{graphics}

\usepackage[hyperref]{degt}

\makeatletter
\def\HyPsd@CatcodeWarning#1{}
\makeatother

\def\2{\color{magenta}}
\def\3{\color{red}}

\def\Dg:{\endgraf{\bf Dg:}\enspace\ignorespaces}

\def\pencil{\Cal P}

\let\ppencil\pencil

\def\Fano{\Cal F}
\def\Fn{\operatorname{Fn}}
\def\Span{\operatorname{span}}
\def\CS{\Cal S}

\def\pseudo{\Cal K}
\def\tS{\smash{\tilde S}}
\def\tP{\smash{\tilde P}}

\def\bG{\smash{\tilde G}}
\def\L{\bold L}
\def\mult{\operatorname{mult}}
\def\val{\operatorname{val}}
\def\fb{\operatorname\#}

\def\CG#1{\Z/#1}

\def\spcheck{\vee}

\def\impr{^\circ}
\def\prim{^\bullet}

\let\poly\xi

\let\leq\le

{\catcode`\*\active \catcode`\[\active
\gdef\deMaple{%
 \catcode`\*\active \catcode`\[\active
 \def*{}\def[##1]{_{##1}}}%
}
\let\e\epsilon

\let\Gm\mu
\let\Gn\nu
\def\CK{\Cal K}

\def\aD{\tilde D}

\def\lk{\operatorname{lk}}
\def\Ps{\operatorname{\frak{ps}}}
\def\Ls{\operatorname{\frak{ls}}}
\def\Pic{\operatorname{Pic}}

\def\G{\Bbb G}
\let\Gr\varrho
\let\Ge\epsilon
\let\Gu\upsilon
\let\Go\omega
\def\bs{\bar s}
\def\barS{\bar S}
\def\CU{\Cal U}
\def\CV{\Cal V}

\def\SV{\frak V}
\let\om\ominus
\let\op\oplus
\let\+\op
\def\*{^{\scriptscriptstyle\vee}}
\def\num#1(#2){\ls|{#2}|_{#1}}
\def\sdif#1(#2,#3){\{{#2\om#3}\}_{#1}}
\def\dif#1(#2,#3){\ls|{#2\om#3}|_{#1}}
\let\cs*
\def\scom#1(#2){\{{#2}\}_{#1}}
\def\com#1(#2){\ls|{#2}|_{#1}}
\def\one{\Bbb I}

\def\coord{\Cal C}
\def\cpre{\bold S}
\def\fib{\operatorname{fb}}
\def\sect{\frak S}
\def\bS{\bold S}
\def\bX{\bold X}
\def\bY{\bold Y}
\def\bZ{\bold Z}
\def\bQ{\bold Q}
\def\F{\Bbb F}
\def\CA{\Cal A}
\def\qq{\frak q}

\let\quarticbar\relax
\def\newconfig#1#2#3{%
 \expandafter\def\csname config=#1\endcsname{#2}%
 \expandafter\def\csname quartic=#1\endcsname{#3}}
\def\Xconfig#1#2{\newconfig{#1}{\bX#2}{X#2}}
\def\Yconfig#1#2{\newconfig{#1}{\bY#2}{Y#2}}
\def\Zconfig#1#2{\newconfig{#1}{\bZ#2}{Z#2\sp{t}}}
\def\Qconfig#1#2{\newconfig{#1}{\bQ#2}{Q#2}}
\def\theconfig#1{\hyperlink{config=#1}
 {\csname config=#1\endcsname}}
\def\confighook#1{\smash{\raise9pt\hbox{\hypertarget{config=#1}{}}}\csname config=#1\endcsname}
\pdef\aconfig#1{\ifmmode\confighook{#1}\else$\confighook{#1}$\fi}
\pdef\config#1{\ifmmode\theconfig{#1}\else$\theconfig{#1}$\fi}
\def\thequartic#1{\expandafter\expandafter\expandafter\quarticbar\csname quartic=#1\endcsname}
\pdef\quartic#1{\ifmmode\thequartic{#1}\else$\thequartic{#1}$\fi}
\pdef\bquartic#1{{\let\quarticbar\bar\quartic{#1}}}
\def\configno#1{^{#1}}
\def\configno#1{^{\scriptscriptstyle\mathrm{%
 \expandafter\uppercase\expandafter{\romannumeral#1}}}}
\def\configno#1{^{\mathrm{\romannumeral#1}}}
\def\q{*}
\def\pmin{p_{\min}}
\def\pmax{p_{\max}}
\def\vmin{v_{\min}}
\def\vmax{v_{\max}}
\def\ln{\ell}
\def\plane{\pi}

\makeatletter
\def\PS#1{\expandafter\@PSi#1\end@PS}
\def\@PSi#1[#2[#3]#4,#5]{(#3)^{#5}\futurelet\next\@PSii}
\def\@PSii{\ifx\next,\expandafter\@PSi\fi}
\let\end@PS\relax
\makeatother

\Xconfig{X}{}
\Yconfig{Y}{}
\Zconfig{Z}{}
\Qconfig{Q}{}
\Xconfig{64}{\sb{64}}
\Xconfig{60}{\sb{60}'}
\Xconfig{60.2}{\sb{60}''}
\Xconfig{56}{\sb{56}}
\Yconfig{56.real}{\sb{56}}
\Qconfig{q56}{\sb{56}}
\Xconfig{54}{\sb{54}}
\Qconfig{q54}{\sb{54}}
\Zconfig{52.0}{\sb{52}}
\Xconfig{52.1}{\sb{52}'}
\Xconfig{52.2}{\sb{52}''}
\Xconfig{52.5-3}{\sb{52}\configno5}
\Xconfig{52.4}{\sb{52}'''}
\Yconfig{52.5}{\sb{52}'}
\Yconfig{52.real}{\sb{52}''}
\Qconfig{q52.1}{\sb{52}'}
\Qconfig{q52.2}{\sb{52}''}
\Xconfig{51}{\sb{51}}
\Xconfig{50}{\sb{50}'''}
\Xconfig{50.1}{\sb{50}'}
\Xconfig{50.2}{\sb{50}''}
\Zconfig{50.0}{\sb{50}}
\Zconfig{49}{\sb{49}}
\Yconfig{48}{\sb{48}'}
\Yconfig{48.2}{\sb{48}''}
\Qconfig{q48.1}{\sb{48}'}
\Qconfig{q48.2}{\sb{48}''}
\Qconfig{q48.3}{\sb{48}'''}
\Qconfig{q48.4}{\sb{48}\configno4}
\Qconfig{q48.5}{\sb{48}\configno5}
\Qconfig{q48.6}{\sb{48}\configno6}
\Qconfig{q48.7}{\sb{48}\configno7}
\newconfig{33.1}{\bS\sb{33}'}{S\sb{33}'}
\newconfig{33.2}{\bS\sb{33}''}{S\sb{33}''}

\title{Lines on quartic surfaces}

\author{Alex Degtyarev}

\address{%
Bilkent University\\
Department of Mathematics\\
06800 Ankara, Turkey}

\thanks{%
The first author was supported by the JSPS grant L15517
and T\"{U}B\DOTaccent{I}TAK grant 114F325}

\email{degt@fen.bilkent.edu.tr}

\author{Ilia Itenberg}

\address{%
Institut de Math\'{e}matiques de Jussieu--Paris Rive Gauche\\
Universit\'{e} Pierre et Marie Curie\\
4 place Jussieu, 75252 Paris Cedex 5, France \\
and
D\'{e}partement de Math\'{e}mati\-ques et Applications,
Ecole Normale Sup\'{e}rieure \\
45 rue d'Ulm, 75230 Paris Cedex 5, France}

\thanks{%
The second author
was supported in part by
the FRG Collaborative Research grant
DMS-1265228 of the U.S. National Science Foundation}

\email{ilia.itenberg@imj-prg.fr}

\author{Al{\accent95\i} S{\accent95\i}nan Sert\"{o}z}

\address{%
Bilkent University\\
Department of Mathematics\\
06800 Ankara, Turkey}

\email{sertoz@bilkent.edu.tr}

\thanks{%
The third author was supported by the T\"{U}B\DOTaccent{I}TAK grant 114F325}

\keywords{%
$K3$-surface, elliptic pencil, integral lattice, discriminant form%
}

\subjclass[2000]{%
Primary: 14J28;
Secondary: 14J27, 14N25%
}

\begin{document}

\begin{abstract}
We show that the maximal number of (real) lines in a (real) nonsingular spatial
quartic surface is 64 (respectively, 56). We also give a complete projective
classification of all quartics containing more than 52 lines: all such
quartics are
projectively rigid.
Any value not exceeding 52 can appear as the number of lines of an
appropriate quartic.
\end{abstract}

\maketitle

\section{Introduction}

\subsection{Principal results}\label{s.results}
Throughout the paper,
all algebraic varieties are defined over~$\C$.
Given an algebraic surface $X\subset\Cp3$, we denote by $\Fn(X)$ the set of
projective lines contained in~$X$. If $X$ is real (see definition below),
$\Fn_\R(X)$ stands for the set of real lines contained in~$X$.

\theorem[see~\autoref{proof.64}]\label{th.unique}
Let $X\subset\Cp3$ be a nonsingular quartic, and assume that
$\ls|\Fn(X)|>52$.
Then $X$ is projectively equivalent to either
\roster*
\item
Schur's
quartic~$\quartic{64}$, see \autoref{s.Schur}, or
\item
one of the three
quartics~\quartic{60},
\quartic{60.2},
\bquartic{60.2}
described in \autoref{ss.th.unique}, or
\item
the quartic \quartic{56.real}, see \autoref{s.56}, or quartics \quartic{56},
\bquartic{56},
\quartic{q56} described in \autoref{ss.th.unique}, or
\item
one of the two quartics \quartic{54}, \quartic{q54}
described in \autoref{s.60}.
\endroster
In particular, one has
$\ls|\Fn(X)|=64$, $60$, $56$, or $54$, respectively.
\endtheorem

\corollary[see Segre~\cite{Segre} and Rams, Sch\"{u}tt~\cite{rams:2012}]\label{th.64}
Any nonsingular quartic in $\Cp3$ contains at most $64$ lines.
\done
\endcorollary

Note that the field of definition~$\C$ is essential for all statements. For
example, over $\F_9$, the quartic given by the equation
$z_0 z_3^3+z_1 z_2^3 + z_1^3 z_2+ z_0^3 z_3=0$
contains $112$ lines. According to Rams, Sch\"{u}tt~\cite{rams:2012},
the bound $\ls|\Fn(X)|\le64$ holds over any field of characteristic other
than~$2$ or~$3$.

As
was observed by T.~Shioda, $\quartic{56}$  and \bquartic{56} are
alternative projective models of
the Fermat quartic: this fact follows from the description of their
transcendental lattice, see \autoref{lem.unique}.
I.~Shimada has recently
found an explicit defining equation of these surfaces.
Other similar examples are discussed in \autoref{rem.abstract.K3}.

Recall that a \emph{real variety} is a complex algebraic variety~$X$ equipped
with a \emph{real structure}, \ie, an
anti-holomorphic involution $\conj\:X\to X$.
The \emph{real part} of~$X$ is
the fixed point set $X_\R:=\Fix\conj$. A subvariety (\eg, a line)
$Y\subset X$ is called \emph{real} if it is $\conj$-invariant.
When speaking about a \emph{real quartic} $X\subset\Cp3$, we assume that the
real structure on~$X$
is the restriction of the standard coordinatewise complex conjugation
$z\mapsto\bar z$ on $\Cp3$.

\corollary[see~\autoref{proof_th.56}]\label{th.56}
Let $X\subset\Cp3$ be a nonsingular \rom(over~$\C$\rom) real quartic,
and assume that $\ls|\Fn_\R(X)|>52$.
Then $X$ is projectively equivalent \rom(over $\R$\rom) to
the quartic $\quartic{56.real}$ given by \eqref{eq.56}.
In particular, one has $\ls|\Fn_\R(X)|=56$,
and this is the maximal number of real lines that can be contained in
a nonsingular real quartic.
\endcorollary

\addendum[see~\autoref{proof.counts}]\label{ad.counts}
For any number
\[*
n\in\{0,1,\ldots,51,52,54,56,60,64\},
\]
there exists a nonsingular
quartic $X\subset\Cp3$ such that $\ls|\Fn(X)|=n$.
For any
number
\[*
m\in\{0,1,\ldots,47,48,52,56\},
\]
there exists a nonsingular real
quartic $X\subset\Cp3$ such that $\ls|\Fn_\R(X)|=m$.
\endaddendum

Thus, for the moment
we are not certain about the values $\ls|\Fn_\R(X)|=49,50,51$.
We know three families of real quartics with $52$ real lines;
for a list of currently known large configurations of lines, see
\autoref{tab.list} in \autoref{obs.pencils}.

The quartic \quartic{56.real} can be defined over~$\Q$;
however, some of the lines are still defined only over $\Q(\sqrt2)$
(see \autoref{rem.Y.rational}).
At present, we do not know how many lines defined over~$\Q$ a quartic defined
over~$\Q$ may have; since $\Q\subset\R$ and
\quartic{56.real} has been ruled out,
\autoref{th.56} implies that this maximal number is at
most~$52$, the first candidates being the configurations
\config{52.5}, \config{52.real}, \config{52.0}.
Though, see \autoref{rem.min.field}.

Another open question is the maximal number of lines
contained in a triangle free
configuration, see \autoref{th.Segre} and \autoref{rem.triang.free}.

\subsection{Contents of the paper}
In \autoref{S.history},
we start with a brief introduction to the history of the subject.
In \autoref{S.reduction}, we recall basic notions and facts related to
integral lattices and $K3$-surfaces and use the theory of $K3$-surfaces to
reduce the original geometric problem to a purely arithmetical question
about \emph{configurations}; the
main results of this section are stated in \autoref{arithmetical_reduction}.
The simplest properties of
configurations, not related directly to
quartic surfaces, are treated in \autoref{S.geometry}, whereas
\autoref{S.arithm} deals with the more subtle arithmetic properties of the main
technical tool of the paper, the so-called \emph{pencils}.
The technical part is \autoref{S.counting}: we outline the algorithm used for
counting lines in a pair of obverse pencils and state the counts obtained in
the output. \autoref{tab.list} lists most known large configurations
of lines.
In \autoref{S.triang.free}, we digress to the so-called \emph{triangle free}
configurations, for which one can obtain a stronger bound
on the number of lines, see \autoref{th.Segre}.
The principal results of the paper stated in \autoref{s.results} are proved
in \autoref{S.pencils}. Finally, in \autoref{S.examples}, we discuss the
properties of quartics with many lines (in particular, \autoref{s.56}
contains an explicit equation of \quartic{56.real}) and make a few concluding
remarks.

\subsection{Acknowledgements}
A large part of the work on this project was accomplished during our visits to a
number of institutions worldwide. We are grateful to these organizations for
their hospitality and support:
\roster*
\item
\emph{\'{E}cole Normale Sup\'{e}rieure}
(first author),
\item
Hiroshima University, supported by the Japan Society for the Promotion of Science
(first author),
\item
\emph{Institut des Hautes \'{E}tudes Scientifiques}
(first and second authors),
\item
International Centre for Theoretical Physics
(first and second authors),
\item
\emph{Max-Planck-Institut f\"{u}r Mathematik}
(first and second authors),
\item
\emph{Universit\'{e} Pierre et Marie Curie} - Paris 6
(first author).
\endroster
We extend our gratitude to
Dmitrii Pasechnik,
S{\l}awomir Rams, Matthias Sch\"{u}tt, Ichiro Shimada,
Tetsuji Shioda, and
Davide Veniani for the motivation and fruitful discussions.

\section{History of the subject}\label{S.history}

The fact that there always exist exactly twenty seven lines on every smooth
cubic surface in the complex projective space $\PP^3$ naturally leads to
inquiries about higher degree  surfaces in $\PP^3$. The situation however
seems to be more involved for higher degree surfaces since it follows
immediately from a standard dimension count that a generic surface in $\PP^3$
of degree four or higher does not contain any lines at all, whereas each
Fermat surface of the form
\begin{equation*} z_0^d+z_1^d+z_2^d+z_3^d=0, \end{equation*}
where $[z_0:z_1:z_2:z_3]$ are homogeneous coordinates of $\PP^3$, contains
exactly $3d^2$ lines, for all $d\ge 3$. This then prompts the more relevant
question about how many lines a surface of degree $d\ge 4$ can have. In
particular, for a fixed $d\ge 4$, is there an upper bound for the number of
lines that a surface of degree $d$ can contain?

At this point it is appropriate to call attention to the difference between
the existence of rational curves on a surface and the existence of lines. A
line in $\PP^3$ is defined as the intersection of two hyperplanes but a
rational curve is an isomorphic image of $\PP^1$, which need not be a line.
Whereas we are expecting a finite number of lines on a surface the situation
is drastically different for the existence of rational curves.  While a
generic quartic does not contain a line, it is shown by Mori and Mukai
\cite{mori:1982} that every projective $K3$-surface, in particular every smooth
quartic in $\PP^3$, contains at least one rational curve. Moreover Bogomolov,
Hassett and Tschinkel showed in \cite{bogomolov:2011} that a generic
$K3$-surface, including a generic quartic surface in $\PP^3$,  contains infinitely
many rational curves. Going away from generic case to specific examples,
Bogomolov and Tschinkel showed in \cite{bogomolov:2000} that if a $K3$-surface
admitting an elliptic fibration has Picard number at most 19, then it
contains infinitely many rational curves.

Xi Chen showed in \cite{chen:1999} that for a generic quartic in $\PP^3$,
every linear system $\mathcal{O}(n)$, for any $n>0$, contains a nodal
rational curve. In fact Yau and Zaslow in \cite{yau:1996}, inspired by string
theory, counted those rational curves for the $n=1$ case.

Existence of smooth curves on quartic surfaces in $\PP^3$ is also relatively
well understood. Mori showed in \cite{mori:1984} that a quartic surface in
$\PP^3$ contains a smooth curve of degree $n>0$ and genus $g\ge 0$ if and
only if
either
\roster*
\item
$g=(n^2/8)+1$, or
\item
$g<(n^2/8)$ and $(n,g)\not=(5,3)$.
\endroster

The problem of counting lines on smooth surfaces in $\PP^3$ is on the other
hand a totally different game.

The first work which we can trace about this problem is Schur's article
\cite{schur:1882} where he exhibits a certain quartic surface which contains
64 lines. This surface is now known as Schur's quartic and is given by the
equation
\begin{equation*}
z_0(z_0^3-z_1^3)=z_2(z_2^3-z_3^3).
\end{equation*}
In \autoref{s.Schur} we give an account of the 64 lines on this quartic.

Apparently no progress was made on this result for about half a century until
1943 when Segre published some articles on the arithmetic and geometry of
surfaces in $\PP^3$. In one of these articles, in \cite{Segre}, he
claimed  that the number of lines which can lie on a quartic surface cannot
exceed 64. Since Schur's quartic already contains 64 lines, this result of
Segre would close the question for quartics were it not for a flaw in his
arguments which was only recently detected and corrected by Rams and Sch\"{u}tt
in \cite{rams:2012}. Rams and Sch\"{u}tt showed that the theorem is correct but
the proof needs some modifications using techniques which were not available
to Segre at that time.

Segre article \cite{Segre} contains an upper bound for the number of
lines which can lie on a surface of degree $d\ge 4$. His upper bound, which
is not affected by his erroneous argument about quartics, is $(d-2)(11d-6)$.
This bound is not expected to be sharp. For quartics it
predicts 76,
larger than the actual 64.

There is one curious fact about Segre's work of 1943. Most of the techniques
he uses were already in Salmon's book \cite{salmon:1862} which was originally
published in 1862. It would be reasonable to expect that a work similar to
Segre's be published much earlier than 1943. We learn from a footnote in
\cite{Segre} that the problem was mentioned by Meyer in an encyclopedia
article \cite{meyer:1908}  as early as 1908 but even that was not enough to
spur interest in the subject at the time.

After Segre's work there was again a period of long silence on the problem of
lines on surfaces. In 1983 Barth mentioned this problem in  \cite{barth:1983}
which turned out to be an influential manuscript on the subject. There he
also noted that since a smooth quartic in $\PP^3$ is a $K3$-surface and since
by Torelli theorems a $K3$-surface is nothing but its Picard
lattice, all
results of Segre on quartics could possibly be reproduced in the lattice
language. This teaser was one of the challenges which prompted us to work on this
problem
thirty
years later.

In 1995, Caporaso, Harris and Mazur, in \cite{caporaso:1995}, while
investigating  the number of rational points on a curve over an algebraic
number field, attacked the problem of finding a lower bound for the maximal
number $N_d$ of lines lying on a surface of the form
$\Gf(z_0,z_1)=\Gf(z_2,z_3)$, where $\Gf$ is a homogeneous form of degree
$d$. Their arguments being purely geometric, their findings made sense in the
complex domain. They found that in general for all $d\ge 4$,
\begin{equation*} N_d\ge 3d^2, \text{ but } N_4\ge 64, \quad N_6\ge 180, \quad N_8\ge 256, \quad N_{12}\ge 864, \quad N_{20}\ge 1600. \end{equation*}
Here the equality $N_4=64$ follows from Segre's work \cite{Segre}.

In 2006 Boissi\'{e}re and Sarti attacked this problem in \cite{boissiere:2007}
using group actions. They studied the maximal number of lines on symmetric
surfaces in $\PP^3$, where we called a surface symmetric
if its equation is of the
form
\begin{equation*}
\Gf(z_0,z_1)=\psi(z_2,z_3),
\end{equation*}
where $\Gf$ and $\psi$ are homogeneous forms of degree $d$, as
studied by
Caporaso, Harris and Mazur. This approach may seem restrictive at first;
nonetheless,
it is reasonable since Schur's surface which contains the maximal
possible number of lines a quartic surface can contain is itself of this
form. Boissi\'{e}re and Sarti first showed that for symmetric surfaces, the
inequalities about $N_d$ which Caporaso, Harris and Mazur obtained are
actually equalities. This increased the hope that the symmetric surfaces are
candidates to carry the most number of lines among other surfaces of the same
degree. However, Boissi\'{e}re and Sarti showed in the same work that this
expectation fails.  They showed that the non-symmetric surface given by
\begin{align*}  z_0^8&+z_1^8+z_2^8+z_3^8+168 z_0^2z_1^2z_2^2z_3^2 \\
                    &+14(z_0^4z_1^4+z_0^4z_2^4+z_0^4z_3^4+z_1^4z_2^4+z_1^4z_3^4+z_2^4z_3^4)=0
\end{align*}
contains 352 lines,
which is far greater than the upper bound of 256 for the symmetric
surfaces of the same degree. Notice that the number 352 is within the limits
allowed by Segre's upper bound, which gives 492  in this case.

Finally, almost thirty years after Barth's teaser, two teams started to work
on this problem, unaware of each other, from two different points of
approach. While we concentrated on understanding the ``lines on surfaces"
problem for $K3$-surfaces in $\PP^3$ and aimed at transliterating Segre's
results into the lattice language, Rams and Sch\"{u}tt decided to re-attack  the
problem by using elliptic fibration techniques
in \cite{rams:2012}.
They
discovered a flow in Segre's arguments which rendered his
proof void;
nonetheless, his theorem proved to be correctly stated.
Moreover,
Rams and Sch\"{u}tt's proof works on any
algebraically closed field of any characteristic
$p\not=2,3$.
Schur's quartic
becomes singular
when $p=2$ (still containing $64$ lines); when $p=3$,
it is shown in \cite{rams:2012}
that the surface contains 112 lines.

It is interesting to note that the concept of an elliptic fibration is
inevitable in studying the lines on a quartic. If $X$ is a smooth quartic in
$\PP^3$ and $L$ is a line lying on $X$,
one can parametrize the space of planes
$\Lambda_t$ in $\PP^3$ passing through $L$ by $t\in\PP^1$. Then any point
$p\in X$ determines a unique plane $\Lambda_t$, and the map sending $p$ to
$t$ is an elliptic fibration. If $p\in L$, we take $\Lambda_t$ as the plane
tangent to $X$ at~$p$. Segre starts with this observation but, using intuitive
geometric arguments,
he
erroneously claims that the maximal number of lines in the fibers of the
pencil is~$18$. The true bound is~$20$, see~\cite{rams:2012} or
\eqref{tablichka}, which calls for more work to establish the ultimate
bound~$64$ for the total number of lines in~$X$.

\section{The reduction}\label{S.reduction}

Throughout the paper, we consider various abelian groups $A$ equipped
with bilinear and/or quadratic forms. Whenever the form is fixed, we use
the abbreviation $x \cdot y$ (respectively, $x^2$)
for the value of the bilinear form on $x \otimes y$ (respectively, the quadratic form on $x$).
Given a subset $B \subset A$, its \emph{orthogonal complement} is $B^\perp =
\bigl\{x \in A \bigm|
\text{$x\cdot y = 0$ for all $y \in B$}\bigr\}.
$

\subsection{Integral lattices}\label{s.lattices}
An \emph{\rom(integral\rom) lattice} is a finitely generated free abelian group~$S$
supplied with a symmetric bilinear form $b\:S\otimes S\to\ZZ$.
A lattice~$S$ is
{\it even\/} if $x^2=0\bmod2$ for all $x\in S$.
As the transition matrix between two integral bases
has determinant $\pm1$, the
determinant $\det S\in\ZZ$
({\it i.e.}, the determinant
of the Gram matrix of~$b$ in
any
basis of~$S$)
is well defined.
A lattice~$S$ is called
{\it nondegenerate\/} if $\det S\ne0$; it is called {\it
unimodular\/} if $\det S=\pm1$.
Alternatively, $S$ is nondegenerate if and only if its \emph{kernel} $\ker S := S^\perp$
is trivial. An \emph{isometry} $\psi\: S \to S'$
between two lattices
is a group homomorphism respecting the bilinear forms;
obviously, one always has $\Ker \psi \subset \ker S$.
The group of auto-isometries of a
nondegenerate
lattice $S$ is denoted by $\OG(S)$.
Given a collection of subsets/elements $A_1,\ldots$ in~$S$, we
use the notation $\OG(S,A_1,\ldots)$ for the subgroup of $\OG(S)$
preserving each~$A_i$ as a set.

Given a lattice~$S$,
the bilinear form extends to $S\otimes\Q$ by linearity.
The inertia indices $\Gs_\pm S$, $\Gs_0 S$ and the signature $\Gs S$ of $S$
are defined as those of $S \otimes \Q$.
The orthogonal projection establishes a linear isomorphism between any two maximal positive definite subspaces of $S \otimes Q$,
thus providing a way for comparing their orientations.
A coherent choice of orientations of all
maximal positive definite subspaces is called a \emph{positive sign structure}.
Assuming $S$ nondegenerate, we denote by $\OG^+(S) \subset \OG(S)$ the subgroup
formed by the auto-isometries preserving a positive sign structure.

A {\it $d$-polarized lattice} is a lattice $S$ with a distinguished
vector $h \in S$,
referred to as the {\it polarization},
such that $h^2 = d$.
We use the abbreviation $\OG_h(A_1, \ldots)$ for $\OG(h, A_1, \ldots)$;
a similar convention applies for $\OG^+$.

If
$S$ is nondegenerate, the dual group $S^\spcheck=\Hom(S,\ZZ)$ can
be identified with the subgroup
$$
\bigl\{x\in S\otimes\Q\bigm|
 \text{$x\cdot y\in\ZZ$ for all $y\in S$}\bigr\}.
$$
In particular, $S\subset S^\spcheck$ and the quotient $S^\spcheck\!/S$
is a finite group; it is called the {\it discriminant group\/}
of~$S$ and is denoted by $\discr S$ or~$\CS$. The discriminant
group~$\CS$ inherits from $S\otimes\Q$ a symmetric bilinear form
$\CS\otimes\CS\to\Q/\ZZ$,
called the {\it discriminant form},
and, if $S$ is even, a quadratic
extension $\CS\to\Q/2\ZZ$ of this form.
When
speaking about the discriminant groups, their
(anti-)isomorphisms, \etc., we always assume that the discriminant
form (and its quadratic extension if the lattice is even) is taken
into account.
The number of elements in $\CS$ is equal to
$\mathopen|\det S\mathclose|$; in
particular, $\CS=0$ if and only if $S$ is unimodular.

Given a prime number $p$, we denote by $\CS_p$ or $\discr_p S$
the $p$-primary part of $\CS = \discr S$.
The form $\CS$ is called \emph{even} if there is no order $2$ element $\alpha \in \CS_2$
with $\alpha^2 = \pm \frac{1}{2} \bmod 2\Z$.
We use the notation $\ell(\CS)$ for the minimal
number of generators of~$\CS$, and we put $\ell_p(\CS) = \ell(\CS_p)$.
The quadratic form on $\CS$ can be described by means
of an analog $(\varepsilon_{ij})$
of the Gram matrix: assuming that $d_1 \mathbin | d_2 \mathbin | \ldots \mathbin | d_\ell$
are the invariant factors of $\CS$, we pick a basis $\alpha_1, \alpha_2, \ldots, \alpha_\ell \in \CS$
so that the order of $\alpha_i$ is $d_i$, and let $\varepsilon_{ij} = \alpha_i \cdot \alpha_j \bmod \Z$ for $i \ne j$
and $\varepsilon_{ii} = \alpha^2_i \bmod 2\Z$.
A similar construction applies to $\CS_p$.
Furthermore, according to R. Miranda and D. Morrison \cite{Miranda.Morrison:book},
unless $p = 2$ and $\CS_2$ is odd, the determinant of the resulting matrix is
a unit in $\Z_p$ well defined
modulo $(\Z^*_p)^2$; this determinant is denoted by
$\det_p \CS\in\Z_p^*/(\Z^*_p)^2$.

Two nondegenerate lattices are said to have the same {\it genus} if
their localizations
at all primes and at infinity
are pairwise
isomorphic.
The genus of
an even lattice is determined by its signature and the isomorphism
class of the quadratic extension of the discriminant
form,
see~\cite{Nikulin:forms}.

In what follows, we denote by $[s]$ the rank one lattice $\Z w$, $w^2 = s$.
The notation $\bU$ stands for the {\it hyperbolic plane},
{\it i.e.},
the lattice generated by a pair of vectors~$u$,~$v$
(referred to as a {\it standard basis\/} for~$\bU$)
with
$u^2=v^2=0$ and $u\cdot v=1$. Furthermore, given a lattice~$S$, we
denote by~$nS$, $n\in\NN$, the orthogonal direct sum of $n$~copies
of~$S$, and by~$S(q)$, $q\in\QQ$, the lattice obtained from~$S$ by
multiplying the form by~$q$ (assuming that the result is still an
integral lattice). The notation $n\CS$ is also used for the
orthogonal sum of $n$~copies of a discriminant group~$\CS$.

A {\it root\/} in an even lattice~$S$ is a vector $r\in S$ of
square~$-2$. A {\it root system\/} is an even negative definite
lattice generated by its roots. Recall that each root system
splits (uniquely up to order of the summands) into orthogonal sum
of indecomposable root systems, the latter being those of types
$\bA_n$, $n\ge1$, $\bD_n$, $n\ge4$, $\bE_6$, $\bE_7$, or~$\bE_8$,
see~\cite{Bourbaki:Lie:French}.

From now on, we fix an even unimodular lattice $\L$ of rank~$22$ and signature~$-16$.
All such lattices are isomorphic to $2\bE_8\oplus3\bU$.
It can easily be shown that, up to
the action $\OG^+(S)$,
this lattice has a unique $4$-polarization $h$;
thus, $\L$ is always considered equipped with a distinguished $4$-polarization $h$
and a positive sign structure.

We also fix the notation for certain discriminant forms.
Given coprime integers $m$, $n$ such that one of them is even,
$\<\frac{m}{n}\>$ is the quadratic form $1 \mapsto \frac{m}{n} \bmod 2\Z$
on $\CG{n}$.
Given a positive integer $k$, consider the group $\CG{2^k} \times \CG{2^k}$
generated by $\alpha = (1, 0)$ and $\beta=(0, 1)$;
denote by $\CU_{2^k}$ (respectively, $\CV_{2^k}$) the quadratic form
on the above group such that $\alpha \cdot \beta = \frac{1}{2^k} \bmod \Z$
and $\alpha^2 = \beta^2 = 0 \bmod 2\Z$ (respectively, $\alpha^2 = \beta^2 = \frac{1}{2^{k - 1}} \bmod 2\Z$).

An {\it extension\/} of a nondegenerate lattice~$S$ is another lattice~$M$
containing~$S$. An \emph{isomorphism} between two extensions $M', M'' \supset S$
is a bijective isometry $M' \to M''$ identical on $S$.
More generally, given a subgroup $G \subset \OG(S)$,
a $G$-isomorphism is a bijective isometry $M' \to M''$
whose restriction to $S$ is an element of $G$.

The two extreme cases are those of \emph{finite index} extensions
(where $S$ has finite index in $M$) and
\emph{primitive} ones
(where $M/S$
is torsion free).
The general case
$M\supset S$ splits into the finite index
extension $\tilde S\supset S$ and primitive extension
$M\supset\tilde S$, where
$$
\tilde S=\bigl\{x \in M\bigm|nx\in S\ \text{for some $n\in\ZZ$}\bigr\}
$$
is the {\it primitive hull\/} of~$S$ in~$M$.

If $S$ is nondegenerate and $M \supset S$ is a finite index extension, we have a chain of inclusions
$$
S \subset M \subset M^\spcheck \subset S^\spcheck,
$$
and, hence, a subgroup $\CK = M/S \subset \CS$;
this subgroup is called the \emph{pivot} of $M \supset S$.
The pivot $\CK$ is \emph{$b$-isotropic}, that is, the restriction to $\CK$
of the discriminant form $\CS \otimes \CS \to \Q/\Z$ is trivial.
Furthermore, the lattice $M$ is even if and only if $S$ is even
and $\CK$ is \emph{isotropic}, that is, the restriction to $\CK$
of the quadratic extension $\CS \to \Q/2\Z$ of the discriminant form is trivial.

\begin{theorem}[V. Nikulin \cite{Nikulin:forms}]\label{thm:Nik1} Given a nondegenerate
lattice~$S$,
the map sending $M \supset S$ to the pivot $\CK = M/S \subset \CS$
establishes a
one-to-one correspondence
between the set of isomorphism classes of finite index extensions of $S$
and the set of $b$-isotropic subgroups of $\CS$.
Under this correspondence, one has
$\discr M=\CK^\perp\!/\CK$ and
$M=\bigl\{x\in S^\spcheck\bigm|x\bmod S\in\CK\bigr\}$.
\end{theorem}

In the other extreme case, we confine ourselves to primitive extensions $M \supset S$
to an even unimodular lattice $M$. Assuming $S$
nondegenerate,
these are equivalent to
appropriate finite index extensions of $S \oplus S^\perp$, the pivot of the latter
giving rise to an anti-isomorphism $\CS \to \discr S^\perp$ and thus determining the genus of $S^\perp$.
It follows that, given a subgroup $G \subset \OG(S)$ and the signature of $M$,
a $G$-isomorphism class of even unimodular primitive extensions $M \supset S$
is determined by a choice of
\roster*
\item
an even lattice $T$ such that
$\discr T \cong -\CS$ and $\sigma_\pm T = \sigma_\pm M - \sigma_\pm S$, and
\item a bi-coset in $G\backslash\!\Aut\discr T/\!\OG(T)$.
\endroster
For details see \cite{Nikulin:forms}.
The following theorem is a combination of the above observation and Nikulin's existence theorem \cite{Nikulin:forms}
applied to the genus of $S^\perp$.

\theorem[V. Nikulin \cite{Nikulin:forms}]\label{th.Nikulin}
A nondegenerate even lattice $S$ admits a primitive extension to the lattice $\L$
if and only if the following conditions are satisfied\rom:
\roster
\item $\sigma_+ S \leq 3$, \ $\sigma_- S \leq 19$, and
$\rank S + \ell(\CS) \leq 22$\rom;
\item $(-1)^{\sigma_+ S - 1}\ls|\CS| = \det_p \CS \bmod (\Z^*_p)^2$ for all odd prime numbers $p$
such that  $\rank S + \ell_p(\CS) = 22$\rom;
\item either $\rank S + \ell_2(\CS) < 22$, or $\CS_2$ is odd, or
$\ls|\CS| = \pm \det_2 \CS \bmod (\Z^*_2)^2$.
\endroster
\endtheorem

\subsection{$K3$-surfaces}\label{s.K3}
Let $X\subset\Cp3$ be a nonsingular quartic.
It is a minimal $K3$-surface. Introduce the following
objects:
\roster*
\item $L_X=H_2(X)=H^2(X)$,
regarded as a lattice \via\ the intersection form
(we always identify homology and cohomology \via\ Poincar\'{e}
duality);
\item
$h_X\in L_X$, the class of a generic plane section
of~$X$;
\item
$\Fano(X) \subset H_2(X; \Z)$,
the primitive sublattice
spanned over $\Q$ by $h_X$ and the classes of lines $l \subset X$
(the {\it Fano configuration}
of $X$);
\item
$\Go_X\subset L_X\otimes\R$, the oriented $2$-subspace
spanned by the real and imaginary parts of the class of a
holomorphic $2$-form on~$X$ (the \emph{period} of~$X$).
\endroster
Note that $\Go_X$ is positive definite and orthogonal to $h_X$;
furthermore, the Picard group $\Pic X$ equals $\omega_X^\perp \cap L_X$.


The following statement is an immediate consequence of the above description of $\Pic X$
and the Riemann--Roch theorem.

\lemma\label{line_classes}
A vector $a \in L_X$ is realized by a line $l\subset X$
if and only if $a \cdot \omega_X = 0$, $a^2 = -2$, and $a \cdot h_X = 1$.
Distinct lines represent distinct classes in $L_X$.
\done
\endlemma

In view of the uniqueness part of this statement, we identify lines in $X$
and their classes in $L_X$.

As is well known, the lattice $L_X$ is isomorphic to $\L$;
a \emph{marking} of $X$ is a choice of
a particular isomorphism $\psi\: L_X \to \L$
such that $\psi(h_X) = h \in \L$
and the maximal positive definite subspace $\psi(\R h_X \oplus \omega_X)$
is positively oriented.
Consider a period $\omega$, {\it i.e.}, an oriented positive definite
$2$-subspace $\omega \subset \L \otimes \R$
orthogonal to $h$.
The following statement provides a criterion for the realizability of the
triple $(\L, h, \omega)$
by a quartic, \ie, the existence of a marked nonsingular quartic $(X, \psi)$
such that $\psi$ takes $\omega_X$ to $\omega$.
It is a combination of the surjectivity of the period map for $K3$-surfaces
(see Vik\. Kulikov
\cite{Kulikov:periods})
and Saint-Donat's description~\cite{Saint-Donat} of projective models of $K3$-surfaces.

\proposition\label{Saint-Donat}
A triple $(\L ,h , \omega)$ is realizable by a quartic $X \subset \Cp3$ if and only if
$\L$ contains no vector $e$ such that $e \cdot \omega = 0$ and either
\roster
\item\label{ex-div1} $e^2 = -2$ and $e \cdot h= 0$, or
\item\label{ell-pencil1} $e^2 = 0$ and $e \cdot h = 2$.
\pni
\endroster
\endproposition

Denote by $\Omega$ the space of oriented positive definite $2$-subspaces $\omega \subset \L \otimes \R$
orthogonal to $h$ and such that $\R h \oplus \omega$ is positively oriented.
By \autoref{Saint-Donat}, the image of the \emph{period map}
$(X, \psi) \mapsto \psi(\omega_X)$
is the subset $\Omega^\circ \subset \Omega$ obtained by removing the locally finite collection
of codimension two subspaces
\[*
\Omega_e = \{\omega \in \Omega \; | \;  \omega \cdot e = 0 \},
\]
where $e \in \L$ runs over all vectors as in \autoref{Saint-Donat}\iref{ex-div1} or \iref{ell-pencil1}.
Restricting to $\Omega^\circ$ Beauville's universal family \cite{Beauville:moduli} of marked polarized $K3$-surfaces,
we obtain the following statement on marked quartics.

\proposition\label{moduli}
The subset $\Omega^\circ \subset \Omega$ is a fine moduli space of marked nonsingular quartics in $\Cp3$. \pni
\endproposition

Now, let $X \subset \Cp3$ be a real nonsingular quartic.
The complex conjugation induces an involutive isometry $c_X\: L_X \to L_X$
taking $h_X$ to $-h_X$, preserving $\omega_X$ as a subspace and reversing its orientation.
In particular, it follows that the positive inertia index of the skew-invariant eigenlattice of $c_X$
equals $2$.

Consider an involutive isometry $c\: \L \to \L$ and denote by $L_{\pm c}$
its ($\pm 1$)-eigen\-lattices.
The involution $c$ is called \emph{geometric} if $h \in L_{-c}$ and $\sigma _+ L_{-c} = 2$.
As explained above, a marking of a nonsingular real quartic $X \subset \Cp3$ takes $c_X$ to a geometric involution
on $\L$.
This involution
is called
the \emph{homological type} of $X$; it is determined by $X$ up to
the action of $\OG^+_h(\L)$.
Conversely, according to Nikulin \cite{Nikulin:forms}, any geometric involution $c\: \L \to \L$
is the homological type of a marked nonsingular real quartic,
and the periods of such quartics constitute the whole space
\[
\Omega^\circ \cap \{\R\omega_+ \oplus \R\omega_- \; | \; \omega_\pm \in L_{\pm c} \otimes \R\}. \label{real_periods}
\]

\subsection{Configurations}\label{s.configurations}
Motivated by \autoref{line_classes}, we define
a {\it line} in a $4$-polarized lattice $S$ as a vector $a \in S$
such that $a^2 = -2$ and $a \cdot h = 1$.
The set of all lines in~$S$ is denoted by $\Fn(S)$.

\definition\label{pre-conf}
A {\it pre-configuration} is a $4$-polarized lattice $S$
generated over $\Q$ by its polarization $h$ and all lines $a \in S$.
A pre-configuration $S$ is called {\it hyperbolic} if $\sigma_+(S) = 1$.
A \emph{configuration} is a
nondegenerate
hyperbolic pre-configuration~$S$
that
contains no vector $e$ such that either
\roster
\item\label{ex-div} $e^2 = -2$ and $e \cdot h= 0$, or
\item\label{ell-pencil} $e^2 = 0$ and $e \cdot h = 2$
\endroster
(\cf. Proposition \ref{Saint-Donat}).
For a pre-configuration $(S, h)$ and a subset $A \subset \Fn(S)$,
the notation $\Span_h(A)$ stands for
the pre-configuration $S'\subset S$ generated (over~$\Z$) by $A$ and $h$.
\enddefinition

\remark\label{subconf}
Let $S$ be a nondegenerate hyperbolic pre-configuration. Then
\roster*
\item
$S$ contains finitely many lines, and
\item
any pre-configuration $S' \subset S$ is also
nondegenerate and hyperbolic.
\endroster
In particular, if $S$ is a configuration, then so is~$S'$.
\endremark

Let
$L \subset \L$ be a nondegenerate primitive polarized sublattice.
An {\it $L$-configuration}
is a
configuration $S\subset L$ \emph{primitive} in~$L$.
Two $L$-configurations $S',S''\subset L$
are said to be \emph{isomorphic}, or \emph{strictly isomorphic},
if
there exists an element of the group $\OG^+_h(\L,L)$ sending~$S'$ to~$S''$.
An {\it $L$-realization}
of a pre-configuration $S$ is a polarized isometry $\psi\: S \to L$
such that the image $\Im\psi$ is non-degenerate, \ie,
$\Ker\psi = \ker S$.
If the primitive hull $(\Im(\psi) \otimes \Q) \cap L$ is an
$L$-configuration,
the realization~$\psi$ is called \emph{geometric}.
A configuration admitting
a primitive geometric $L$-realization is called
\emph{$L$-geometric}
(or just \emph{geometric} if $L=\L$).

Note that there is a subtle difference between
$\L$-configurations and geometric ones: typically,
the former are considered up to
the action of $\OG_h^+(\L)$,
whereas the latter, up to abstract automorphisms of polarized lattices
(\cf. \autoref{lem.unique}).

To
simplify the classification of configurations,
we introduce also the notion of weak
isomorphism. Namely, two $\L$-configurations are said to be \emph{weakly
isomorphic} if they are taken to each other by an element of the group
$\OG_h(\L)$; in other words, we disregard the positive sign structure
on~$\L$.
Respectively, an $\L$-configuration $S\subset\L$
is called \emph{symmetric}
if it is preserved by an element $a \in \OG_h(\L) \sminus \OG^+_h(\L)$;
if such an element $a$ can be chosen involutive (respectively, involutive and identical on $S$),
the configuration $S$ is called \emph{reflexive} (respectively, \emph{totally reflexive}).
Putting $c = -a$, one concludes that $S$ is totally reflexive if and only if $S \subset L_{-c}$
for some geometric involution $c$. It is also clear that each weak
isomorphism class consists of one or two strict isomorphism classes,
depending on whether the configurations are symmetric or not, respectively.

\lemma\label{totally_reflexive}
An $\L$-configuration $S$ is totally reflexive if and only if the orthogonal complement $S^\perp$
contains either $[2]$ or $\bU(2)$.
\endlemma

\proof
We use the classification of geometric involutions found in \cite{Nikulin:forms}.
On the one hand, any sublattice isomorphic to $[2]$ or $\bU(2)$ in $h^\perp \subset \L$
is of the form $L_{+c}$ for some geometric involution $c$.
On the other hand, for any geometric involution $c$ the sublattice $L_{-c}$ is totally reflexive.
\endproof

\subsection{The arithmetical reduction}\label{arithmetical_reduction}
Let $X \subset \PP^3$ be a nonsingular quartic surface.
Choosing a marking
$\psi\: L_X \to \L$,
we obtain an $\L$-configuration $\psi(\Fano(X))$
(see Proposition \ref{Saint-Donat}).
Since any two markings differ by an element of $\OG^+_h(\L)$,
the surface $X$ gives rise to a well-defined isomorphism class
$[\Fano(X)]$ of $\L$-configurations.

Two nonsingular quartics $X_0$ and $X_1$ in $\Cp3$
are said to be \emph{equilinear deformation equivalent}
if there exists a path $X_t$, $t \in [0,1]$, in the space of nonsingular quartics
such that the number of lines in $X_t$ remains constant.

\theorem\label{th.complex}
The map $X \mapsto [\Fano(X)]$
establishes
a bijection between the set of equilinear deformation
classes of nonsingular quartics in $\PP^3$
and that of strict isomorphism classes of
$\L$-configurations.
\endtheorem

\proof
For the surjectivity, we choose a period $\omega \in \Omega^\circ$
so that $\omega^\perp \cap \L$
represents the chosen class of $\L$-configurations
and apply \autoref{Saint-Donat} and \autoref{line_classes}.
For the injectivity, we prove a stronger statement, \viz. the connectedness
of the space $\Omega'(S)$
of marked nonsingular quartics whose lines are taken by the marking
to the lines of a fixed $\L$-configuration $S \subset \L$.
To this end, consider the spaces
\[*
\Omega(S)=\{\Go\in\Omega\,|\,S\subset\Go^\perp\},\qquad
\Omega^\circ(S)=\Omega(S)\cap\Omega^\circ.
\]
By \autoref{moduli}, the latter is a fine moduli space of marked nonsingular
quartics $(X,\psi)$ such that $\psi(\Pic X)\supset S$; hence, by
\autoref{line_classes}, the space $\Omega'(S)$ is obtained from $\Omega^\circ(S)$ by
removing the union of the subspaces $\Omega_e$, where
\roster[3]
\item\label{ex-line1}
$e\in\L\sminus S$ is such that $e^2=-2$ and $e\cdot h=1$.
\endroster
In other words, $\Omega'(S)$ is obtained from a connected (in a sense, convex)
manifold $\Omega(S)$ by removing the codimension~$2$ subspaces $\Omega_e$
with $e$ as in \autoref{Saint-Donat}\iref{ex-div1}, \iref{ell-pencil1} or as
in~\iref{ex-line1} above. This family of subspaces is obviously locally
finite, and this fact implies the connectedness of the complement.
\endproof

\proposition\label{symmetric-reflexive}
Let $S$ be an $\L$-configuration, and
denote by $\Cal X$ the equilinear deformation class
corresponding to $S$ under the bijection of \autoref{th.complex}.
Then\rom:
\roster*
\item $\Cal X$ is invariant under the complex conjugation if and only if $S$ is symmetric\rom;
\item $\Cal X$ contains a real quartic if and only if $S$ is reflexive.
\endroster
\endproposition

\proof
Since $\omega_{\bar X}$ is $\omega_X$ with the orientation reversed,
the statement follows from the description of the moduli space $\Omega'(S)$
given in the proof of \autoref{th.complex}.
\endproof

A nonsingular quartic $X \subset \PP^3$ is called \emph{$\Fano$-maximal}
if $\rank\Fano(X) = 20$.

\addendum\label{app.complex}
The map $X \mapsto [\Fano(X)]$ establishes
a bijection between the set of projective equivalence classes of $\Fano$-maximal quartics in $\PP^3$
and that of isomorphism classes of
$\L$-configurations of rank $20$.
\endaddendum

\proof
Such quartics have maximal Picard rank, and for $S \subset \L$ of rank $20$,
the moduli space $\Omega'(S)/\!\PGL(4, \C)$ (\cf. the proof of \autoref{th.complex})
is discrete.
\endproof

Now, consider a nonsingular real quartic $X \subset \PP^3$ of
a certain homological type $c\: \L \to \L$.
The real structure on $X$ reverses the orientation of any real algebraic curve $C \subset X$,
thus reversing the class $[C] \in L_X$.
Hence, as above, considering real lines only,
we can define the \emph{real Fano configuration} $\Fano_\R(X)$ and
the isomorphism class
$[\Fano_\R(X)]$ of $L_{-c}$-configurations.

The following statements are straightforward, \cf. \eqref{real_periods}.

\theorem\label{th.real}
The real Fano configuration
of a nonsingular real quartic $X \subset \PP^3$
of homological type
$c\: \L \to \L$
is $L_{-c}$-geometric.
Conversely, any isomorphism class of
$L_{-c}$-configurations
is of the form $[\Fano_\R(X)]$ for some nonsingular real quartic $X \subset \PP^3$
of homological type~$c$.
\done
\endtheorem

\corollary\label{cor.real}
An $\L$-configuration~$S$ is in the class $[\Fano_\R(X)]$
for some nonsingular real quartic $X \subset \PP^3$
if and only if
$S$ is totally reflexive.
\done
\endcorollary

A nonsingular real quartic $X$ is called \emph{$\Fano_\R$-maximal}
if  $\rank\Fano_\R(X) = 20$. Even though we do not study equivariant
equilinear deformations
of
real quartics, in the case of the maximal Picard rank, where the moduli
spaces are discrete,
we still have projective equivalence; the precise statement is as follows.

\addendum\label{app.real}
The map $X \mapsto [\Fano_\R(X)]$ establishes
a bijection between the set of
real
projective equivalence classes of $\Fano_\R$-maximal real quartics in $\PP^3$
of a given homological type $c\: \L \to \L$
and that of isomorphism classes of
$L_{-c}$-configurations of rank~$20$.
\done
\endaddendum

\section{Geometry of configurations}\label{S.geometry}

In this section, we study the simplest properties of configurations, \viz.
those with a simple geometric interpretation. Most statements hold
without the assumption that the configuration should be geometric.

\subsection{Planes}\label{s.planes}
Fix a configuration~$S$ and denote by $h \in S$
its polarization.

\lemma\label{matrix}
For any two distinct lines $a_1,a_2\in S$ one has
$a_1 \cdot a_2 = 0$ or $1$.
\endlemma

\proof
Let $a_1 \cdot a_2 = x$, and consider the subconfiguration
$S':=\Span_h(a_1,a_2)$ (see \autoref{subconf}).
From $\det S' > 0$, one has $-1 \le x \le 2$.
If $x = -1$, then $a_1 - a_2$ is as in \autoref{pre-conf}\iref{ex-div};
if $x = 2$, then $a_1 + a_2$ is as in \autoref{pre-conf}\iref{ell-pencil}.
\endproof

Two distinct lines $a_1,a_2\in S$
are said to \emph{intersect} (respectively, to be \emph{disjoint}, or
\emph{skew}) if $a_1 \cdot a_2 = 1$
(respectively, $a_1 \cdot a_2 = 0$).
We regard the set of lines $\Fn(S)$ as a graph, with a pair of lines
(regarded as vertices)
connected by an edge if and only if the lines intersect.
A \emph{subgraph} of $\Fn(S)$ is always assumed induced.

A \emph{plane} in a configuration~$S$ is a collection
$\{a_1, a_2, a_3, a_4\} \subset S$ of four pairwise intersecting lines.

\lemma\label{plane}
For any plane $\{a_1, a_2, a_3, a_4\}\subset S$ one has
$a_1 + a_2 + a_3 + a_4 = h$.
\endlemma

\proof
The difference $h - (a_1 + a_2 + a_3 + a_4)$ is in the kernel
of $\Span_h(a_1,a_2,a_3,a_4)$;
hence, this difference is zero, see \autoref{subconf}.
\endproof

\corollary[of Lemmas \ref{matrix} and \ref{plane}]\label{cor-plane}
Let $\alpha = \{a_1, a_2, a_3, a_4\} \subset S$ be
a plane and $b \in S$ a line not contained in $\alpha$.
Then $b$ intersects exactly one line of $\alpha$.
\done
\endcorollary

The \emph{valency} $\val l$ of a line $l\in S$
is the number of lines in $S$ that intersect $l$.

\corollary[of \autoref{cor-plane}]\label{cor-Segre}
For any plane $\alpha = \{a_1, a_2, a_3, a_4\} \subset S$,
one has
\[*
\ls| \Fn(S) | = \val a_1 + \val a_2 + \val a_3 + \val a_4 - 8.
\]
\endcorollary

\lemma\label{two-planes}
Let $a_1$, $a_2 \in S$ be two intersecting lines, and
assume that there is a line $b_1 \in S$ that intersects both $a_1$ and $a_2$.
Then, there exists exactly one other line $b_2 \in S$ intersecting $a_1$ and $a_2$.
Furthermore, the lines $a_1,a_2,b_1,b_2$ form a plane.

As a consequence, if two planes $\alpha_1$, $\alpha_2 \subset S$
share two lines,
then $\alpha_1=\alpha_2$.
\endlemma

\proof
For the existence, let $b_2 = h - (a_1 + a_2 + b_1)$ (\cf. \autoref{plane}).
For the uniqueness, consider a line $c$ as in the statement.
If $b_1 \cdot c = 0$, then the difference $h - (a_1 + a_2 + b_1 + c)$
is as in \autoref{pre-conf}\iref{ex-div}.
Otherwise, one has $b_1 \cdot c = 1$ by \autoref{matrix},
and $\{a_1,a_2,b_1,c\}$ is a plane.
Hence, $c = b_2$ by \autoref{plane}.
\endproof

If two distinct lines lie in a (unique) plane $\Ga\subset S$, they are said to
\emph{span}~$\Ga$.

\subsection{Skew lines}\label{s.skew}
We keep the notation $(S,h)$ from the previous section.
The next lemma states some properties of skew lines.

\lemma\label{lem.skew}
Consider a number of lines $a_1,\ldots,a_m,b_1,\ldots,b_n\in S$ such that
all $a_i$ are pairwise disjoint, all $b_j$ are pairwise distinct, and
$a_i\cdot b_j=1$ for all $i=1,\ldots,m$, $j=1,\ldots,n$.
Then the following holds\rom:
\roster
\item\label{skew.plane}
if $m\ge2$, then all lines $b_j$ are pairwise disjoint\rom;
\item\label{skew.2}
if $m=2$, then $n\le10$\rom;
if $n=9$, then there exists a unique other line $b_{10}$ such that
$a_i\cdot b_{10}=1$ for $i=1,2$\rom;
\cf. also \autoref{cor.10.fibers} below\rom;
\item\label{skew.4}
if $m=4$, then $n\le4$\rom;
if $n=3$, then there exists a unique other line $b_{4}$ such that
$a_i\cdot b_{4}=1$ for $i=1,2,3$\rom; for this line, also
$a_4\cdot b_4=1$\rom;
\item\label{skew.4'}
if $m=n=4$, then any other line $c\in S$ intersects exactly two of the
given lines $a_1,\ldots,a_4$, $b_1,\ldots,b_4$\rom;
\item\label{skew.5}
if $m\ge3$, then $n\le4$\rom; if $m\ge5$, then $n\le2$.
\endroster
\endlemma

\proof
Item~\iref{skew.plane} is a partial restatement of \autoref{two-planes}.
The next two statements are proved similarly, with
\[*
b_{10}=4h-3(a_1+a_2)-(b_1+\ldots+b_9)
\]
in item~\iref{skew.2} and
\[*
b_4=2h-(a_1+\ldots+a_4+b_1+b_2+b_3)
\]
in item~\iref{skew.4}.
In the latter case, if $a_4\cdot b_4$ were~$0$, the vector
$a_1+\ldots+b_4-2h$ would be as in \autoref{pre-conf}\iref{ex-div}.
The expression for $b_4$ proves also item~\iref{skew.4'}, and
item~\iref{skew.5} is a simple consequence of item~\iref{skew.4}.
\endproof

Recall that our ultimate goal is the study of the configuration~$S$ of lines
in a nonsingular quartic surface~$X$.
From this perspective, as the name suggests, a plane is the subconfiguration
cut on~$X$ by a plane in~$\Cp3$, provided that the intersection splits
completely into components of degree one. A collection
$a_1,\ldots,a_4,b_1,\ldots,b_4$ as in \autoref{lem.skew}\iref{skew.4}
and~\iref{skew.4'} can similarly be interpreted as the intersection of~$X$
with a quadric (the lines~$a_i$ and~$b_j$ lying in the two distinct families of
generatrices),
and a subconfiguration as in
\autoref{lem.skew}\iref{skew.2} is (probably, a special case of) the
intersection of~$X$ with another quartic.
The following lemma, not used in the paper, is in the same spirit: it
describes the intersection of~$X$ with a cubic.
For the statement, define a \emph{double sextuple} as a collection of lines
$a_1,\ldots,a_6,b_1,\ldots,b_6$ in a configuration~$S$ intersecting as
follows:
\[
a_i\cdot b_j=1-\delta_{ij}
\label{eq.double.6}
\]
(where $\Gd_{ij}$ is the Kronecker symbol).

\lemma\label{lem.cubic}
Let $A':=\{a_1,\ldots,a_6,b_1,\ldots,b_5\}\subset S$ be a collection of lines
which
satisfy~\eqref{eq.double.6}.
Then there is a unique line $b_6\in S$ completing $A'$ to a double
sextuple $A$.
Furthermore,
all elements of~$A$ are pairwise distinct,
the lines~$a_i$ are pairwise disjoint,
the lines $b_j$ are pairwise disjoint, and
any other line $c\in S$
intersects exactly three elements of $A$.
\endlemma

\proof
The twelfth line is
\[*
b_6=3h-(a_1+\ldots+a_6+b_1+\ldots+b_5),
\]
and the other statements are immediate, \cf. the proof of \autoref{lem.skew}.
\endproof

\subsection{Pencils}\label{s.pencils}
Let $X \subset \Cp3$ be a nonsingular quartic such that
$\rank\Fano(X)\ge2$. Fix a line $l \subset X$. The pencil of planes through $l$ gives rise to
an elliptic pencil $X \to \Cp1$. Each fiber containing a line is
reducible:
it splits either into three lines
or a line and a conic;
in the former case, the three lines and~$l$ form a plane in~$\Fano(X)$.
Clearly, the lines in $X$ contained in the fibers of the pencil defined by $l$
are precisely those intersecting $l$.
Motivated by this observation,
we define a \emph{pencil} $\pencil$ in a
configuration $(S, h)$ as a set of lines
satisfying the following properties:
\roster*
\item
all lines in $\pencil$ intersect a given line $l$, called the
\emph{axis} of $\pencil$;
\item
if $a_1$, $a_2 \in \pencil$ and $a_1 \cdot a_2 = 1$, then
$h-l-a_1-a_2 \in \pencil$ (\cf. \autoref{plane}).
\endroster
\autoref{two-planes} implies that
\[*
\text{$a \sim b$ if $a = b$ or $a \cdot b = 1$}
\]
is an equivalence relation on $\pencil$. The equivalence classes are called
the \emph{fibers} of $\pencil$. The number $m$ of lines in a fiber may take
values $3$ or $1$; a fiber consisting of $m$ lines is called an
\emph{$m$-fiber}, and the number of such fibers is denoted by
$\fb_m(\pencil)$.
By \autoref{cor-plane}, $\pencil$ has a unique axis whenever
$\fb_3(\pencil)\ge1$ and $\fb_3(\pencil)+\fb_1(\pencil)\ge2$.

Each line $l\in S$ gives rise to a well-defined pencil
\[*
\pencil(l):=\{a\in\Fn S\,|\,a\cdot l=1\};
\]
such a pencil is called \emph{maximal}.
Any line $a \in S$ disjoint
from $l$ is called a \emph{section} of $\pencil(l)$ or any subpencil thereof.
The set of sections of~$\pencil$ depends on the ambient
(pre-)configuration~$S$; it
is denoted by $S(\pencil)$. By definition,
\[*
S(\pencil)=\{a\in\Fn(S)\,|\,a\cdot l=0\}.
\]
Clearly, for any line~$l\in S$,
one has
\[*
\val l=\ls|\pencil(l)|=3\fb_3(\pencil(l))+\fb_1(\pencil(l)).
\]
The number $\mult l:=\fb_3(\pencil(l))$
is called the \emph{multiplicity}
of~$l$. Alternatively,
$\mult l$ is the number of distinct planes containing~$l$.

Two pencils $\pencil_1$, $\pencil_2$ are called \emph{obverse} if their axes
are disjoint; otherwise, the pencils are called \emph{adjacent}.
The following lemma is an immediate consequence of Lemmas~\ref{two-planes}
and \ref{lem.skew}\iref{skew.2}.

\lemma\label{lem.2-pencils}
Let $\pencil_1\ne\pencil_2$ be two pencils. Then
\roster
\item
$\ls|\pencil_1\cap\pencil_2|\le10$ if $\pencil_1$, $\pencil_2$ are obverse,
and
\item
$\ls|\pencil_1\cap\pencil_2|\le2$ if $\pencil_1$, $\pencil_2$ are adjacent.
\done
\endroster
\endlemma

\subsection{Combinatorial invariants}\label{s.invariants}
A pencil~$\pencil$ is often said to be of \emph{type~$(p,q)$},
where $p:=\fb_3(\pencil)$ and $q:=\fb_1(\pencil)$. If an $\L$-realization~$\psi$
is fixed, the pencil is called \emph{primitive} or \emph{imprimitive} if so is
the sublattice $\Span_h\psi(\pencil)\subset\L$.
In this case, the type is further
refined to $(p,q)\prim$ and $(p,q)\impr$, respectively.
A geometric configuration containing a maximal pencil~$\pencil$
of type $(p,q)^*$ is
called a \emph{$(p,q)^*$-configuration}, and the pair $(S,\pencil)$ is called
a \emph{$(p,q)^*$-pair}.
The multiset
\[*
\Ps(S):=\bigl\{\text{type of $\pencil(l)$}\bigm|l\in\Fn(S)\bigr\}
\]
is called the \emph{pencil structure} of a configuration~$S$.
We usually represent $\Ps(S)$
in the partition notation (see, \eg, \autoref{obs.pencils} below):
a ``factor'' $(p,q)^a$ means that $S$ has $a$ pencils of type $(p,q)$.

The \emph{linking type} $\lk(\pencil_1,\pencil_2)$ of a pair of obverse
pencils is the pair $(\mu_1,\mu_3)$, where
$\mu_1:=\ls|\pencil_1\cap\pencil_2|$ and $\mu_3$ is the number of lines in
$\pencil_1\cap\pencil_2$ that belong to a $3$-fiber both in~$\pencil_1$
and~$\pencil_2$.
If $\pencil_i=\pencil(l_i)$, $i=1,2$, we also use the notation
$\lk(l_1,l_2)$.
The multiset
\[*
\Ls(S):=\bigl\{\lk(l_1,l_2)\bigm|
 l_1,l_2\in\Fn(S),\ l_1\cdot l_2=0\bigr\}
\]
is called the \emph{linking structure} of~$S$.

Clearly, both $\Ps(S)$ and $\Ls(S)$ are invariant under isomorphisms.

\section{The arithmetics of pencils}\label{S.arithm}

In this section, we study the more subtle properties of geometric
configurations related to their primitive embeddings to~$\L$.

\subsection{Notation and setup}\label{s.setup}
Throughout this section, we consider a pencil~$\pencil$ of a certain
type~$(p,q)$.
Thus, we have the sets $\fib_3\pencil=\{1,\ldots,p\}$ and
$\fib_1\pencil=\{1,\ldots,q\}$ of the $3$- and $1$-fibers
of~$\pencil$, respectively, and the full set
$\fib\pencil:=\fib_3\pencil\sqcup\fib_1\pencil$ of fibers
is their disjoint
union.
We regard~$\pencil$ as a pencil in the ``minimal'' configuration $P:=P_{p,q}$,
which is
generated over~$\Z$ by $\pencil$ itself, the axis~$l$, and the
polarization~$h$. We also keep in mind a geometric realization
$\psi\:P\to\L$, identifying $\pencil$ and~$P$ with their images in~$\L$ and
denoting by~$\tP$ the primitive hull
$(P\otimes\Q)\cap\L$.

When speaking about sections of~$\pencil$, we assume $\pencil$ embedded to a
configuration $S$, which is usually not specified.
(One can consider the minimal configuration generated by~$P$ and the
sections in question.) However, {\em we always assume that the realization
of~$P$ extends to a geometric realization $S\to\L$.}

The group of symmetries of $\pencil$ is obviously
\[*
\G_{p,q}:=(\SG3^p\rtimes\SG{p})\times\SG{q}.
\]
In addition to~$h$ and~$l$,
consider the following classes in $P_{p,q}$:
\roster*
\item
$m_{i,j}$, $i\in\fib_3\pencil$, $j\in\CG3$, the lines in the $3$-fibers;
\item
$n_k$, $k\in\fib_1\pencil$, the lines in the $1$-fibers.
\endroster
Then $P_{p,q}$ is the hyperbolic lattice
freely generated by $h$, $l$, $m_{i,j}$, $i\in\fib_3\pencil$,
$j=\pm1$, and $n_k$, $k\in\fib_1\pencil$.
For the lines $m_{i,\pm1}$, we will also use the shortcut $m_{i.\pm}$.

\observation\label{obs.3-torsion}
One has $\det P_{p,q}=-3^{p+2}(-2)^q$.
The $3$-primary part $\discr_3P_{p,q}$ contains the classes
represented by the following mutually orthogonal vectors:
\roster*
\item
$\Gl:=\frac13(l-h)$: one has $\Gl^2=0$ and $\Gl\cdot h=\Gl\cdot l=-1$;
\item
$\Gm_i=\Gm_{i,0}:=\frac13(m_{i,+}-m_{i,-})$, $i\in\fib_3\pencil$:
one has $\Gm_i^2=-\frac23$ and $\Gm_i\cdot h=0$.
\endroster
If $r:=p+q-1\ne0\bmod3$, then $\discr_3P_{p,q}$ is generated by $\Gm_i$,
$i\in\fib_3\pencil$, and
the order~$9$ class of the vector
\roster*
\item
$\Gu:=\frac13\bigl(l-r\Gl-\sum_{k=1}^qn_k\bigr)$;
note that $3\Gu=-r\Gl\ne0\bmod P$.
\endroster
Hence, in this case the subgroup of elements of order~$3$ is
generated by~$\Gl$ and~$\Gm_i$.
If $p+q=1\bmod3$, then $\discr_3P_{p,q}$ is generated by $\Gl$, $\Gm_i$,
and the order~$3$ class of
\roster*
\item
$\Go:=\frac13\bigl(l+\sum_{i=1}^p(m_{i,+}+m_{i,-})-\sum_{k=1}^qn_k\bigr)$.
\endroster
The $2$-primary part $\discr_2P_{p,q}$
is generated by the classes of $3\Gn_k$, where
\roster*
\item
$\Gn_k:=n_k^*=-\frac12(\Gl+n_k)$, $k\in\fib_1\pencil$:
one has $\Gn_k^2=-\frac12$ and $\Gn_k\cdot h=0$.
\endroster
The class $\Gm_i\in\discr P_{p,q}$ is also represented by the vector
$\bar\Gm_i^+:=\frac13(m_{i,+}+2m_{i,-})$, so that
one has $\bar\Gm_i^2=-\frac23$ and
$\bar\Gm_i\cdot h=1$.
The class $-\Gm_i\in\discr P_{p,q}$ is also represented by
$\bar\Gm_i^-:=\frac13(2m_{i,+}+m_{i,-})$.
For any line $a\in\pencil$, the class $\Gl$ is represented by the vector
$\Gl+a\in h^\perp$, so that one has $(\Gl+a)^2=-2$.
\endobservation

The following two statements are immediate.

\lemma\label{lem.sum.3}
For any triple of distinct indices $i,j,k\in\fib_3\pencil$ and any $u\in\CG3$, the
classes $\pm\Gl$ and $u\Gl\pm\Gm_i\pm\Gm_j\pm\Gm_k$ are represented by vectors of
square~$(-2)$ in $h^\perp\subset P_{p,q}$.
Hence, these classes cannot belong to the pivot~$\tP/P$.
\done
\endlemma

\lemma\label{lem.sum.4}
The sum of any four distinct elements of the form $3\Gn_k$, $k\in\fib_1\pencil$,
is represented by a vector of
square~$(-2)$ in $h^\perp\subset P_{p,q}$.
Hence, the class of such a sum cannot belong to the pivot~$\tP/P$.
\done
\endlemma

\subsection{Euler's bound}\label{s.large}
We start with eliminating very large pencils.

\proposition\label{prop.Euler}
The type $(p,q)$ of a pencil
contained in a geometric
configuration
satisfies the inequalities
\[*
3p+2q\le24\quad\text{and}\quad 3p+q\le20.
\]
\endproposition

\corollary[\cf. Rams, Sch\"{u}tt~\cite{rams:2012}]\label{cor.le20}
The valency of any line~$l$ in a geometric configuration~$S$
does not exceed $20$.
\done
\endcorollary

In the real case, there is an additional restriction to the types of pencils.

\proposition\label{prop.5-2.real}
A pencil~$\pencil$ contained in a totally reflexive geometric
configuration cannot be of type~$(6,0)\prim$
or~$(5,q)$, $q\ge2$.
\endproposition

\proof[Proof of Propositions~\ref{prop.Euler} and~\ref{prop.5-2.real}]
Assume that $(p,q)=(7,0)$.
By \autoref{obs.3-torsion}, the isotropic elements in $\discr_3P_{7,0}$ are:
\roster
\item\label{i.(7,0).3}
the classes mentioned in \autoref{lem.sum.3};
\item\label{i.(7,0).6}
classes of the form $u\Gl+\sum_{i\in I}\pm\Gm_i$, where $u\in\CG3$ and
$I\subset\fib_3\pencil$, $\ls|I|=6$;
all these classes form a single orbit of $\G_{7,0}$;
\item\label{i.(7,0).any}
classes of the form (up to sign) $\Go+u\Gl-\sum_{i\in I}\pm\Gm_i$, where
$I\subset\fib_3\pencil$ is any subset and $u=(5-\ls|I|)\bmod3$.
\endroster
Each class as in \autoref{i.(7,0).any} is represented by a vector of
square~$(-2)$ orthogonal to~$h$, \viz.
$\Go+(5-\ls|I|)\Gl-\sum_{i\in I}\bar\Gm_i^\pm$.
Hence, neither~\iref{i.(7,0).3} nor~\iref{i.(7,0).any} can belong to
the pivot $\tP/P$.
On the other hand, by \autoref{th.Nikulin}, one has $\ell_3(\tP/P)\ge2$
and $\tP/P$ must contain two distinct nontrivial orthogonal vectors
$\Gb_1$, $\Gb_2$ as
in~\iref{i.(7,0).6}.
On the other hand,
if both vectors are as in~\iref{i.(7,0).6}, then at least one of their linear
combinations is as in~\iref{i.(7,0).3}, \cf.~\cite{degt:Shapiro}.

Similar arguments apply to the other border cases: by \autoref{th.Nikulin},
one has
\roster*
\item
$\ell_3(\tP/P)\ge1$ if $(p,q)=(5,4)$ (use \autoref{lem.sum.3}),
\item
$\ell_2(\tP/P)\ge1$ if $(p,q)=(3,8)$,
\item
$\ell_2(\tP/P)\ge2$ if $(p,q)=(1,11)$ (use \autoref{lem.sum.4}), and
\item
$\ell_2(\tP/P)\ge3$ if $(p,q)=(0,13)$ (use \autoref{lem.sum.4}).
\endroster
In the case $(p,q)=(3,8)$, the only isotropic element allowed by
\autoref{lem.sum.4} is the characteristic element $\Gn:=\sum_{k=1}^8\Gn_k$.
The discriminant form $\Gn^\perp/\Gn$ is even, and the new lattice does not
embed to~$\L$ by \autoref{th.Nikulin}.

For \autoref{prop.5-2.real}, one
uses \autoref{obs.3-torsion} and \autoref{th.Nikulin}; the latter should
be applied to either $P\oplus[2]$ or an appropriate finite
index extension of $P\oplus[2]$ or $P\oplus\bU(2)$, see
\autoref{totally_reflexive}.
\endproof

The conclusion of \autoref{prop.Euler} can be recast as follows:
for any line~$l$ in a geometric configuration~$S$, one
has $\val l\le20$ and $\mult l\le6$; furthermore,
\[
\def\0{\phantom{0}}
\aligned
\text{if $\mult l$}&\le\00,\,\01,\,\02,\,\03,\,\04,\,\05,\,\06={\max},\\
\text{then $\val l$}&\le12,\,13,\,15,\,16,\,18,\,18,\,20={\max},\text{ respectively}.
\label{tablichka}
\endaligned
\]
It follows from~\eqref{tablichka} that $\max\{\val l\,|\,l\in\Fn(S)\}\le18$
if and only if $S$ does not contain a pencil of type $(6,q)$, $q\ge1$.

\remark\label{rem.Euler}
Interpreting pencil geometrically as in \autoref{s.pencils},
one can easily see that
the first inequality $3p+2q\le24$ in \autoref{prop.Euler} is nothing but the
well-known bound on the number and types of singular fibers in
an elliptic pencil.
\endremark

\subsection{Coordinates}
Consider a section~$s$ of a pencil~$\pencil$.
By \autoref{cor-plane}, for each index
$i\in\fib_3\pencil$, the section~$s$ intersects exactly one of the three lines
$m_{i,j}$, $j\in\CG3$; the corresponding index $\Ge_i:=j\in\CG3$
is called the \emph{$i$-th $3$-coordinate} of~$s$.
Introduce also the \emph{$k$-th $1$-coordinate} as the residue
$\Gr_k:=(s\cdot n_k)\bmod2\in\CG2$, $k\in\fib_1\pencil$.

We will treat the coordinate space
$\coord_{p,q}:=(\CG3)^p\times(\CG2)^q$ as an abelian group,
even though
only few linear combinations of coordinate vectors have
invariant meaning.
To avoid confusion with the operations in lattices, we will use $\op$ and
$\om$ for the addition and subtraction in~$\coord_{p,q}$, respectively.

\convention\label{conv.coord}
Given sections $s,s_1,s_2,\dots$ of~$\pencil$ and $u=1,3$,
we will use the following notation:
\roster*
\item
$\Ge_i:=\Ge_i(s)$ and $\Gr_k:=\Gr_k(s)$ are, respectively, the $3$- and
$1$-coordinates of~$s$;
\item
$[s]$ or $\bs:=[\Ge_1,\ldots,\Ge_p;\Gr_1,\ldots,\Gr_q]$ is the sequence of all
coordinates of~$s$;
\item
$\num u(s)$ is the number of non-vanishing $u$-coordinates of~$s$;
\item
$\dif u(s_1,s_2)$ is the number of positions where the $u$-coordinates
of~$s_1$, $s_2$ differ;
\item
$\scom3(s_1\cs s_2\cs\ldots):=\{i\in\fib_3\pencil\,|\,\Ge_i(s_1)=\Ge_i(s_2)=\ldots\}$;
\item
$\scom1(s_1\cs s_2\cs\ldots):=\{k\in\fib_1\pencil\,|\,\Gr_k(s_1)=\Gr_k(s_2)=\ldots=1\}$;
\item
$\scom(\ldots):=\scom3(\ldots)\sqcup\scom1(\ldots)$ (regarded as a set of
fibers of~$\pencil$);
\item
$\com*(\ldots)$ is the cardinality of the set $\scom*(\ldots)$ for $*=1,3$,
or empty;
\item
$\one:=\one_{p,q}=[0,\ldots,0;1,\ldots,1]\in\coord_{p,q}$.
\endroster
The same notation applies if all or some of $s$, $s_1$, $s_2$ are elements of
the coordinate space $\coord_{p,q}$.
The $3$-coordinates $\Ge_i(s)$, numbers $\num3(s)$,
and element $\one\in\coord$ depend on the indexing
of the lines in the $3$-fibers;
however, the sets $\scom3(\ldots)$, numbers $\dif3(s_1,s_2)$, and expressions
of the form
\[*
\one\op\bs,\quad
\bs_1\op\bs_2\op\bs_3=\one,\quad\text{or}\quad\bs_3=\one\om\bs_1\om\bs_2
\]
have invariant meaning.
Note also the difference between the definitions of $\scom3(\ldots)$ and
$\scom1(\ldots)$: in the former case, we count \emph{all} equal coordinates,
whereas
in the latter, only the \emph{non-vanishing} ones.
\endconvention

The following statements are immediate consequences of Lemmas~\ref{two-planes}
and~\ref{lem.skew}.

\lemma\label{lem.s.s=1}
Let $s_1$, $s_2$ be two sections of~$\pencil$ and $s_1\cdot s_2=1$.
Then
$\com(s_1\cs s_2)\le1$. If $\com(s_1\cs s_2)=1$, then there
is a
section~$s$ satisfying
$\bs\+\bs_1\+\bs_2=\one$\rom;
the
sections
$s$, $s_1$, $s_2$ and the only line $a\in\pencil$
intersecting all three of them constitute a plane.
\done
\endlemma

\lemma\label{lem.s.s=0}
Let $s_1$, $s_2$, $s_3$ be distinct sections of~$\pencil$.
Then\rom:
\roster
\item\label{com.5}
one has $\com(s_1\cs s_2)\le4$\rom;
\item\label{com.4}
if $\com(s_1\cs s_2)=4$, there is a unique
section~$s$ such that
$\bs\+\bs_1\+\bs_2=\one$\rom;
\item\label{com.3}
if $\com(s_1\cs s_2\cs s_3)=3$, the pencil $\pencil$ is not maximal.
\done
\endroster
\endlemma

\remark\label{rem.linear.combination}
In Lemmas~\ref{lem.s.s=1} and~\ref{lem.s.s=0}, as well as in the other
similar places below, the existence statement means that $s$ is a certain
(explicit, but not specified)
integral linear combination of the other sections involved and
generators of~$\tP$.
\endremark

\corollary\label{cor.coord}
If $p\ge5$, then, for any configuration $S\supset P$, the coordinate
map $c\colon S(\pencil)\to\coord_{p,q}$, $s\mapsto[s]$, is injective.
\done
\endcorollary

The injectivity of~$c$ for types $(4,\q)$ and $(3,7)$ is
discussed in \autoref{ss.4-q} below.

%

The next corollary deals with an obverse pencil in a configuration
$S\supset\pencil$.

\corollary\label{cor.obverse}
Given a section $s_0\in S(\pencil)$, consider
$s,s_1,s_2\in\pencil(s_0)\cap S_k(\pencil)$
and assume that $s_1\cdot s_2=1$. Then\rom:
\roster
\item\label{obverse.<=1}
one has $\com(s\cs s_0)\le1$\rom;
\item\label{obverse.common}
$\scom(s_0\cs s_1)=\scom(s_0\cs s_2)=\scom(s_1\cs s_2)=\scom(s_0\cs s_1\cs s_2)$\rom;
\item\label{obverse.maximal}
if $\pencil$ is maximal, then $\com(s_1\cs s_2)=1$\rom;
\item\label{obverse.3-fiber}
if $\pencil$ is maximal, then $s$ is in a $1$-fiber of~$\pencil(s_0)$ if and only if
$\com(s\cs s_0)=0$.
\endroster
\endcorollary

\proof
Statement~\iref{obverse.<=1} is a paraphrase of \autoref{lem.s.s=1}.
For~\iref{obverse.common} and~\iref{obverse.maximal}, just observe that
$s_0,s_1,s_2$ span a plane, and the forth line~$a$ of this plane must
intersect~$l$, see \autoref{cor-plane}; hence, either $a\in\pencil$ or
$\pencil$ is not maximal.
Finally, Statement~\iref{obverse.3-fiber} is a paraphrase
of~\iref{obverse.maximal}.
\endproof

Denote $D:=2p+\frac12q-2$ and, given a collection of sections $s_1,\ldots,s_k$,
let
\[*
r_{ij}:=(s_1\cdot s_2)+\tfrac19D+\tfrac12\com1(s_1\cs s_2)
 -\tfrac16(\num1(s_1)+\num1(s_2))-\tfrac13\dif3(s_1,s_2),\quad
 1\le i,j\le k,
\]
and define the \emph{determinant}
\[*
\det(s_1,\ldots,s_k):=\det[-r_{ij}]_{1\le i,j\le k}.
\]
The following lemma is a simple sufficient condition for the existence
of a collection of sections in terms of their coordinates
and pairwise intersections:
the orthogonal complement $P^\perp$ in any configuration $S\supset\pencil$
must be negative definite.

\lemma\label{lem.det}
For any collection $s_1,\ldots,s_k$ of sections one has
$\det(s_1,\ldots,s_k)\ge0$.
If $\det(s_1,\ldots,s_k)=0$, then the sections are linearly dependent.
\done
\endlemma

\subsection{Combinatorial rigidity}\label{s.rigidity}
The group $\G_{p,q}$ acts on the coordinate space $\coord_{p,q}$.
Furthermore, given two configurations $S,S'\supset\pencil$, any isometry
$(S,\pencil)\to(S',\pencil)$
induces an injection $\barS\into\barS'$, which is the restriction of an element
of $\G_{p,q}$. (Here, $\barS$ and $\barS'$ are the images of $S(\pencil)$ and
$S'(\pencil)$, respectively, under the coordinate map.)
A configuration $S\supset\pencil$
or, more precisely, pair $(S,\pencil)$
is called \emph{\rom(combinatorially\rom) rigid} if,
for any configuration $S'\supset\pencil$, any
bijection $g(\barS)=\barS'$ restricted from an element $g\in\G_{p,q}$
is induced by an isometry
$(S,\pencil)\to(S',\pencil)$.

We say that $S$ or $(S,\pencil)$ is
\emph{generated by a subset $\bar A\subset\barS$} if
$S=\bigl(\tP+\sum_{\bs\in\bar A}\Z s\bigr)/\!\ker$;
if $\bar A=\barS$, then $S$ is said to be \emph{generated by
sections}.
For such a configuration,
an obvious
sufficient condition for the combinatorial rigidity is that
the intersection $s_1\cdot s_2$ of a pair of sections
$s_1,s_2$ such that $\bs_1,\bs_2\in\bar A$
is determined by their images $\bs_1,\bs_2$, \ie, for any other
configuration $S'\supset\pencil$ and pair of sections
$s_1',s_2'\in S'(\pencil)$ such that $\barS'=\barS$ and $\bs_1'=\bs_1$,
$\bs_2'=\bs_2$, one has $s_1'\cdot s_2'=s_1\cdot s_2$.
By \autoref{lem.s.s=1}, an ambiguity may arise only if $\com(s_1\cs s_2)\le1$.
The following statement is a partial converse of \autoref{lem.s.s=1};
we do not need to assume that the configuration
$S\supset\pencil$ is geometric.

\lemma\label{lem.com2=1}
Let $p=6$, $(p,q)=(5,3)$, $p=4$ and $q\ge4$, or $(p,q)=(3,7)$.
Consider a
pair of sections $s_1,s_2\in S(\pencil)$
such that $\com(s_1\cs s_2)=1$.
Then, $\pencil$ has a pair of sections $s_1',s_2'$ such that
$s_1'\cdot s_2'=1$ and $[s_i']=\bs_i$, $i=1,2$, if and only if there is a
section~$s$ such that $\bs\+\bs_1\+\bs_2=\one$.
\endlemma

\proof
The necessity is given by \autoref{lem.s.s=1}. For the converse, it suffices
to show that three sections
$s,s_1,s_2$ as in the statement cannot be pairwise disjoint.
Most such triples are eliminated by
\autoref{lem.det}, and the few remaining ones violate condition~\iref{ex-div}
in \autoref{pre-conf}.
\endproof

\subsection{Primitivity and rigidity for type $(6,\q)$}\label{ss.6-q}
Primitive and imprimitive pencils of type $(6,\q)$ exhibit very different
behaviour. Here, we start with a few common observations; imprimitive pencils
are treated separately in the next section.

\proposition\label{prop.6-0}
Assume
that $p = 6$. Then the following holds\rom:
\roster
\item\label{6.0_one}
if $\pencil$ is not maximal or $q \ge 1$, then $\pencil$ is imprimitive\rom;
\item\label{6.0.pivot}
if $\pencil$ is imprimitive, then $\tP/P=\<\Gb\>$,
$\Gb:=\sum_{i=1}^6\Gm_i$, up to automorphism.
\endroster
\endproposition

\proof
The imprimitivity follows from \autoref{th.Nikulin}, and the only possible
nontrivial pivot is given by \autoref{obs.3-torsion} and \autoref{lem.sum.3}.
\endproof

\lemma\label{lem.6.0-sections}
Let $(p,q)=(6,0)$. Consider a geometric configuration $S\supset\pencil$,
let $\barS$ be the image of $S(\pencil)$ under the coordinate map,
and,
for a pair
$s_1,s_2\in S(\pencil)$, denote
$\bs:=\one\om\bs_1\om\bs_2\in\coord_{6,0}$.
Then the following holds\rom:
\roster
\item\label{i.(6,p).0}
if $\com(s_1\cs s_2)=0$ and
$s_1\cdot s_2=0$, then $\pencil$ is imprimitive
and $\frac13(s_1-s_2)\in\tP$\rom;
\item\label{i.(6,p).0-intr}
if $\com(s_1\cs s_2)=0$ or~$3$ and $\bs\in\barS$,
then $\pencil$ is imprimitive\rom;
\item\label{i.(6,p).1-intr}
if $\com(s_1\cs s_2)=1$, then $\bs\in\barS$ if and only if
$s_1\cdot s_2=1$\rom;
\item\label{i.(6,p).4}
if $\com(s_1\cs s_2)=4$, then $\bs\in\barS$.
\endroster
\endlemma

\proof
Statement~\iref{i.(6,p).0}:
the two vectors are linearly dependent by
\autoref{lem.det}; then $s_1-s_2=\Gb$ up to automorphism.

Statement~\iref{i.(6,p).0-intr}, $\com(s_1\cs s_2)=0$:
if $\pencil$ is primitive, then $s\cdot s_1=s\cdot s_2=s_1\cdot s_2=1$ by
Statement~\iref{i.(6,p).0}; hence, the three sections span a plane, and the
forth line of this plane is in $\pencil(l)\sminus\pencil$, which contradicts
\autoref{prop.6-0}\iref{6.0_one}.

Statement~\iref{i.(6,p).0-intr}, $\com(s_1\cs s_2)=3$: the imprimitivity of
$\Span_h(\pencil,s_1,s_2,s)$ is given by \autoref{th.Nikulin},
and the enumeration of isotropic elements
not realized by vectors $e$ as in \autoref{pre-conf}\iref{ex-div}
shows that the pivot is generated by~$\Gb$ (up to isomorphism).

Statements~\iref{i.(6,p).1-intr} and~\iref{i.(6,p).4} follow from
Lemmas~\ref{lem.com2=1}
and~\ref{lem.skew}\iref{skew.4}, respectively.
\endproof

\corollary\label{cor.6-0.combinatorial}
Any $(6,0)\prim$-configuration
generated by sections is rigid.
\done
\endcorollary

\subsection{Triplets of sections}\label{s.triplets}
In this section, we study in more detail an imprimitive pencil of
type $(6,0)\impr$. Thus, we fix a pencil~$\pencil$ and
number the lines $m_*$ in the fibers so that the pivot $\tP/P$ is generated
by the element~$\Gb$ introduced in \autoref{prop.6-0}.
Then, for any section~$s$,
\[
\Ge_1(s)+\ldots+\Ge_6(s)=0\bmod3.
\label{eq.6-0.coord}
\]
The group $\OG_h(\tP,l)$ is obviously the subgroup
\[
\tilde\G:=((\CG3)^5\rtimes\CG2)\rtimes\SG{6}\subset\G_{6,0};
\label{eq.6-0.group}
\]
indeed, the choice of~$\Gb$ gives rise to a distinguished cyclic order in
each fiber, which is well defined up to simultaneous reversal.
This group has a distinguished subgroup of order~$3$: it is generated by the
permutations
$\Gs^{\pm1}\:m_{i,j}\mapsto m_{i,j\pm1}$, $i\in\fib_3\pencil$, $j\in\CG3$.
A choice of one of these two generators makes $\coord_{6,0}$ an
$\F_3$-affine space.

Consider a configuration $S\supset\tP$ and let $\barS\subset\coord_{6,0}$
be the image of
$S(\pencil)$ under the coordinate map.

\lemma\label{lem.triplet}
The set $\barS$ is $\Gs$-invariant, \ie, $\bs_\pm:=\Gs^{\pm1}\bs\in\barS$
whenever $\bs\in\barS$. The three sections $\bs,\bs_\pm$ are
pairwise disjoint.
\endlemma

\proof
Up to automorphism, one can assume that $\bs=[1,\ldots,1]$.
Then the two other sections are $s+\Gb$ and
$s-2h+2l+\sum_{i=1}^6(m_{i,1}+m_{i,-1})-\Gb$.
\endproof

A subset $\{s,s_\pm\}\subset S(\pencil)$ or
$\{\bs,\bs_\pm\}\subset\barS$ as in \autoref{lem.triplet} is called a
\emph{triplet}. Two sections $s_1,s_2\in S(\pencil)$ are said to be
\emph{equivalent}, $s_1\sim s_2$, if they belong to one triplet.
Note that $\com(s_1\cs s_2)=0$ whenever $s_1\sim s_2$.

\lemma\label{lem.6-0.intr}
For a pair of sections $s_1,s_2\in S(\pencil)$, one has
$s_1\cdot s_2=1$ if and only
if $\com(s_1\cs s_2)\le1$ and $s_1\not\sim s_2$.
\endlemma

\proof
If $\com(s_1\cs s_2)=0$ and $s_1\cdot s_2=0$,
\autoref{lem.6.0-sections}\iref{i.(6,p).0} and the fact that
$\ell_3(\tP/P)=1$ imply that $s_1\sim s_2$.
If $\com(s_1\cs s_2)=1$, then, using~\eqref{eq.6-0.coord} and
Lemmas~\ref{lem.6.0-sections}\iref{i.(6,p).4} and~\ref{lem.triplet}, one can
easily show that $\one\om\bs_1\om\bs_2\in\barS$; then, $s_1\cdot s_2=1$ by
\autoref{lem.6.0-sections}\iref{i.(6,p).1-intr}.
\endproof

\corollary\label{cor.6-q.combinatorial}
Any $(6,\q)\impr$-configuration
generated by sections is rigid.
\done
\endcorollary

Note that, for $(6,\q)\impr$-configurations, the rigidity
holds in a very strong sense: the intersection of two
sections is completely determined by their coordinates.

The set of triplets can be coordinatized by the affine space
\[*
\CA:=\{\bs\in\coord_{6,0}\,|\,\text{$\bs$ satisfies~\eqref{eq.6-0.coord}}\}/\Gs.
\]
In fact, $\CA$ is naturally a principal homogeneous space over
the subquotient $\Gl^\perp\!/\Gl$ of the discriminant $\discr\tP$.
Denote by~$\qq$ the descent of the discriminant form of~$\tP$
reduced modulo~$\Z$; then, clearly,
\smash{$\qq(\bs_1-\bs_2)=\frac13(\dif(s_1,s_2)\bmod3)$}.
Comparing the orders, one can see that the group $\tilde\G/\Gs$ is isomorphic
to the full group $\OG(\qq)\rtimes\CA$ of $\qq$-isometries of~$\CA$.
In other words, any \smash{$\tilde\G$}-invariant property of a set of
sections $\barS\subset\coord_{6,0}$
satisfying~\eqref{eq.6-0.coord} and \autoref{lem.triplet} can be stated as a
``metric'' (with respect to~$\qq$) property
of the projection $\bar\CS$ of this set to~$\CA$.

Below, we state two properties that hold for any configuration~$S$,
not necessarily geometric.
Recall that the lines in $\Gl^\perp\!/\Gl$ can be subdivided into
\roster*
\item
$15$ \emph{positive} lines~$\ln^+$ and $15$ \emph{negative} lines~$\ln^-$,
with $q|_{\ln^\pm}\cong\<\pm\frac13\>$, and
\item
$10$ \emph{isotropic} lines $\ln^0$, with $q|_{\ln^0}\equiv0$.
\endroster
The planes in $\Gl^\perp\!/\Gl$ can be subdivided into
\roster*
\item
$20$ \emph{positive} planes~$\plane^+$ and $20$ \emph{negative} planes~$\plane^-$,
with $\plane^\pm\cong \ln^\pm\oplus \ln^0$,
\item
$45$ \emph{hyperbolic} planes, isomorphic to $\ln^+\oplus \ln^-$, and
\item
$45$ \emph{definite} planes, isomorphic to
$\ln^+\oplus \ln^+\cong \ln^-\oplus \ln^-$.
\endroster
(There are no isotropic planes.)
The same terminology applies to the lines/planes in~$\CA$, according to the
underlying vector space.
The group $\OG(\qq)$ acts transitively on
the set of lines/planes of the same type.

\lemma\label{lem.convex}
For any configuration $S\supset\tP$, the set $\bar\CS\subset\CA$ is ``convex''\rom:
whenever a negative line $\ln^-\subset\CA$ has two common points
with~$\bar\CS$, it is contained in $\bar\CS$.
\endlemma

\lemma\label{lem.plane}
Let $S\supset\tP$ be a configuration and $\plane^-\subset\CA$ a negative plane.
Then
the intersection $\bar\CS\cap\plane^-$ is contained in a line\rom;
equivalently, $\plane^-\not\subset\bar\CS$.
\endlemma

\proof[Proof of Lemmas~\ref{lem.convex} and~\ref{lem.plane}]
\autoref{lem.convex} is a restatement of
\autoref{lem.6.0-sections}\iref{i.(6,p).4}.
By \autoref{lem.convex}, the two restrictions in \autoref{lem.plane} are
equivalent: $\bar\CS\supset \plane^-$ if and only if $S$ contains three
non-collinear points of~$\plane^-$. If this is the case, the points can be chosen to
form an equilateral triangle with side~$-\frac13$; by
\autoref{lem.triplet}, we can find three sections $s_1,s_2,s_3$ so that
$\com(s_i\cs s_j)=1$ for all $i\ne j$ but $\com(s_1\cs s_2\cs s_3)=0$.
Then $s_i\cdot s_j=1$, see
\autoref{lem.6-0.intr}, and the three sections span a plane. This plane must
contain three more lines, \viz. the elements of $\pencil$ intersecting
the three pairs
$s_i$, $s_j$, $1\le i<j\le3$. This is a contradiction to \autoref{two-planes}.
\endproof

Remarkably, Lemmas~\ref{lem.convex} and~\ref{lem.plane} almost characterize
the sets of sections in configurations (not necessarily geometric) containing
a pencil of type $(6,0)\impr$: this fact is established experimentally during
the proof of \autoref{th.6-q}. There is but one extra restriction, stated
below without proof.

\lemma\label{lem.two.lines}
Let $S\supset\tP$ be a configuration
and $\ln',\ln''$ two parallel isotropic lines in a positive plane in $\CA$.
If $\ln'\subset\bar\CS$ and
$\bar\CS$ contains two points of~$\ln''$,
then $\ln''\subset\bar\CS$.
\done
\endlemma

\subsection{Primitivity and rigidity for types $(4,\q)$ and $(3,\q)$}\label{ss.4-q}
As above,
we fix a configuration $S\supset\pencil$ and denote by
$\barS\subset\coord_{p,q}$
the image of the set of sections $S(\pencil)$ under the coordinate map.
It follows from
\autoref{obs.3-torsion} and \autoref{lem.sum.3} that
any pencil of type $(4,q)$, $q\le5$, or $(3,q)$, $q\le6$, is primitive.
Below, we consider in detail the two extremal cases.

\proposition\label{prop.4-6}
If $(p,q) = (4, 6)$, then $\ppencil$ is imprimitive
and has a unique, up to isomorphism, geometric finite index extension.
Furthermore,
\roster
\item\label{4.6-l*}
$\ppencil$ has a unique section $l^* \in S$ intersecting all ten fibers\rom;
\item\label{4.6-max}
as a consequence, $\pencil$ is maximal in any
configuration\rom;
\item\label{4.6-plane}
if a section~$s$ intersects~$l^*$, then the lines $s$ and~$l^*$ span a
plane\rom;
\item\label{4.6-minus}
the set $\barS$
is invariant under the
involution $\bs\mapsto\bs\*:=\one\om l^*\om s$.
\endroster
If $(p,q)=(3,7)$, then $\pencil$ is imprimitive if and only if there is a
section~$l^*$ as in~\iref{4.6-l*} above\rom;
if this is the case, Statements~\iref{4.6-plane} and~\iref{4.6-minus} also
hold.
\endproposition

\proof
Let $(p,q)=(4,6)$.
The pivot $\tP/P$ must have $3$-torsion by
\autoref{th.Nikulin}, whereas its $2$-torsion is trivial by
\autoref{lem.sum.4}. In addition to
the classes mentioned in \autoref{lem.sum.3},
the isotropic elements in $\discr_3P_{4,6}$ are those
constituting the $\G_{4,6}$ orbits of the classes of $\pm\Go$,
see \autoref{obs.3-torsion}.
Hence, up to automorphism, $\tP/P$ is generated by~$\Go$,
and it is immediate that $\Go$ is a section~$l^*$ as in
\iref{4.6-l*}. A section with these properties is unique due to
\autoref{lem.s.s=0}\iref{com.5}.

If $(p,q)=(3,7)$, the only nontrivial elements that may be contained in the
pivot are the orbits of the classes of $\pm(\Go-\Gl)$, and $\Go-\Gl$ is a
section $l^*$ as in \iref{4.6-l*}.

With the above choice of $l^*$, we have $[l^*]=[0,\ldots,0;1,\ldots,1]$ and
\[
 \com(s\cs l^*)=4-3(s\cdot l^*)
\label{eq.l*}
\]
for any other section~$s$. (In particular, this relation restricts the
coordinate vectors realized by sections.) Clearly, $s\cdot l^*=1$
if and only if $\com(s\cs l^*)=1$, in which case $s$ and~$l^*$
intersect a third common line $a\in\pencil$ and thus span a plane;
in fact, this plane is $\{l^*,a,s,s\*\}$.
Statement~\iref{4.6-minus} follows from
\autoref{lem.s.s=1} or
\autoref{lem.s.s=0}\iref{com.4}
if $s\cdot l^*=1$ or $0$, respectively.
\endproof


\proposition\label{prop.3-7}
Let $(p,q)=(3,7)$.
If $\pencil$ is not maximal, then there is a section~$s$ of~$\pencil$ such
that $\num2(s)\le6$. Conversely, if there is a section~$s$ such that
$\num2(s)=6$, then $\pencil$ is not maximal.
\endproposition

\proof
The only pencil~$\pencil'$ that may properly contain~$\pencil$ is one of
type~$(4,6)$, and the section~$s$ as in the statement is the restriction
of~$l^*$ given by \autoref{prop.4-6}.
If $\pencil$ has a section~$s$ such that $\num2(s)=6$, then $s$ and~$l$
intersect nine disjoint lines; by \autoref{lem.skew}\iref{skew.2},
they must intersect a
tenth line.
\endproof

\proposition\label{prop.4-5}
Let $(p,q)=(4,5)$. Then
$\pencil$ is primitive, and
$\pencil$ is maximal in a geometric configuration~$S$
if and only
$\num1(s)\le4$ for each section $s\in S(\pencil)$.
\endproposition

\proof
The primitivity is essentially given by \autoref{obs.3-torsion} and
Lemmas~\ref{lem.sum.3} and~\ref{lem.sum.4}. By
\autoref{lem.skew}\iref{skew.4}, if there is a section~$s$
with $\num1(s)=5$,
the pencil has a tenth fiber.
Conversely, the only pencil that can
properly contain~$\pencil$ is one of type $(4,6)$,
and its section~$l^*$ given by \autoref{prop.4-6}\iref{4.6-l*}
restricts to~$\pencil$.
\endproof

\proposition\label{prop.coord.p<=4}
Let $p=4$, $q\ge4$ or $(p,q)=(3,7)$, and assume that $\pencil$ is maximal.
Then, for any ambient geometric configuration $S\supset\pencil$,
the coordinate map $c\colon S(\pencil)\to\coord_{p,q}$, $s\mapsto[s]$,
identifies at most one pair of sections.
Furthermore, if such a pair $s_1,s_2$ identified by~$c$ does exist, then
there also is a
\rom(unique\rom) section $l^*\in S(\pencil)$
such that $\bs_1+\bs_2+\bar l^*=\one$,
and, for this section~$l^*$, one has $\num2(l^*)=q$.
\endproposition

\proof
Let $s_1\ne s_2$ be a pair of sections such that $\bs_1=\bs_2$.
By \autoref{lem.skew}\iref{skew.4}, we have $\num2(\bs_i)+p\le4$ and,
if $\num2(\bs_i)+p=4$, there also is a section $l^*$ as in the statement.
The number of sections $l^*$ with $\num2(l^*)=q\ge4$ is
\roster*
\item
one if $(p,q)=(4,6)$ or at most one if $(p,q)=(3,7)$, see \autoref{prop.4-6},
\item
zero if $(p,q)=(4,5)$, see \autoref{prop.4-5}, and
\item
zero, one, or three if $(p,q)=(4,4)$, see \autoref{lem.skew}\iref{skew.4}.
\endroster
Furthermore, a given section~$l^*$ cannot share all $3$-coordinates with any
section other than $s_1,s_2$, see \autoref{lem.skew}\iref{skew.4} again.

If $(p,q)=(4,4)$ and $\pencil$ has three sections $l_1^*,l_2^*,l_3^*$ with
$\num2(l_i^*)=4$, one can easily show that only one pull-back
$c\1(\bar l_i^*+\one)$ may be nonempty, as otherwise $S$ does not admit a
geometric $\L$-realization.

In the remaining case $(p,q)=(3,7)$ and $\num2(\bs_i)=0$, one can use
\autoref{th.Nikulin} to show that the image of any geometric realization
of~$S$ must contain a section~$s$ of~$\pencil$ such that $\num2(s)=6$;
hence, $\pencil$ is not maximal, see \autoref{prop.3-7}.
\endproof

Till the rest of this section, we assume that $(p,q)=(4,6)$.

Denote $S^*(\pencil):=\{s\in S(\pencil)\,|\,s\cdot l^*=1\}$.
According to~\eqref{eq.l*}, the image of this set in $\coord_{4,6}$ can be
characterized as
\[
\barS^*=\bigl\{\bs\in\barS\bigm|\com(\bs\cs l^*)=1\bigr\}.
\label{eq.barS*}
\]
Let also
\[*
\barS^\circ:=\bigl\{s\in\barS\bigm|
 \text{$\com(s\cs s')=0$ and $\num1(s)+\num1(s')=1$ for some $s'\in\barS^*$}\bigr\}.
\]
The following statement complements \autoref{lem.com2=1};
we do not need to assume that
the configuration $S\supset\tP$ is geometric.

\lemma\label{lem.4-6.com2=0}
Let $(p,q)=(4,6)$. Consider a
pair of sections $s_1,s_2\in S(\pencil)$
such that $\com(s_1\cs s_2)=0$ and let
$\bs_1':=\one\om\bs_1\*\om\bs_2$ and
$\bs_2':=\one\om\bs_1\om\bs_2\*=(\bs_1')\*$. Then\rom:
\roster
\item\label{4-6.bounds}
one has $1\le\num1(s_1)+\num1(s_2)\le5$\rom;
\item\label{4-6.=5}
if $\num1(s_1)+\num1(s_2)=5$, then also $\bs_1',\bs_2'\in\barS$.
\endroster
If the pair $s_1,s_2$ is ``homogeneous'', then\rom:
\roster[\lastitem]
\item\label{4-6.S*}
if $\bs_1,\bs_2\in\barS^*$, one has $s_1\cdot s_2=0$, and
\item\label{4-6.not.S*}
if $\bs_1,\bs_2\notin\barS^*$, one has $s_1\cdot s_2=1$.
\endroster
If the pair is ``mixed'', $\bs_1\in\barS^*$ and $\bs_2\notin\barS^*$, then\rom:
\roster[\lastitem]
\item\label{4-6.S*-3}
if $\num1(s_1)+\num1(s_2)\ge3$, one has $s_1\cdot s_2=1$, and
\item\label{4-6.S*-2}
if $\num1(s_1)+\num1(s_2)=2$,
one has $s_1\cdot s_2=0$ if and only if $\bs_1',\bs_2'\in\barS$.
\endroster
\endlemma

If $S$ is required to be geometric, then one can also state that
$\num1(s_1)+\num1(s_2)\le4$ whenever $\bs_1\in\barS^*$.
We do not use this restriction explicitly.

\proof[Proof of \autoref{lem.4-6.com2=0}]
%
Statement~\iref{4-6.S*} is obvious, as $s_1,s_2$ are in distinct fibers
of the pencil $\pencil(l^*)$.
In all other cases, by \autoref{lem.skew}\iref{skew.4'}, the
section~$s_1$ must intersect exactly one (if $\bs_1\in\barS^*$) or two (if
$\bs_1\notin\barS^*$) of the lines $s_2$, $s_2\*$; with \eqref{eq.l*} taken
into account,
the intersection $s_1\cdot s_2\*$ is given by
Lemmas~\ref{lem.s.s=1} and~\ref{lem.com2=1}.
\endproof

\corollary[of Lemmas~\ref{lem.com2=1} and~\ref{lem.4-6.com2=0}]\label{lem.4-6.combinatorial}
Any $(4,6)$-configuration~$S$ 
generated by $\barS\sminus\barS^\circ$ is
rigid.
\done
\endcorollary

\remark\label{rem.4-6.combinatorial}
For many configurations,
the hypotheses of \autoref{lem.4-6.combinatorial} can also be verified
combinatorially,  using
Lemmas \ref{lem.s.s=0}\iref{com.4} and~\ref{lem.com2=1}:
assuming that $S\supset\tP$ is generated by sections, it is generated by
$\barS\sminus\barS^\circ$ if, for any $\bs\in\barS^\circ$, there is a pair
$\bs_1,\bs_2\subset\barS\sminus\barS^\circ$ such that
$\bs\op\bs_1\op\bs_2=\one$ and $\com(\bs_1\cs\bs_2)=1$ or~$4$.
\endremark


\subsection{Rigidity for type $(5,3)$}
As an immediate consequence of
\autoref{obs.3-torsion} and \autoref{lem.sum.3},
any pencil of type $(5,\q)$ is primitive.

In the next two statements, $S$ does not need to be geometric.

\lemma\label{lem.5-q.maximal}
Let $p=5$, $q\ge1$, and assume that $\pencil$ has a section. Then $\pencil$
is contained in a pencil $\pencil'$ of type $(6,\q)\impr$ if and only if
$\pencil$ has a pair of sections $s_1,s_2$ such that
$s_1\cdot s_2=0$, $\com3(s_1\cs s_2)=0$, and $\dif1(s_1,s_2)>0$.
\endlemma

\proof
If $\pencil\subset\pencil'$, then $s_1,s_2$ are two
appropriate equivalent sections
of~$\pencil'$, see \autoref{lem.triplet}. For the sufficiency, assume that
$(p,q)=(5,1)$ and
\[*
\bs_1=[0,0,0,0,0;1],\quad \bs_2=[1,1,1,1,1;0].
\]
Then an extra member of~$\pencil'$ is
$h-l+\sum_{i=1}^5(m_{i,+}-m_{i,0})-2n_1-3s_1+3s_2$.
\endproof

\corollary\label{cor.5-3.com=0}
Let $p=5$, $q\ge1$, and assume that $\pencil$ is maximal.
Then, for any pair $s_1,s_2\in S(\pencil)$
such that $\com(s_1\cs s_2)=0$ and $\num1(s_1)+\num1(s_2)>0$, one has
$s_1\cdot s_2=1$.
\done
\endcorollary

Let $(p,q)=(5,3)$ and assume that $\pencil$ is maximal
(see \autoref{lem.5-q.maximal} for a criterion).
Then, according to \autoref{lem.com2=1} and
\autoref{cor.5-3.com=0}, the intersection $s_1\cdot s_2$ may not be
determined by the coordinates $\bs_1,\bs_2\in\barS$ only if
\roster*
\item
one has $\com(\bs_1\cs\bs_2)=\num1(\bs_1)=\num1(\bs_2)=0$ and
\item
for any $\bs\in\barS$, if $\com(\bs\cs\bs_1)=\com(\bs\cs\bs_2)=0$, then
$\num1(\bs)=0$.
\endroster
(For the latter condition, if $\num1(\bs)>0$, then $s\cdot s_1=s\cdot s_2=1$
by \autoref{cor.5-3.com=0} and,
hence, $s_1\cdot s_2=0$,
see \autoref{cor.obverse}\iref{obverse.maximal}.)
Denote by $\barS^\circ\subset\barS$ the union of all such pairs
$(\bs_1,\bs_2)$.

\corollary\label{cor.5-3.combinatorial}
Any $(5,3)$-configuration~$S$
generated by $\barS\sminus\barS^\circ$ is rigid.
\done
\endcorollary

For another sufficient rigidity condition, consider a section
$s_0\in S(\pencil)$
and let $\sect(s_0):=\pencil(s_0)\cap S_k(\pencil)$.
If $\num1(s_0)>0$, this set is determined by the coordinates:
by \autoref{lem.com2=1} and \autoref{cor.5-3.com=0}, one has $s\in\sect(s_0)$
if and only if $\com(\bs\cs\bs_0)=0$ or $\com(\bs\cs\bs_0)=1$ and
$\one\om\bs\om\bs_0\in\barS$.
Furthermore, the intersections $s_1\cdot s_2$, $s_1,s_2\in\sect$,
are also known: they are given by
\autoref{cor.obverse}.

\corollary\label{cor.5-3.pencil}
Any $(5,3)$-pair $(S,\pencil)$ generated by the union
$\{\bs_0\}\cup\bar\sect(s_0)$ for some section
$s_0\in S(\pencil)$ such that $\num1(s_0)>0$
is rigid.
\done
\endcorollary

\subsection{Other types}\label{s.other.pencils}
For completeness, we discuss the primitivity of the other types of
pencils.
We treat the $3$- and $2$-torsion of the pivot separately.

\proposition\label{prop.3.pivot}
Let $\pencil$ be a pencil of type $(p,q)$ with $p\le2$.
If the
pivot $\tP/P$ has $3$-torsion, then
\roster*
\item
$p+q=10$, \ie,
$(p,q)=(2,8)$, $(1,9)$, or $(0,10)$, and
\item
$\pencil$ has a section $l^*$ as in \autoref{prop.4-6}\iref{4.6-l*}.
\endroster
Conversely, if $\pencil$ has a section $l^*$
as in \autoref{prop.4-6}\iref{4.6-l*}, then $p+q=10$,
one has $\tP/P=\CG3$, and
Statements~\iref{4.6-plane} and~\iref{4.6-minus} of \autoref{prop.4-6} also
hold.

A section $l^*$ as
above
\rom(or, equivalently, a
geometric index~$3$ extension $\tP\supset P$\rom) is unique up to
automorphism.
\endproposition

\proof
The proof repeats literally that of \autoref{prop.4-6}; the section $l^*$ is
the class $\frac13[\Go+(p-4)\Gl]$ (\cf. also \autoref{lem.skew}\iref{skew.2}).
A direct computation shows that,
whenever the pivot $\tP\ni l^*$, one has $\tP/P=\CG3$, \ie, no further finite
index extension satisfies the conditions of \autoref{pre-conf}.
\endproof

By \autoref{obs.3-torsion}, any $2$-torsion element $\Ga\in\tP/P$ is of the
form $\sum3\Gn_k$, where the index $k$ runs over a certain subset
$\supp\Ga\subset\fib_1\pencil$, called the \emph{support} of $\Ga$.
It is clear that $\supp(\Ga+\Gb)$ is the symmetric difference
$(\supp\Ga)\bigtriangleup(\supp\Gb)$.

\proposition\label{prop.2.pivot}
Let $\pencil$ be a pencil of type $(p,q)$, and let $\Ga\in\tP/P$ be a nonzero
$2$-torsion element. Then
\roster
\item\label{2.pivot.=8}
one has $\ls|\supp\Ga|=8$ and, in particular, $q\ge8$\rom;
\item\label{2.pivot.section}
$\ls|\supp\Ga\cap\scom1(s)|=0$, $2$, or~$4$ for any section~$s$ of~$\pencil$.
\endroster
Besides, the $2$-torsion of the pivot is as follows\rom:
\roster*
\item
$(\CG2)\oplus(\CG2)$ if $(p,q)=(0,12)$,
\item
$\CG2$ if $(p,q)=(0,11)$, $(1,10)$, or $(2,9)$,
\item
$0$ or $\CG2$ in all other cases with $q\ge8$.
\endroster
A geometric index~$2$ \rom(index~$4$ in the case $q=12$\rom) extension
$\tP\supset P$ is unique up to automorphism.
\endproposition

\proof
Clearly, $\ls|\supp\Ga|=0\bmod4$;
hence, $\ls|\supp\Ga|=8$ or $12$ by \autoref{lem.sum.4}.
The last statement is proved by a direct computation using
\autoref{th.Nikulin}. In particular, it follows that, in the case $q=12$,
there are three distinct nonzero elements and, hence, none of them can have
support of length~$12$. This proves Statement~\iref{2.pivot.=8}.
For statement~\iref{2.pivot.section}, it suffices to consider the minimal
pencil of type $(0,8)$, so that $\supp\Ga=\fib_1\pencil$.
Then, clearly, $\num1(s)$ is even, as otherwise $s\notin P$, and the values
$\num1(s)=6$ and $8$ are ruled out by \autoref{pre-conf}\iref{ex-div} and
\iref{ell-pencil}, respectively.

The uniqueness is immediate. In the case of index~$2$, an extension is
determined by a choice of the octet $\supp\Ga\subset\fib_1\pencil$.
If $q=12$, three octets $\supp\Ga_i\subset\fib_1\pencil$, $i=1,2,3$, should
be chosen so that $\ls|\supp\Ga_i\cap\supp\Ga_j|=4$ whenever $i\ne j$.
This choice is equivalent to partitioning $\fib_1\pencil$ into three
quadruples.
\endproof

\corollary[\cf. \autoref{lem.skew}\iref{skew.2}]\label{cor.10.fibers}
If a pencil~$\pencil$ has a section~$s$ intersecting ten fibers of~$\pencil$,
then $\pencil$ has no other fibers.
\endcorollary

\proof
Assuming that $\pencil$ is of type $(0,11)$,
\autoref{prop.2.pivot}\iref{2.pivot.section} applied to~$s$ and
the nontrivial element $\Ga\in\tP/P$ leads to a contradiction.
The existence of~$\Ga$ is also
guaranteed by \autoref{prop.2.pivot}.
\endproof

As another consequence of the results of this section, the type $(p,q)$ and
the primitivity bit almost determine a geometric realization $P\to\L$ up to
isomorphism. The pivot $\tP/P$ may (must if $q>10$) have $2$-torsion if and
only if $q\ge8$ (see \autoref{prop.2.pivot}), and it may
(must if $(p,q)=(4,6)$ or $p=6$ and $q>0$) have
$3$-torsion if and only if $p=6$ (see \autoref{prop.6-0}) or $p+q=10$ (see
Propositions \ref{prop.4-6} and~\ref{prop.3.pivot}). The
case $p+q=10$ and $q\ge8$
is exceptional: here, the pivot may be trivial, $\CG2$, or $\CG3$,
\ie, there are three geometric realizations $P\to\L$.
In this latter case, it makes sense to subdivide the type $(p,q)\impr$ into
$(p,q)^2$ and $(p,q)^3$.

\conjecture
The pivot $\tP/P$ has $3$-torsion if and only if the axis of the pencil is a
line of the second kind in the sense of Segre~\cite{Segre}.
\endconjecture

\section{Counting sections of pencils}\label{S.counting}

The goal of this section is a computer aided estimate on the size of a
geometric configuration containing a pair of large obverse pencils.
Even though most extra restrictions in the ``counting'' lemmas seem purely
technical, for the moment we do need them to keep the computation under
control.

\subsection{The algorithm}\label{s.counting}
Fix a pencil $\pencil:=\pencil(l)$ of type~$(p,q)$ and a section~$s_0$ of~$\pencil$.
Let $\bs_0:=[s_0]\in\coord_{p,q}$ and denote by $\G(\bs_0)$ the
stabilizer of~$\bs_0$.
(Up to automorphism, there are $q+1$ possibilities for $\bs_0$;
we usually choose for $\bs_0$ the vector with several last $1$-coordinates equal
to~$1$ and all other coordinates equal to~$0$.)
More sections $s_1,s_2,\ldots$ are added one by one, building the
obverse pencil
$\pencil(s_0)$. Thus, we assume that
\[
s_0\cdot s_i=1\quad\text{and}\quad s_i\cdot s_j=0\quad
 \text{for $i>j\ge1$},
\label{eq.intr}
\]
\ie, all new sections are in separate fibers of~$\pencil(s_0)$.
Our goal is adding sufficiently many sections, so that, in the resulting
configuration, $\pencil$ is still a maximal pencil and
the multiplicity and valency of~$s_0$
satisfy certain prescribed bounds
\[*
\pmin\le\mult s_0\le\pmax,\qquad
 \vmin\le\val s_0\le\vmax.
\]

It is essential that most of the time we deal with coordinates rather than
sections themselves: we choose certain elements $\bs_i\in\coord_{p,q}$
and consider the pre-configuration
\[*
S_k:=P(\bs_0,\ldots,\bs_k)=(\tP+\Z s_0+\ldots+\Z s_k)/\!\ker,
\]
where $[s_i]=\bs_i$ for all $i\ge0$ and the intersection matrix of~$P$ is
extended using~\eqref{eq.intr} and the definitions of sections and
coordinates.
By \autoref{cor.obverse}, for each $i\ge1$ we must have
$\bs_i\in\coord_0(\bs_0)\cup\coord_1(\bs_0)$, where
\[*
\coord_r(\bs_0):=\bigl\{\bs\in\coord_{p,q}\bigm|\com(\bs\cs\bs_0)=r\bigr\};
\]
furthermore,
$s_i$ is contained in a $1$-fiber of $\pencil(s_0)$
if and only if $\bs_i\in\coord_0(\bs_0)$.

Once a lattice~$S_k$ has been constructed, we denote by
\[*
G_{k}:=\OG_h(S_{k},l,s_0)
\]
the group of its isometries preserving~$h$, $l$ and~$s_0$.
The computation of this group is discussed in \autoref{obs.group} below.
(At the expense of a certain overcounting, we compute separately the
stabilizers in $\SG{3}^p$ and $\SG{p}\times\SG{q}$.)

The algorithm runs in several steps.

\subsubsection{Step 1\rom: collecting the candidates}\label{obs.pre}
Assume $S_{k-1}$ known
and denote by
$\barS_{k-1}$ the multiset $\{[s]\,|\,s\in S_{k-1}(\pencil)\}$.
The group $G_{k-1}$
acts on
$\coord_0(\bs_0)\cup\coord_1(\bs_0)\sminus\barS_{k-1}$ and, when passing
to~$S_k$,
it suffices to take for~$\bs_k$ one representative from each orbit of this
action. We can also assume that all explicit $3$-fibers are added first
and avoid adding too many $3$-fibers:
\roster
\item\label{pre.plane}
$\bs_k\in\coord_1(\bs_0)$ if $\mult s_0<\pmin$ and
$\bs_k\in\coord_0(\bs_0)$ if $\mult s_0\ge\pmax$.
\endroster
There is an obvious injective map from the set of $3$-fibers of
$\pencil(s_0)$ to $\fib\pencil$ (each $3$-fiber contains a unique line
$a\in\pencil$); this map should remain injective:
\roster[\lastitem]
\item\label{pre.injective}
if $s\in S_{k-1}(\pencil)$ is contained in a $3$-fiber of $\pencil(s_0)$,
then $\com(\bs_0\cs\bs_k\cs s)=0$.
\endroster
Other restrictions taken into account when choosing $\bs_k$ are as
follows:
\roster[\lastitem]
\item\label{pre.det}
\autoref{lem.det} (in fact, we check that $[-r_{ij}]$ is negative
semi-definite);
\item\label{pre.rank}
$\rank[-r_{ij}]+2p+q\le18$ (as $S_k$ should admit an embedding to~$\L$);
\item\label{pre.4}
$\cpre_4(\bs_k)$:
$\com(\bs_k\cs\bs)\le4$ for any $\bs\in\barS_{k-1}$,
see \autoref{lem.s.s=0}\iref{com.5};
\item\label{pre.3}
$\cpre_3(\bs_k)$:
$\com(\bs_k\cs\bs'\cs\bs'')\le3$ for all $\bs'\ne\bs''\in\barS_{k-1}$,
see \autoref{lem.s.s=0}\iref{com.4},~\iref{com.3};
\item\label{pre.h}
$\cpre_h(\bs_k)$: if $\com(\bs_k\cs\bs)=4$ for some $\bs\in\barS$, then
$\cpre_4(\bs')$ and $\cpre_3(\bs')$ hold,
where $\bs':=\one\om\bs\om\bs_k$, see \autoref{lem.s.s=0}\iref{com.4};
\item\label{pre.p}
if $\bs_k\in\coord_1(\bs_0)$, then
$\cpre_4(\bs')$, $\cpre_3(\bs')$, and $\cpre_h(\bs')$ hold
for $\bs':=\one\om\bs_0\om\bs_k$.
\endroster
In cases~\iref{pre.h} and~\iref{pre.p}, we also exclude from further
consideration the $G_{k-1}$-orbit of the respective section~$\bs'$, as
its presence in $\barS_k$ would imply the presence of~$\bs_k$.

\subsubsection{Step 2\rom: validating a section $\bs_k$}\label{obs.post}
Now, for each candidate~$\bs_k$
collected at the previous step,
we compute
the pre-configuration $S_k=(S_{k-1}+\Z s_k)/\!\ker$, consider the orthogonal
complement $h^\perp$ in $S_k$, and use \GAP~\cite{GAP4}
function {\tt ShortestVectors} to compute the sets
$\SV_2(S_k)$ and $\SV_4(S_k)$, where
\[*
\SV_r(S_k):=\{v\in h^\perp\subset S_k\,|\,v^2=-r\}.
\]
(Note that
the lattice $S_k$ is hyperbolic,
hence $h^\perp$ is elliptic,
by \autoref{obs.pre}\iref{pre.det}.)

A candidate $\bs_k$ is rejected as invalid (not leading to a
geometric configuration) if one of the following holds:
\roster
\item\label{post.ex-div}
$\SV_2\ne\varnothing$, see \autoref{pre-conf}\iref{ex-div};
\item\label{post.ell-pencil}
there is $v\in\SV_4$ such that $v+h\in2S_k$,
see \autoref{pre-conf}\iref{ell-pencil}.
\endroster
Otherwise, the new set of sections $S_k(\pencil)$ is computed \via
\[*
S_k(\pencil)=\{v+l\,|\,v\in\SV_4,\ v\cdot l=2\}.
\]
At this point, the full intersection matrix is known, and we can compute and
record
the set
\[*
\sect_k:=\sect(\bs_0,\ldots,\bs_k)=\pencil(s_0)\cap S_k(\pencil),
\]
including types of
the fibers. This set is used for the further validation. Namely, we
reject~$\bs_k$ if
\roster[\lastitem]
\item\label{post.planes}
$\mult s_0>\pmax$ (too many $3$-fibers), or
\item\label{post.lines}
$\val s_0>\vmax$ (too many lines in $\pencil(s_0)$), or
\item\label{post.fiber}
there is a pair $s'\ne s''\in\sect_k$ such that $s'\cdot s''=1$ and
$\com(s_0\cs s'\cs s'')=0$,
see \autoref{cor.obverse}, or
\item\label{post.maximal}
any other type specific restriction is not satisfied
(whenever used, this extra restriction is specified
explicitly in the respective proof).
\endroster
To conserve space, for each candidate~$\bs_k$ that passed the validation, we
record
\roster*
\item
the elements $\bs_0,\ldots,\bs_k\in\coord_{p,q}$,
\item
the multiset $\barS_k$ (sections in terms of coordinates), and
\item
the image $\bar\sect_k\subset\barS_k$ of $\sect_k$ under the coordinate
map,
\endroster
disregarding all other information.

\subsubsection{Step~3\rom: eliminating repetitions}\label{obs.group}
Before further processing, we eliminate the repetitions in the obtained list
of lattices~$S_k$ by retaining a single representative of each orbit of the
$\G(\bs_0)$-action.
To compute the orbits or, equivalently, the stabilisers $G_k$,
we use on of the following two approaches.
\roster
\item\label{orbits.partial}
The stabilizers are computed \via\ $G_k=\stab\bar\sect_k\subset\G(\bs_0)$,
and the lattices are compared by means of the orbits of $\bar\sect_k$.
This approach works if each $S_k$ is exactly as in the construction above,
\ie, generated over $\tP_{p,q}$ by the set
$\{s_0\}\cup\sect_k$, on which {\em the intersection matrix is known}.
\item\label{orbits.full}
The stabilizers are computed \via\ $G_k=\stab\barS_k\subset\G(\bs_0)$,
and the lattices are compared by means of the orbits of $\barS_k$.
This approach applies if each $S_k$ is known to be combinatorially rigid.
\endroster
By default, we use approach~\iref{orbits.partial}.

\subsubsection{Step 4\rom: checking the $\L$-realizability}\label{obs.K3}
For each configuration~$S_k$ obtained at Step~3, we check if it admits
a geometric $\L$-realization.
To this end, we start with the lattice~$S_k$ itself and apply
\autoref{th.Nikulin}
to see if $S_k$ admits a primitive $\L$-realization.
If not, we replace~$S_k$ with a finite index extension
$\tS_k\supset S_k$ defined by an
isotropic vector $v\in\discr S_k$ of prime order.
(This and subsequent steps are repeated
for each isotropic vector found.) The new lattice~$\tS_k$ is rejected if
it fails to satisfy one of the conditions in \autoref{obs.post};
otherwise, we apply
\autoref{th.Nikulin}
again. The algorithm stops when
a primitive embedding is found (and then $S_k$ is accepted)
or all isotropic vectors are exhausted; in the
latter case, the original lattice~$S_k$ is rejected as not
admitting a geometric $\L$-realization.
Admittedly ineffective, this algorithm works reasonably well for the vast
majority of configurations.

\subsubsection{Increasing the rank}\label{obs.rank}
We repeat Steps~1--4 above until either nothing else can be added or the
desired bounds $\mult s_0\ge \pmin$, $\val s_0\ge \vmin$ have been
achieved. Most lattices~$S_k$ obtained have rank~$20$ and, hence, each geometric
configuration containing~$S_k$ is a finite index extension of~$S_k$.
In
the
cases where
$\rank S_k\le19$, we keep~$S_k$ on the
list, but
we allow
also
the addition of an extra
section~$s_{k+1}$ disjoint from~$s_0$. (Certainly, in this case we have to
switch to approach~\iref{orbits.full} in \autoref{obs.group}, \ie, we need
to know that the configurations obtained are combinatorially rigid.
If the latter property cannot be asserted, configurations with
extra
sections are excluded from Step~3.)
This time, we have $s_0\cdot s_{k+1}=0$,
but the intersections $\iota_i:=s_i\cdot s_{k+1}$, $i=1,\ldots,k$,
should be given as part of the input;
for each
pair $(\bs_{k+1},[\iota_i])$,
we check conditions \iref{pre.det}--\iref{pre.p} in
\autoref{obs.pre}, requiring in addition that
$\rank S_{k+1}>\rank S_k$,
\ie, the same lattice cannot be obtained as a finite index
extension of~$S_k$.
Then, Steps~2--4 are repeated and, at Step~2,
we require that
\roster
\item\label{rank.=}
the valency of~$s_0$ in $S_{k+1}$ must be equal to that in~$S_k$,
\endroster
as otherwise the same lattice can be obtained by adding a section
intersecting~$s_0$.

\subsubsection{Final step\rom: computing $\L$-realizations}\label{obs.L}
There remains to enumerate, for each lattice~$S_k$,
its geometric $\L$-realizations.
This is done similar to
\autoref{obs.K3}, except that we do not stop at the first valid realization;
on the other hand, we
require that
\roster
\item\label{L.=}
the valency of~$s_0$ in $\tS_{k}$ must be equal to that in~$S_k$,
\cf. \autoref{obs.rank}\iref{rank.=}.
\endroster
At this step,
for all consecutive extensions $S_k=\tS_k^0\subset\tS_k^1\subset\ldots$ of
prime index, we can also check that $\ls|\Fn(\tS_k^i)|>\ls|\Fn(S_k^{i-1})|$;
this fact implies that all configurations
found are generated by sections.

For each finite index extension $\tS_k\supset S_k$ found in this way,
{\em assuming that $\pencil$ is maximal in $\tS_k$}, we have
\[
\ls|\Fn(\tS_k)|=\ls|\tS_k(\pencil)|+3p+q+1.
\label{eq.lines}
\]
In extreme cases (when too many lines have been found),
we recompute the maximal pencil \via
\[*
\pencil(l)=\{v+l\,|\,v\in\SV_6(\tS_k),\ v\cdot l=3\}
\]
and compute the pencil structure of~$\tS_k$.
(The computation of $\SV_6$ is rather expensive and we try to avoid it as
much as possible.)

\subsection{A list of configurations}\label{obs.pencils}
For further references, we collect in
\autoref{tab.list}
a list of large
configurations found in the experiments.
(We list all known configurations with more that $48$ lines;
for the moment,
we do not assert that the list is complete.)
The notation
refers to certain particular configurations
found in the
computation.
We will also speak
about configurations of \emph{type $\config{X}_*$, $\config{Y}_*$}, \etc.,
meaning that the pencil
structures of the two configurations are equal.
The configurations
marked with a $^*$ in the table (most notably, the \aconfig{Y}-series)
admit totally reflexive $\L$-realizations; the others do
not.
One has
$\rank\aconfig{Z}_*=19$;
the other configurations listed in the table are of
rank~$20$. There is no particular difference between \aconfig{X} and
\aconfig{Q}.

\table
\caption{Known large geometric configurations}\label{tab.list}
\def\*{\relax\llap{$^*$}}
\hbox to\hsize{\hss\vbox{\halign{\strut\quad#\quad\hss&\hss#\hss\quad&#\hss\ \cr
\noalign{\hrule\vspace{2pt}}%
 \strut\hss$S$\hss&$\ls|\Fn|$&Pencil structure (see \autoref{s.invariants}), reference, remarks
\cr
\noalign{\vspace{1pt}\hrule\vspace{2pt}}
\aconfig{64}&64&
$\PS{[ [ 6, 0 ], 16 ], [ [ 4, 6 ], 48 ]}$,
see \autoref{th.6-q} and \autoref{th.4-6}
\cr
\aconfig{60}&60&
$\PS{[ [ 6, 2 ], 10 ], [ [ 4, 4 ], 30 ], [ [ 3, 7 ], 20 ]}$,
see \autoref{th.6-q}
\cr
\aconfig{60.2}&60&
$\PS{[ [ 4, 5 ], 60 ]}$,
see \autoref{lem.4-5}
\cr
\aconfig{56}&56&%
$\PS{[ [ 4, 6 ], 8 ], [ [ 4, 4 ], 32 ], [ [ 2, 8 ], 16 ]}$,
see \autoref{th.4-6}
\cr
\*\aconfig{56.real}&56&%
$\PS{[ [ 4, 4 ], 32 ], [ [ 3, 7 ], 24 ]}$,
see \autoref{lem.3-7}
\cr
\aconfig{q56}&56&
$\PS{[ [ 4, 4 ], 24 ], [ [ 3, 7 ], 32 ]}$,
see \autoref{lem.3-7}
\cr
\aconfig{54}&54&
$\PS{[ [ 6, 2 ], 4 ], [ [ 4, 6 ], 6 ], [ [ 4, 4 ], 6 ], [ [ 4, 2 ], 24 ],
 [ [ 2, 8 ], 12 ], [ [ 0, 10 ], 2 ]}$,
see \autoref{th.6-q}
\cr
\aconfig{q54}&54&
$\PS{[ [ 4, 4 ], 24 ], [ [ 4, 3 ], 24 ], [ [ 0, 12 ], 6 ]}$,
see \autoref{lem.4-4}
\cr
\aconfig{52.1}&52&%
$\PS{[ [ 6, 0 ], 1 ], [ [ 4, 4 ], 12 ], [ [ 4, 3 ], 12 ], [ [ 4, 2 ], 3 ],
 [ [ 3, 5 ], 18 ], [ [ 0, 12 ], 6 ]}$,
see \autoref{th.6-q}
\cr
\aconfig{52.2}&52&%
$\PS{[ [ 6, 0 ], 1 ], [ [ 4, 4 ], 9 ], [ [ 4, 3 ], 18 ], [ [ 3, 5 ], 18 ],
 [ [ 0, 12 ], 6 ]}$,
see \autoref{th.6-q}
\cr
\aconfig{52.4}&52&%
$\PS{[ [ 4, 6 ], 10 ], [ [ 3, 5 ], 40 ], [ [ 0, 10 ], 2 ]}$,
see \autoref{th.4-6}
\cr
\aconfig{52.5-3}&52&%
$\PS{[ [ 5, 3 ], 8 ], [ [ 3, 5 ], 32 ], [ [ 2, 8 ], 12 ]}$,
see \autoref{th.5-3}
\cr
\*\aconfig{52.5}&52&%
$\PS{[ [ 4, 6 ], 2 ], [ [ 4, 4 ], 16 ], [ [ 3, 5 ], 20 ], [ [ 2, 8 ], 14 ]}$,
see \autoref{th.4-6}
\cr
\*\aconfig{52.real}&52&%
$\PS{[ [ 4, 5 ], 8 ], [ [ 4, 3 ], 12 ], [ [ 3, 6 ], 16 ], [ [ 2, 7 ], 16 ]}$,
see \autoref{lem.4-5}
\cr
\*\aconfig{52.0}&52&%
$\PS{[ [ 6, 0 ], 4 ], [ [ 4, 4 ], 12 ], [ [ 4, 2 ], 24 ], [ [ 2, 8 ], 12 ]}$,
see \autoref{th.6-q}; $\rank\config{52.0}=19$
\cr%
\aconfig{q52.1}&52&
$\PS{[ [ 4, 4 ], 16 ], [ [ 4, 3 ], 16 ], [ [ 4, 2 ], 16 ], [ [ 0, 12 ], 4 ]}$,
see \autoref{lem.4-4}
\cr
\aconfig{q52.2}&52&
$\PS{[ [ 4, 4 ], 8 ], [ [ 4, 3 ], 32 ], [ [ 4, 2 ], 8 ], [ [ 0, 12 ], 4 ]}$,
see \autoref{lem.4-4}
\cr
\aconfig{51}&51&%
$\PS{[ [ 6, 2 ], 1 ], [ [ 5, 3 ], 6 ], [ [ 4, 3 ], 3 ], [ [ 3, 6 ], 6 ],
         [ [ 3, 4 ], 8 ], [ [ 2, 7 ], 27 ]}$%
, see \autoref{th.6-q}
\cr%
\aconfig{50.1}&50&%
$\PS{[ [ 6, 1 ], 1 ], [ [ 4, 4 ], 9 ], [ [ 4, 3 ], 9 ], [ [ 4, 2 ], 9 ],
         [ [ 3, 4 ], 18 ], [ [ 0, 12 ], 3 ], [ [ 0, 10 ], 1 ]}$%
, see \autoref{th.6-q}
\cr%
\aconfig{50.2}&50&%
$\PS{[ [ 6, 1 ], 1 ], [ [ 4, 4 ], 6 ], [ [ 4, 3 ], 15 ], [ [ 4, 2 ], 6 ],
         [ [ 3, 4 ], 18 ], [ [ 0, 12 ], 3 ], [ [ 0, 10 ], 1 ]}$%
, see \autoref{th.6-q}
\cr%
\aconfig{50}&50&%
$\PS{[ [ 5, 3 ], 4 ], [ [ 4, 4 ], 8 ], [ [ 3, 5 ], 16 ],
  [ [ 2, 8 ], 4 ], [ [ 2, 6 ], 18 ]}$,
see \autoref{th.5-3}
\cr%
\aconfig{50.0}&50&
$\PS{[ [ 4, 4 ], 10 ], [ [ 3, 5 ], 40 ]}$,
see \autoref{lem.4-4}; $\rank\config{50.0}=19$
\cr
\aconfig{49}&49&
$\PS{[ [ 6, 0 ], 1 ], [ [ 4, 3 ], 18 ], [ [ 4, 2 ], 9 ], [ [ 3, 4 ], 18 ],
[ [ 0, 12 ], 3 ]}$%
, see \autoref{th.6-q}%
; $\rank\config{49}=19$
\cr
\*\aconfig{48}&48&%
$\PS{[ [ 5, 1 ], 2 ], [ [ 3, 7 ], 6 ], [ [ 3, 5 ], 24 ], [ [ 2, 6 ], 12 ],
[ [ 1, 9 ], 4 ]$},
see \autoref{lem.3-7}
\cr
\*\aconfig{48.2}&48&
$\PS{[ [ 4, 4 ], 4 ], [ [ 4, 2 ], 16 ], [ [ 3, 6 ], 8 ], [ [ 2, 7 ], 12 ],
 [ [ 2, 6 ], 8 ]}$,
see \autoref{lem.4-4}
\cr
\noalign{\vspace{2pt}\hrule}\crcr}}\hss}
\endtable

\theorem\label{th.uniqueness}
A
geometric configurations of
each type
listed in \autoref{tab.list}
is unique up to isomorphism.
\endtheorem

\proof
Each configuration~$S$
satisfies the hypotheses of the
respective
classification statement cited in the table (with pencils of type
$(6,0)\prim$ ruled out by \autoref{th.6-0.primitive}),
and the uniqueness follows from the classification.

Indeed,
the essential part of the hypotheses is the existence of a certain pair of
obverse pencils.
Let $v:=\max\{\val l\,|\,l\in\Fn(S)\}$, and denote ny~$n$ the number of lines
of valency~$v$.
If $v>18$, then, in view of \eqref{tablichka},
the configuration is covered by \autoref{th.6-q}.
If $n\ge5$ or $n\ge4$ and $\ls|\Fn(S)|<4v-8$, then,
referring in the latter case to \autoref{cor-Segre}, we obtain a pair of
skew lines of valency~$v$, which suffices for all statements.
In the remaining four cases (\config{52.1}, \config{52.2}, \config{52.5},
and \config{49}), a similar argument gives us a pair of lines of valency
$v=18$ and $\ge15$.
\endproof

Among others,
\autoref{tab.list} lists
all geometric configurations~$S$ containing a pair of obverse
pencils $\pencil_1$, $\pencil_2$ such that
\[*
\ls|\Fn(S)|>48\quad\text{and}\quad
\ls|\pencil_1|+\ls|\pencil_2|\ge33.
\]

\subsection{Pencils of type $(6,\q)$}\label{s.6.0}
For the moment, $(6,\q)\impr$-configurations is the only class that is
sufficiently well
understood. The properties of such configurations
are summarized in the next
theorem.

\theorem\label{th.6-q}
There are $300$ isomorphism classes of $(6,q)\impr$-pairs\rom:
\roster*
\item
for $q=0$\rom: $62$ classes, of which $43$ are totally reflexive\rom;
\item
for $q=1$\rom: $107$ classes, none totally reflexive\rom;
\item
for $q=2$\rom: $131$ classes, none totally reflexive.
\endroster
Let $(S,\pencil)$ and $(S',\pencil')$ be two $(6,\q)\impr$-pairs. Then\rom:
\roster
\item\label{6-q.rigid}
$S$ is generated by sections and combinatorially rigid\rom;
\item\label{6-q.ls}
with one exception,
one has
$(S',\pencil')\cong(S,\pencil)$
if and only if
$\Ls(S')=\Ls(S)$\rom;
\item\label{6-q.known}
either one has $\ls|\Fn(S)|<52$ or
$S\cong\config{64}$, \config{60}, \config{54}, \config{52.1}, \config{52.2},
or~\config{52.0}.
\endroster
Furthermore, for any
$n\in\{19,\ldots,52,54,60,64\}$, there is a
$(6,\q)\impr$-configuration $S$ such that $\ls|\Fn(S)|=n$.
\endtheorem

As an addendum to \autoref{th.6-q}\iref{6-q.ls}, note that, with the
exception of eleven pairs, any two distinct $(6,\q)\impr$-configurations are
distinguished by the pencil structure.

\proof[Proof of \autoref{th.6-q}]
We start with a pencil~$\pencil$ of type $(6,0)\impr$ and apply the algorithm
of \autoref{s.counting}, introducing a number of modifications:
\roster*
\item
we do not fix a section~$\bs_0$ and use
the group $\tilde\G$ instead of $\G(\bs_0)$,
see~\eqref{eq.6-0.group};
the intersection matrices are computed by means of \autoref{lem.6-0.intr};
\item
at Step~1, all restrictions are lifted: instead, we construct the
``convex hull'' (in the sense of \autoref{lem.convex}) of the set
$\barS_{k-1}\cup\bs_k$ and check whether the resulting set~$\bar\CS_k$
satisfies \autoref{lem.plane};
certainly,
$\bs_k$ must satisfy~\eqref{eq.6-0.coord};
\item
at Step~2, all restrictions except~\iref{post.ex-div}
and~\iref{post.ell-pencil} are lifted;
\item
at Step~3, approach~\iref{orbits.full} can be used due to
\autoref{cor.6-q.combinatorial};
\item
since all sets of sections are to be tried,
we replace condition~\iref{L.=} in \autoref{obs.L} with
$\ls|\Fn(\tS_k)|=\ls|\Fn(S_k)|$.
It turns out that such extensions do not exist;
hence, any geometric configuration is generated by sections.
\endroster
As a result, we obtain $84$ configurations (of which $25$ are extremal with
respect to inclusion) generated by sections of~$\pencil$; in these
configurations, $\pencil$ is not always maximal.
Then, we try to add up to two extra $1$-fibers. The procedure is similar to
\autoref{obs.rank}: we specify the intersection of the fiber added with
sections generating~$S_k$ and repeat Steps~1--4 of the algorithm; a new
configuration $S_k'$ is accepted only if $\ls|\Fn(S_k')|=\ls|\Fn(S_k)|$.
Repetitions are eliminated using approach~\iref{orbits.full} of
\autoref{obs.group} and appropriate subgroup \smash{$\tilde\G\subset\G_{6,q}$}.

All other statements of the theorem follow directly from the classification.
\endproof

\theorem\label{th.6-0.primitive}
There are $69$ isomorphism classes of $(6,0)\prim$-pairs $(S,\pencil)$
admitting a section
$s_0\in S(\pencil)$ such that $15\le\val s_0\le18$.
Let $(S,\pencil)$ and $(S',\pencil')$ be two such pairs. Then\rom:
\roster
\item\label{6-0.rigid}
$S$ is generated by sections and combinatorially rigid\rom;
\item\label{6-0.ls}
$(S',\pencil')\cong(S,\pencil)$
if and only if
$\Ls(S')=\Ls(S)$\rom;
\item\label{6-0.known}
one has $\ls|\Fn(S)|<44$.
\endroster
\endtheorem

\proof
The sections are enumerated using the algorithm of
\autoref{s.counting}, letting
\[
\pmin=2,\quad \pmax=6,\quad \vmin=15,\quad \vmax=18.
\label{eq.param.18}
\]
Here, the lower bound $\pmin=2$ follows from~\eqref{tablichka}, and the
seemingly redundant upper bound $\pmax=6$ helps us eliminate a number of
configurations before any further processing.
We introduce also a few modifications to the algorithm.
First,
by \autoref{cor.6-0.combinatorial}, we can use approach~\iref{orbits.full} in
\autoref{obs.group}: this is necessary since some of the configurations~$S_k$
with $\val s_0\ge16$ have rank~$19$, see \autoref{obs.rank}.
Besides, we can
\roster*
\item
use \autoref{lem.6.0-sections}\iref{i.(6,p).0-intr} for
condition~\iref{post.maximal} in \autoref{obs.post}, and
\item
in \autoref{obs.rank}, consider only the candidates~$\bs_{k+1}$ satisfying
$1\le\com(\bs_{k+1}\cs\bs_0)\le4$,
see \autoref{lem.6.0-sections}\iref{i.(6,p).0} and
\autoref{lem.s.s=0}\iref{com.5}.
\endroster
we obtain $81$ configurations, each with a distinguished section~$s_0$.
Switching to the full automorphism group~$\G_{6,0}$ and resorting reduces the
list down to $69$ classes.
The maximal number of lines in the configurations obtained is~$44$.
\endproof

\subsection{Pencils of type $(4,\q)$}\label{s.4-6}
A complete classification of $(4,6)$-configurations also seems feasible;
however, for the moment we confine ourselves to a partial statement.
similar
to \autoref{th.6-0.primitive}.

\theorem\label{th.4-6}
There are $195$
isomorphism classes of $(4,6)$-pairs $(S,\pencil)$ admitting
a section~$s_0\in S(\pencil)$
such that $15\le\val s_0\le18$.
If $(S,\pencil)$ is such a pair, then\rom:
\roster
\item\label{4-6.rigid}
$S$ is generated by sections and combinatorially rigid\rom;
\item\label{4-6.known}
either one has $\ls|\Fn(S)|\le48$ or $S\cong\config{64}$, \config{56},
\config{54}, \config{52.4}, or \config{52.5}.
\endroster
\endtheorem

\proof
First, assume that $\mult l^*\le2$, hence $s_0\ne l^*$. We need to
consider seven cases: $\num1(s_0)\in\{0,\ldots,4\}$ and
$s_0\cdot l^*=0$ or~$1$ for the
first two values $\num1(s_0)=0,1$.
In each case, we employ the algorithm of \autoref{s.counting},
using parameters~\eqref{eq.param.18},
restricting the candidates in \autoref{obs.pre} to satisfy~\eqref{eq.l*}, and
imposing the restriction
$\ls|\barS^*|\le4$, see~\eqref{eq.barS*}, as condition~\iref{post.maximal} in
\autoref{obs.post}.
All pairs obtained are rigid by \autoref{lem.4-6.combinatorial}, and
resorting the list with the full automorphism group $\G_{4,6}$ reduces it to
$20$ classes.

Let $s_0=l^*$.
To avoid complications with
large pivots,
we start with a manual classification of configurations $S\supset\pencil$
generated by up to four sections $s_i$ such that $s_i\cdot l^*=1$ and
$\num1(s_i)=0$.
It is easily shown that, in addition to~$P$ itself, there are six
isomorphism classes of such
configurations~$S$, each admitting a unique, up to automorphism, geometric
finite index extension $\tS\supset S$.
Briefly, they are as follows:
\roster*
\item
$1$ class with $\mult l^*=1$, $\rank S=17$, and $\ell_3(\tS/S)=0$,
\item
$2$ classes with $\mult l^*=2$, $\rank S=18$, and $\ell_3(\tS/S)=1$,
\item
$1$ class with $\mult l^*=3$, $\rank S=19$, and $\ell_3(\tS/S)=2$,
\item
$1$ class with $\mult l^*=4$, $\rank S=19$, and $\ell_3(\tS/S)=2$,
\item
$1$ class with $\mult l^*=4$, $\rank S=20$, and $\ell_3(\tS/S)=3$.
\endroster
Starting, instead of~$\tP$,
with one of these geometric configurations~$\tS$, we build a
separate list, replacing $\G(\bs_0)$ with $\OG_h(\tS,l)$ and inhibiting
sections~$\bs$ with $\num1(\bs)=0$ at Step~1.
Running the algorithm, we obtain
a large number of configurations (due to the lack of sorting in
\autoref{obs.rank} and \autoref{obs.L}).
All but one are rigid by
\autoref{lem.4-6.combinatorial},
and the remaining one
has an ``ambiguous'' pair of sections $s_1,s_2$, but the assumptions
$s_1\cdot s_2=0$ or~$1$ result in configurations~$S_0$, $S_1$
with non-isomorphic sets of
sections (in fact, $S_0$ is generated by $\barS_0\sminus\barS_0^\circ$,
whereas $S_1$ is not; this phenomenon is similar to \autoref{lem.com2=1}).
Thus, \emph{a posteriori},
all configurations are rigid;
switching to approach~\iref{orbits.full}
in \autoref{obs.group} and
resorting the list
reduces it to $175$ classes.
\endproof

\lemma\label{lem.4-5}
If a $(4,5)$-pair $(S,\pencil)$ admits a section
$s_0$ such that $16\le\val s_0\le17$, then
either one has $\ls|\Fn(S)|\le48$
or $S\cong\config{60.2}$ or~\config{52.real}.
Furthermore, a geometric configuration
of type \config{60.2} is unique up to isomorphism.
\endlemma

\proof
We apply the algorithm of \autoref{s.counting}, letting
\[
\pmin=3,\quad \pmax=5,\quad \vmin=16,\quad \vmax=17
\label{eq.param.17}
\]
and using for~\iref{post.maximal} in \autoref{obs.post}
the extra requirement
that
$\num1(\bs)\le4$ for any $\bs\in\barS_k$,
see
\autoref{prop.4-5}. We also suppress the sorting
in \autoref{obs.L}, which results in a rather large number of classes in the
case where $\num1(\bs_0)=4$.
Disregarding the pairs $(S,\pencil)$ with $\ls|S(\pencil)|\le30$, we arrive at
a number of configurations of type \config{52.real}
and
several dozens of those of type \config{60.2};
crucial is the fact that
{\em only two
configurations of type \config{60.2}
appear in the case where $\num1(\bs_0)=0$.}

For the uniqueness,
we compute the linking structure of each configuration~$S$ of type
\config{60.2}.
The result is the same for all configurations:
\[*
\Ls(S)=\PS{[ [ 4, 4 ], 150 ], [ [ 5, 3 ], 360 ], [ [ 6, 2 ], 360 ],
  [ [ 7, 2 ], 240 ], [ [ 8, 0 ], 30 ], [ [ 8, 3 ], 120 ] }.
\]
Since $(4,4)\in\Ls(S)$, it follows that $S$
has a pair of skew lines~$l$, $s_0$ such that
$\num1(s_0)=0$ with respect to $\pencil(l)$; in particular, there are at most
two isomorphism classes.

A further computation in (any) one of the
configurations shows that there are at least
two classes of
pairs $\pencil_1,\pencil_2$ such that $\ls|\pencil_1\cap\pencil_2|=4$.
Namely, in each $3$-fiber of~$\pencil_2$, consider the two lines $s',s''$
that are sections of~$\pencil_1$ and compute
$\num1(s')$, $\num1(s'')$ with respect to $\pencil_1$.
The resulting multiset of four unordered pairs is obviously an invariant of
$\pencil_1,\pencil_2$; it turns out to be symmetric, and it can take values
\[
\PS{[ [ 1, 4 ], 2 ], [ [ 2, 3 ], 2 ]}\ \text{($120$ pairs)}\quad\text{or}\quad
\PS{[ [ 1, 4 ], 4 ]}\ \text{($30$ pairs)}.
\label{eq.invariant}
\]
Thus, we conclude that the two classes obtained in the
case $\num1(s_0)=0$ correspond, in fact, to two distinct pairs of obverse
pencils in the same configuration.

All configurations of type \config{52.real}
(obtained in the computation) are isomorphic,
as only one configuration is obtained when $\num1(s_0)=2$ and
each configuration has a pair of obverse pencils
$\pencil_1$, $\pencil_2$ of type $(4,5)$ and such that
$\ls|\pencil_1\cap\pencil_2|=6$.
\endproof

\corollary[of the proofs]\label{cor.4-5.counts}
For any
$n\in\{18,\ldots,48, 52, 54, 56, 60, 64\}$,
there
exists a $(4,\q)$-configuration~$S$
such that $\ls|\Fn(S)|=n$.
If
$n\in\{18,\ldots,47,52\}$,
this configuration~$S$ can be chosen totally reflexive.
\endcorollary

\proof
By Propositions~\ref{prop.4-6} and~\ref{prop.4-5}, we can reliably detect the
maximality of a pencil~$\pencil$ of type $(4,6)$ or $(4,5)$
in a configuration~$S$ by the set of sections
$S(\pencil)$, without recomputing the full set $\Fn(S)$.
Hence, \eqref{eq.lines} applies to any geometric
finite index extension $\tS_k\supset S_k$
accepted in \autoref{obs.K3}; recording the values obtained, we obtain the
first statement of the corollary. The second one is
obtained by using, in addition, \autoref{totally_reflexive}, \cf. the proof
of \autoref{prop.5-2.real}.
\endproof

\subsection{Pencils of type $(5,\q)$}\label{s.5-3}
As in the case $(4,\q)$, we have a partial classification for the maximal
type $(5,3)$ and certain bounds for the submaximal type $(5,2)$.

\theorem\label{th.5-3}
There are $421$
isomorphism classes of $(5,3)$-pairs $(S,\pencil)$ admitting
a section~$s_0\in S(\pencil)$
such that $15\le\val s_0\le18$.
If $(S,\pencil)$ is such a pair,
then either one has $\ls|\Fn(S)|\le48$ or
$S\cong\config{52.5-3}$, \config{51}, or~\config{50}.
\endtheorem

\proof
The
computation runs exactly as outlined in \autoref{s.counting},
with the parameters as in~\eqref{eq.param.18}
and \autoref{lem.5-q.maximal} used for
condition~\iref{post.maximal} in \autoref{obs.post}. (Note that, since the
only pencil that can properly contain~$\pencil$ is that of type $(6,2)$,
\autoref{lem.5-q.maximal} gives us a criterion of maximality of~$\pencil$.)
With two exceptions, all configurations obtained are rigid by
Corollary~\ref{cor.5-3.combinatorial} or~\ref{cor.5-3.pencil},
and we can resort the combined list
(the union over all four values $\num1(\bs_0)=0,\ldots,3$)
using
approach~\eqref{orbits.full} in \autoref{obs.group} and the full
group $\G_{5,3}$.
Each of the two configurations whose rigidity could not be established
differs from all others by its linking structure.
\endproof

\lemma\label{lem.5-2}
If a $(5,2)$-pair $(S,\pencil)$ admits a section
$s_0$ such that $16\le\val s_0\le17$, then
one has $\ls|\Fn(S)|\le48$.
\endlemma

\proof
The computation runs as outlined in \autoref{s.counting},
using parameters as in~\eqref{eq.param.17} and \autoref{lem.5-q.maximal} for
condition~\iref{post.maximal} in \autoref{obs.post}.
There are a few configurations~$S_k$ of rank~$19$,
to which we add extra sections
(see \autoref{obs.rank}) but {\em do not sort the results}, \ie, skip
Step~3.
Apart from several configurations of type~\config{52.5-3} or~\config{50},
one has
$\ls|\tS_k(\pencil)|\le30$ and the statement follows from~\eqref{eq.lines}.
\endproof

\subsection{Pencils of size $16$}\label{s.3-7}
In this section
we deal with
geometric
configurations containing
a pair of obverse
maximal pencils $\pencil:=\pencil(l)$ and $\pencil':=\pencil(s_0)$
such that $\ls|\pencil|=\ls|\pencil'|=16$.
Since we are interested in the configurations themselves rather than triples
$(S,\pencil,\pencil')$, we make several additional assumptions.

First of all, we assume that $\mult l\le\mult s_0$; hence, when applying the
algorithm  outlined in \autoref{s.counting},
we can use the parameters
\[*
\pmin=p:=\mult l,\quad \pmax=5,\quad \vmin=\vmax=16.
\]

The next few restrictions are considered as part of the type specific
condition~\iref{post.maximal} in \autoref{obs.post}; the necessary
computation uses the set $\SV_4(S_k)$.
\roster
\item\label{3-7.val}
We require that $\max\{\val l\,|\,l\in\Fn(S)\}\le17$.
\endroster
This restriction is part of all statements: on the one hand, it helps us
eliminate a number of configurations covered by other theorems and, on the
other hand, it is sufficient for the proof of
\autoref{prop.56} in its current form.

Besides, we list all pairs $l_1,l_2\in\Fn(S)$ of skew lines such that
$\val l_1=\val l_2=16$ and compute the refined types of the pencils
$\pencil_i:=\pencil(l_i)$, $i=1,2$, and the linking types $\lk(l_1,l_2)$.
For each
pair $l_1,l_2$,
assuming that $\mult l_1\le\mult l_2$,
we require that
\roster[\lastitem]
\item\label{3-7.>=p}
$\mult l_1\ge p$, and
\item\label{3-7.=p}
if $\mult l_1=p$, then
$\ls|\pencil_1\cap\pencil_2|\ge\num2(s_0)+p$.
\endroster
(If these two conditions are not satisfied, we can obtain the same
configuration~$S$ replacing $l,s_0$ with the ``smaller'' pair $l_1,l_2$.)

In \autoref{obs.group}, approach~\iref{orbits.partial} is used for
sorting. In \autoref{obs.rank}, we may need to add up to two extra sections;
since the combinatorial rigidity is not known,
the configurations containing extra sections are excluded from the
sorting algorithm.
Finally, at the final step
we only keep the configurations~$S$ such that $\ls|\Fn(S)|>48$ or
$\ls|\Fn(S)|=48$ and $S$ is totally reflexive.

\lemma\label{lem.3-7}
Let $(S,\pencil)$ be a $(3,7)$-pair and $s_0\in S(\pencil)$ a section such
that
\[*
\max\{\val l\,|\,l\in\Fn(S)\}\le17\quad\text{and}\quad\val s_0=16.
\]
Then either one has $\ls|\Fn(S)|\le48$ or
$S\cong\config{56.real}$ or $\config{q56}$.
If $S$ is totally reflexive, then either
$\ls|\Fn(S)|<48$ or
$S\cong\config{56.real}$ or \config{48}.
\endlemma

\proof
The computation runs as outlined above.
In addition to \iref{3-7.val}--\iref{3-7.=p},
we inhibit all configurations in which $\pencil$ has a
section~$s$ such that $\num2(s)=6$, see \autoref{prop.3-7}.
We obtain several
configurations of type \config{56.real}, \config{q56}, or \config{48};
furthermore,
\roster*
\item
if $\num2(s_0)=0$, there is a single configuration~$S$;
this configuration~$S$ is of type \config{q56}, and
the pencils $\pencil$ and $\pencil'$ are of type $(3,7)\prim$;
\item
if $\num2(s_0)=1$,
there is a unique configuration~$S$ of type~\config{56.real} in which
$\pencil$ is of type $(3,7)\impr$ and $\pencil'$ is
of type $(4,4)$;
\item
if $\num2(s_0)=2$,
there is a unique configuration~$S$ of type~\config{48} in which
$\pencil$ is of type $(3,7)\impr$ and $\pencil'$ is of type~$(3,7)$.
\endroster
On the other hand, a direct computation shows that each configuration~$S$
obtained has a pair $\pencil$, $\pencil'$ of obverse pencils whose types and
intersection $\ls|\pencil\cap\pencil'|$ are as above. (Recall that
$\ls|\pencil(l)\cap\pencil(s_0)|=\num2(s_0)+3$.)
Replacing $l$ and~$s_0$ with the axes of these pencils, we conclude that,
up to isomorphism
{\em and under the assumptions of the lemma}, each
type \config{56.real}, \config{q56}, \config{48}
is represented
by a unique configuration.
\endproof

\lemma\label{lem.4-4}
Let $(S,\pencil)$ be a $(4,4)$-pair and $s_0\in S(\pencil)$ a section such
that
\[*
\max\{\val l\,|\,l\in\Fn(S)\}\le17\quad\text{and}\quad\val s_0=16.
\]
Then either one has $\ls|\Fn(S)|\le48$ or
$S\cong\config{56.real}$, \config{q56}, \config{q54},
\config{q52.1}, \config{q52.2}, or \config{50.0}.
If $S$ is totally reflexive, then either
$\ls|\Fn(S)|<48$ or
$S\cong\config{56.real}$ or \config{48.2}.
\endlemma

\proof
The configurations of type (hence, isomorphic to) \config{56.real} or
\config{q56} are given by \autoref{lem.3-7}.
The other types are obtained by a computation outlined above,
which returns several dozens of configurations with $\num2(s_0)\le2$.
Switching to approach~\iref{orbits.full} in \autoref{obs.group} and the full
automorphism group $\G_{4,4}$ {\em and checking explicitly that each
isomorphism $\barS'\to\barS''$ lifts to an isometry $S'\to S''$}, one can
show that, for any two configurations $S'$, $S''$ in the lists obtained,
$S'\cong S''$
if only if $\Ps(S')=\Ps(S'')$.
The pencil structures realized are those listed in the statement.
\endproof

\lemma\label{lem.5-1}
Let $(S,\pencil)$ be a $(5,1)$-pair and $s_0\in S(\pencil)$ a section such
that
\[*
\max\{\val l\,|\,l\in\Fn(S)\}\le17\quad\text{and}\quad\val s_0=16.
\]
Then one has $\ls|\Fn(S)|\le48$.
\endlemma

\proof
The computation runs as outlined at the beginning of this section, with
\autoref{lem.5-q.maximal} used to rule out some non-maximal pencils. This
computation results in an empty list of configurations.
\endproof

\subsection{Triangle free configurations}\label{s.trig.free}
A configuration~$S$ is called \emph{triangle} (respectively,
\emph{quadrangle}) \emph{free} if the graph $\Fn S$ has no cycles of length~$3$
(respectively, $3$ or~$4$). By \autoref{two-planes}, a configuration is
triangle free if and only if it contains no planes. Clearly, all pencils in
such a configurations are of type $(0,\q)$.

\lemma\label{lem.0-q}
Let $\pencil,\pencil'$ be a pair of obverse pencils in a
geometric triangle free
configuration~$S$, and assume that $\ls|\pencil\cap\pencil'|\ge2$. Then one
has
either $\ls|\pencil\cup\pencil'|\le18$ or $\ls|\Fn(S)|\le33$.
\endlemma

\proof
Assuming that $\ls|\pencil|\ge\ls|\pencil'|$, denote by~$s_0$ the axis
of~$\pencil'$; it is a section of~$\pencil$ and $r:=\num1(s_0)\ge2$.
Clearly, $\pencil$ is of type $(0,q)$, and we can assume that $q\ge11$
and $r\le2q-19$, as
otherwise the inequality $\ls|\pencil\cup\pencil'|\le18$ holds immediately.
The structure of the extension $\tP\supset P$ is given by
\autoref{prop.2.pivot} (the pivot has no $3$-torsion by
\autoref{prop.3.pivot}) and, depending on the values of $q$, $r$, there are
up to two (up to automorphism) possibilities for the section $s_0$.

We apply the algorithm outlined in \autoref{s.counting},
using the parameters
\[*
\pmin=\pmax=0,\quad \vmin=19-q,\quad\vmax=q-r
\]
and introducing a few modifications. Namely, at Step~1 we
allow repetition when collecting sections $\bs_i$, as the coordinate map
(\cf. \autoref{cor.coord}) is not injective for~$\pencil$; on the other hand,
only the sections satisfying \autoref{prop.2.pivot}\iref{2.pivot.section} are to
be considered.
At Step~2, as condition~\iref{post.maximal} in \autoref{obs.post},
we check that the configuration is still triangle free.
Adding, if necessary, up to two extra sections disjoint from~$s_0$ (see
\autoref{obs.rank}; such records are not sorted), we arrive at a number of
configurations, each containing at most $33$ lines.
\endproof

\lemma\label{lem.quad.3}
Let $S$ be a
geometric quadrangle free configuration. Consider three lines
$l_0\in\Fn(S)$ and $l_1,l_2\in\pencil(l_0)$
such that $\val l_0\ge\val l_1\ge\val l_2$.
Then either
\roster*
\item
$\val l_0+\val l_1\le14$ and $\val l_2=1$, or
\item
$\val l_0\le7$ and $\val l_2\le\val l_1\le5$, or
\item
$\val l_1\le\val l_0\le6$ and $\val l_2\le 5$.
\endroster
\endlemma

\proof
It is convenient to consider the pencil $\pencil:=\pencil(l_1)$, of which
$l_0$ is a fiber and $l_2$ is a section. Since $S$ is quadrangle free, each
section of~$\pencil$ intersects at most one fiber, and two sections
intersecting~$l_2$ cannot intersect the same fiber. In addition to~$l_2$,
the pencil~$\pencil$ has
$(\val l_0-2)$ sections intersecting~$l_0$ (all disjoint from~$l_2$) and
$(\val l_2-1)$ sections intersecting~$l_2$
(all disjoint from~$l_0$); all these sections are pairwise disjoint.
An extra parameter is the number of
the sections intersecting~$l_2$
that also intersect a fiber of~$\pencil$. A direct computation
(applying
\autoref{th.Nikulin} to the finite index extensions allowed by
\autoref{pre-conf}) rules out the values
\[*
(6, 6, 6),\ (7, 6, 1),\ (8, 5, 1),\ (10, 4, 2),\ (11, 3, 2),\ (11, 4, 1),\ (12, 2, 1)
\]
for the triple $(\val l_0,\val l_1,\val l_2)$.
\endproof

\subsection{Existence and uniqueness}\label{s.uniqueness}
We conclude this section with two
statements related to the uniqueness of large configurations and
the existence of
configurations with a prescribed number of lines.

\table
\caption{$\L$-configurations with more than $52$ lines
(see \autoref{lem.unique})}\label{tab.large}
\def\*{\raise3pt\hbox{$\scriptstyle\sqrt{}$}}
\hbox to\hsize{\hss\vbox{\halign{\strut\quad#\quad\hss&\hss#\hss\quad&
\hss#\hss\quad&\hss#\hss\quad&\hss#\hss\quad&&\hss$#$\hss\quad\cr
\noalign{\hrule}%
 \vphantom{\Big(}\hss$S$\hss&$\ls|\Fn|$&t.r.&ref&sym&\ls|\OG_h(S)|&\discr S&T:=S^\perp
\cr
\noalign{\vspace{-1pt}\hrule\vspace{2pt}}
\config{64}&64&&\*&\*&4608&
\CV_4\oplus\<\frac43\>&[8,4,8]
\cr
\config{60}&60&&\*&\*&480&
\CU_2\oplus\<\frac43\>\oplus\<\frac25\>&[4,2,16]
\cr
\config{60.2}&60&&&&240&
\<\frac65\>\oplus\<\frac{10}{11}\>&[4,1,14]
\cr
\config{56}&56&&&&128&%
\<\frac{15}8\>\oplus\<\frac{15}8\>&[8,0,8]
\cr
\config{56.real}&56&\*&\*&\*&64&%
\<\frac32\>\oplus\<\frac{63}{32}\>&[2,0,32]
\cr
\config{q56}&56&&\*&\*&384&
\CU_2\oplus\<\frac43\>\oplus\<\frac25\>&[4,2,16]
\cr
\config{54}&54&&\*&\*&384&
\<\frac14\>\oplus\<\frac38\>\oplus\<\frac43\>&[4,0,24]
\cr
\config{q54}&54&&\*&\*&48&
\CV_2\oplus\<\frac2{19}\>&[4,2,20]
\cr
\noalign{\vspace{2pt}\hrule}\crcr}}\hss}
\endtable

\lemma\label{lem.unique}
Each pencil structure listed in \autoref{tab.large} is realized by a unique,
up to weak isomorphism, $\L$-configuration~$S$. This $\L$-configuration~$S$ is
totally reflexive if and only if $S=\config{56.real}$\rom;
it is reflexive
unless $S=\config{60.2}$ or \config{56},
whereas \config{60.2} and \config{56} are not
symmetric.
\endlemma

\proof
By \autoref{th.uniqueness},
each pencil structure as in the
statement is realized by a unique geometric configuration~$S$; hence, there
only remains to verify that each of the three configurations
admits a unique primitive
$\L$-realization.

All configurations are known explicitly, and one can compute their
automorphism groups,
discriminants, and perspective transcendental lattices
$T:=S^\perp$;
they are as shown in \autoref{tab.large}.
(The lattice~$T$ is generated by two vectors $u$, $v$ so that
$u^2=a$, $u\cdot v=b$, and $v^2=c$, where $[a,b,c]$ is the triple given
in the table.
Each lattice is unique in its genus, which follows from the
classical theory of binary forms~\cite{Gauss:Disquisitiones}.)

With two exceptions, the homomorphism $\rho\:\OG_h(S)\to\Aut\discr S$ is
surjective. The exceptions are:
\roster*
\item
$S=\config{q54}$, where $\Im\rho=\Aut\discr_2S$,
and
\item
$S=\config{56}$, which will be treated separately.
\endroster
Furthermore, each involution in $\Im\rho$ lifts to an involution in
$\OG_h(S)$. (This is not a common property of
configurations, a counterexample being \config{50.0}, see \autoref{ss.Z}.)

In each case (other than $S=\config{56}$), it is immediate that
the image of $\OG^+(T)$ intersects each coset modulo $\Im\rho$ and,
hence, a primitive $\L$-realization is unique
up to weak isomorphism (see the description of primitive extensions
in \autoref{s.lattices}).
Besides, whenever $T$ has an orientation reversing isometry (\ie, in all
cases except $S=\config{60.2}$, see \autoref{tab.large}), this
isometry, which is necessarily involutive, can be chosen to induce an element in
the image $\Im\rho$ and, thus, lift to an involution in $\OG_h(S)$.
Hence, the $\L$-configuration is
symmetric and reflexive.

In the exceptional case $S=\config{56}$, the image of $\OG_h(S)$
(respectively, $\OG(T)$) is the index~$2$
subgroup of $\Aut\discr S$ generated by the reflections $t_\Ga$,
where $\Ga\in\discr S$ and $\Ga^2=\frac38$ or
$\frac34\bmod2\Z$ (respectively, $\Ga^2=\frac{15}8$ or $\frac74\bmod2\Z$).
The intersection of the two subgroups has index~$4$ and coincides with the
image of $\OG^+(T)$. It follows that there is a single weak isomorphism
class, which is not symmetric.

The only
totally reflexive configuration is \config{56.real},
as $\config{56.real}^\perp$ is the only transcendental lattice
containing a
vector of square~$2$, see \autoref{tab.large}.
\endproof

\remark\label{rem.aut}
The
computation of the automorphism groups make use of the pencil structure:
we list all pencils of a given type (usually, the first one listed in
\autoref{tab.list}) and then enumerate the isometries taking one fixed
pencil to another one similar to the sorting algorithm in
\autoref{obs.group}.
\endremark

\remark\label{rem.nonunique}
Not every configuration~$S$ listed in \autoref{th.uniqueness} admits a unique
$\L$-realizations.
Simplest examples are $\config{Y}_*$, see
\autoref{ss.real} and \autoref{tab.Y} below.
More examples are found in \autoref{tab.Y} in \autoref{S.examples}.
\endremark

\lemma\label{lem.counts}
For any
number
$n\in\{0,\ldots,52,54,56,60,64\}$, there is a geometric configuration~$S$
such that $\ls|\Fn(S)|=n$.
If
$n\in\{0,\ldots,48,52,56\}$, this configuration can be chosen totally reflexive.
\endlemma

\proof
Any count $n\le17$ is easily realized by the span of a single pencil.
Hence, the first statement follows from \autoref{th.6-q},
and the second one
mostly follows from \autoref{cor.4-5.counts}.
The missing values $n=48$, $56$
for totally reflexive
configurations
are
given by \autoref{lem.3-7}.
\endproof

\section{Triangle free configurations}\label{S.triang.free}

Recall
that a configuration~$S$ is said to be triangle free if it contains no
planes. The principal goal of this section is a proof of a bound to the
number of lines in such a configuration, see \autoref{th.Segre}
in \autoref{s.triang.free} below.
Throughout the section, we fix a configuration~$S$ and
a geometric $\L$-realization $\psi\:S\to\L$.

\subsection{Adjacency graphs}
Given a graph~$\Gamma$, we denote by $\Z\Gamma$ the lattice freely generated
by the vertices of $\Gamma$, so that $v^2=-2$ for each vertex~$v$
and $u\cdot v=1$ (respectively,~$0$) if the vertices $u\ne v$ are
(respectively, are not) adjacent in $\Gamma$. If $\Gamma\subset\Fn(S)$, we also
consider the images $S\Gamma:=\Z\Gamma/\ker\subset S$ and
$\L\Gamma:=\psi(S\Gamma)\subset\L$
of this lattice in~$S$ and~$\L$,
denoting by $\psi_\Gamma\:\Z\Gamma\to\L$ the composed map.

A graph $\Gamma$ is called \emph{elliptic} (respectively,
\emph{parabolic}) if $\Z\Gamma$ is negative
definite (respectively, negative semi-definite).
The \emph{Milnor number} $\mu(\Gamma)$ of an elliptic or parabolic
graph~$\Gamma$ is the rank of the lattice $\Z\Gamma/\ker$.
A connected elliptic (parabolic) graph is called a \emph{Dynkin diagram}
(respectively, \emph{affine Dynkin diagram}).
A Dynkin diagram~$D$ extends to a unique affine Dynkin diagram,
which we denote
by $\aD\supset D$;
we refer to \cite{Bourbaki:Lie:French} for a detailed treaty of Dynkin
diagrams and their affine counterparts.
Recall that any graph~$\Gamma$ such that $\Z\Gamma$ is not negative definite
contains an affine Dynkin diagram as an induced subgraph.
For
any affine
Dynkin diagram~$\aD$, the kernel  $\ker\Z\aD$
is spanned by a single distinguished generator
$k_{\aD}=\sum\kappa(e)e$, $e\in\aD$, with each coefficient
$\kappa(e)$ \emph{strictly positive}.
The coefficient sum
\smash{$\kappa(\aD):=\sum\kappa(e)$} of this linear combination is as follows:
\[
\kappa(\tA_p)=p+1,\quad
\kappa(\tD_q)=2q-2,\quad
\kappa(\tE_6)=12,\quad
\kappa(\tE_7)=18,\quad
\kappa(\tE_8)=30.
\label{eq.kappa}
\]
We extend this $\kappa$-notation to elliptic Dynkin diagrams letting
$\kappa(D):=\kappa(\aD)$.

\lemma\label{lem.Dynkin.mono}
Let $\Gamma\subset\Fn(S)$ be a parabolic subgraph such that
$\rank\ker\Z\Gamma=1$.
Then, the isometry $\psi_\Gamma\:\Z\Gamma\to\L$ is a monomorphism.
\endlemma

\proof
By the assumption, $\Gamma$ is a disjoint union of several Dynkin diagram and a
single affine Dynkin diagram~$\aD$.
Since $\psi_\Gamma$ is an isometry, one has $\Ker\psi_\Gamma\subset\ker\Z\Gamma$,
and, as
explained above, the latter subgroup is spanned by a single vector $k_{\aD}$ so
that \smash{$\psi_\Gamma(k_{\aD})\cdot h=\kappa(\aD)>0$}.
Hence, $\psi_\Gamma(k_{\aD})\ne0$ and $\Ker\psi_\Gamma=0$.
\endproof

\subsection{Pseudo-pencils}
Given a nonzero isotropic vector $v\in S$,
the \emph{pseudo-pencil} defined by $v$ is the set
\[*
\pseudo(v):=\bigl\{a\in\Fn(S)\bigm|a\cdot v=0\bigr\}.
\]
Since $S$ is hyperbolic,
$v\cdot h\ne0$ and we can assume
$v\cdot h>0$. We can also assume~$v$ primitive. Then, the
integer $\deg\pseudo:=v\cdot h$ is called the \emph{degree}
of~$\pseudo$.
The connected components of~$\pseudo$ are called its \emph{fibers}.
A \emph{section} (more generally, \emph{$n$-section}, $n>0$)
of~$\pseudo$ is a line
$s\in\Fn(s)$ such that $s\cdot v=1$ (respectively, $s\cdot v=n$).
The set of sections of~$\pseudo$, depending on the ambient
configuration~$S$, is denoted by $S(\pseudo)$.

Each pencil is a pseudo-pencil of degree~$3$:
one has $\pencil(l)=\pseudo(h-l)$.
Conversely, if $v\cdot h=3$, then $l:=h-v\in\Fn(s)$ and
$\pseudo(v)=\pencil(l)$.

As another example, fix an affine Dynkin diagram $\aD\subset\Fn(S)$ and
let $v\in S$ be the image of $k_{\aD}$; by \autoref{lem.Dynkin.mono},
$\psi(v)\ne0$ and
$\pseudo(\aD):=\pseudo(v)$ is a pseudo-pencil.
Clearly, \smash{$\aD\subset\pseudo(\aD)$}.
Since $k_{\aD}$ is a \emph{positive} linear combinations of the vertices
of~$\aD$ and the intersection of two lines is nonnegative (see
\autoref{matrix}), it follows that
\[
\pseudo(\aD)=\{a\in\Fn(S)\,|\,\text{$a\cdot v=0$ for each vertex $v\in\aD$}\}.
\label{eq.pseudo.count}
\]

\proposition\label{prop.pseudo}
For each pseudo-pencil~$\pseudo$ the following statements hold\rom:
\roster
\item\label{i.pseudo.1}
either $\deg\pseudo=1$ and $\ls|\Fn(S)|=1$, or $\deg\pseudo\ge3$\rom;
\item\label{i.pseudo.2}
as a graph, $\pseudo$ is elliptic or parabolic
and
$\mu(\pseudo)\le18$\rom;
\item\label{i.pseudo.3}
if $D\subset\pseudo$ is a Dynkin diagram and
$(\deg\pseudo)\mathbin|\kappa(D)$, then
$\aD$
is a fiber of $\pseudo$.
\endroster
Furthermore, if $s\in S(\pseudo)$, then, for any parabolic fiber $\aD$
of~$\pseudo$,
\roster[4]
\item\label{i.pseudo.4}
$\sum\kappa(e)(s\cdot e)=\deg\pseudo$, the summation running over
$e\in\aD$\rom;
\item\label{i.pseudo.5}
in particular, if $S(\pseudo)\ne\varnothing$, then
$\kappa(\aD)=\deg\pencil$ and $k_{\aD}=v$.
\endroster
\endproposition

\proof
Let $\pseudo=\pseudo(v)$ with $v\cdot h=\deg\pseudo$.
The possibility $v\cdot h=2$ is excluded by
\autoref{ell-pencil} in \autoref{pre-conf}.
If $v\cdot h=1$, then $a:=h-3v$ is a line.
Consider another line $b\in\Fn(S)$.
If $b\cdot v\ne0$ or~$1$, then $\Gs_+(\Z h+\Z v+\Z b)=2$.
If $b\cdot v=0$, then $e:=b-v$
is as in \autoref{ex-div} in \autoref{pre-conf}.
In the remaining case $b\cdot v=1$ one has $\rank\ker(\Z h+\Z v+\Z b)=2$ and,
hence, $b=a$, \ie, $a$ is the only line.

The assumption that $v\ne0$ implies that $v^\perp$ has
a non-trivial kernel and, hence, is parabolic;
since also $\rank\psi(\pseudo)\le19=\Gs_-\L$,
this proves \autoref{i.pseudo.2}.

For \autoref{i.pseudo.3}, observe that $\kappa(e_0)=1$
for the only vertex $e_0\in\tilde D\sminus D$, see,
\eg,~\cite{Bourbaki:Lie:French}.
Hence, $e_0$ is an integral linear combination of~$v$ and the
vertices of~$D$, \ie, $e_0\in S$. Clearly, $e_0$ is a
line and, thus, $\tilde D\subset\pseudo$.
Finally,
any affine
Dynkin diagram is a whole connected component of any parabolic graph in which it is
contained.

The last two statements follow
from the definitions and the fact
that, for each parabolic fiber~$\aD$ of~$\pseudo$, the vector $k_{\aD}$ is
a multiple of~$v$ (as $k_{\aD}\cdot v=0$);
on the other hand, $\sum_{e\in\aD}\kappa(e)(s\cdot e)=s\cdot k_{\aD}$.
\endproof

\corollary\label{cor.pseudo}
For a pseudo-pencil $\pseudo$, one has $\ls|\pseudo|\le18(1+1/\mu)$,
where $\mu$ is the
minimal Milnor number
of the parabolic fibers of~$\pseudo$.
In particular,
$\ls|\pseudo|\le24$.
\endcorollary

\proof
The first bound follows from the obvious identity
\[*
\ls|\pseudo|=\mu(\pseudo)+\ls|\{\text{parabolic fibers of $\pseudo$}\}|.
\]
If $\pseudo$ has a fiber of type~$\tA_2$, it is an ordinary pencil and
$\ls|\pseudo|\le20$ by \autoref{cor.le20}. Otherwise, $\mu\ge3$ and we have
$\ls|\pseudo|\le24$.
\endproof

Geometrically, if $S=\Fano(X)$ for a nonsingular quartic $X\subset\Cp3$,
a pseudo-pencil~$\pseudo$ can often be interpreted as an elliptic pencil
$\pi\:X\to\Cp1$ whose
fibers are curves of degree $\deg\pseudo$ in~$\Cp3$.
For example, this is so in the important special case where
$\pseudo$ has a parabolic fiber~$\aD$.
Indeed, in this case,
the class $v=\sum\kappa(e)e$, $e\in\aD$, regarded as a divisor,
is obviously
numerically effective
and, hence, does
define a linear system of
arithmetic genus~$1$ without fixed points or components.
From this
geometric point of view, $\pseudo$ is the union of lines contained in the
fibers of~$\pi$. More precisely, if \emph{all} components of a reducible fiber~$F$
of~$\pi$ are lines, these lines form a parabolic fiber of~$\pseudo$;
otherwise, the lines contained in~$F$ constitute
one or
several elliptic fibers of~$\pseudo$.
Furthermore, in this interpretation, the bound $\ls|\pseudo|\le24$
given by \autoref{cor.pseudo} follows from the inequality
\[*
\ls|\{\text{components in the singular fibers of~$F$}\}|\le\chi(X)=24.
\]
Using this geometric interpretation, one can
partially extend
Statements~\iref{i.pseudo.4} and~\iref{i.pseudo.5} of \autoref{prop.pseudo}
to the elliptic fibers of~$\pseudo$. Namely, for each
section $s\in S(\pseudo)$ and each elliptic fiber~$D$ of~$\pseudo$, one has
\roster[4]
\item
$\sum\kappa(e)(s\cdot e)\le\deg\pseudo$, the summation running over
$e\in D$\rom;
\item
in particular, if $S(\pseudo)\ne\varnothing$, then
$\kappa(D)<\deg\pencil$.
\endroster
As we do not use these statements, we will not try to prove them
arithmetically.
(Unlike \autoref{prop.pseudo}, these statements may depend on the requirement
that $S$ should be geometric and involve a case-by-case analysis, \cf. the
discussion below.)

The \emph{type} of a pseudo-pencil~$\pseudo$ is the isomorphism type of the
lattice $\Z\pseudo$; by \autoref{prop.pseudo}, it is an orthogonal
direct sum of elliptic and parabolic root lattices. (For example, in this new
language, an ordinary pencil of type $(p,q)$ has type
\smash{$p\tA_2\oplus q\bA_1$}.)
Using \autoref{prop.pseudo} and arguing as in \autoref{S.arithm}, \ie, applying
Nikulin's \autoref{th.Nikulin} to all finite index extensions of
the lattice $P:=(\Z\pseudo+\Z h)/\ker$ that are not ruled out
by \autoref{pre-conf},
it should not be
difficult to obtain a complete classification of pseudo-pencils appearing in
geometric configurations; in particular, one can probably improve the bound
$\ls|\pseudo|\le24$ given by \autoref{cor.pseudo}. However, we confine
ourselves to just the two special cases used in the proof of
\autoref{th.Segre}.

\lemma\label{lem.tA3}
Assume that $S$ is triangle free, and
let $\pseudo\subset S$ be a pseudo-pencil with a fiber of type~$\tA_3$. Then
either $\ls|\pseudo|\le20$ or $\pseudo$ is of type $5\tA_3\oplus\bA_1$\rom;
in the latter case, one has $\ls|\Fn(S)|\le45$.
\endlemma

\proof
By \autoref{prop.pseudo}, one has $\deg\pseudo=\kappa(\tA_3)=4$ and all
fibers of~$\pseudo$ are of types \smash{$\tA_3$}, $\bA_2$, or $\bA_1$.
Arguing as explained above, we conclude that the only pseudo-pencil~$\pseudo$
such that $\ls|\pseudo|>20$ and
the lattice $P:=(\Z\pseudo+\Z h)/\ker$ admits a geometric
$\L$-realization is that of type \smash{$5\tA_3\oplus\bA_1$}.
Assuming this type, consider the quadrangle
$\aD:=\{l_1,\ldots,l_4\}$ constituting one of
the type~$\tA_3$ fibers. Letting $\pencil_i:=\pencil(l_i)$,
by~\eqref{eq.pseudo.count} we have
\[
\ls|\Fn(S)|=\ls|\pencil_1\cup\pencil_3|
 +\ls|\pencil_2\cup\pencil_4|+\ls|\pseudo|-4.
\label{eq.tA3}
\]
(Since $S$ is triangle free, a line $a\in\Fn(S)$ cannot intersect two
adjacent vertices of the quadrangle.)
Due to \autoref{lem.skew}\iref{skew.2} and \autoref{cor.10.fibers}, for each
of the two pairs $(i,j)=(1,3)$ or $(2,4)$, either
$\ls|\pencil_i\cup\pencil_j|=\ls|\pencil_i\cap\pencil_j|=10$
or $\ls|\pencil_i\cap\pencil_j|\le8$; thus, letting
$n_i:=\ls|\pencil_i\sminus\pencil_j|$, we get
$\ls|\pencil_i\cup\pencil_j|\le\max\{20,16+n_i+n_j\}$
and, if $n_i\le3$ for all $i=1,\ldots,4$, from~\eqref{eq.tA3} we obtain
$\ls|\Fn(S)|\le45$, as stated.

What remains is the case where one of the integers~$n_i$, say, $n_1$, is at
least~$4$, \ie, there are at least four lines intersecting~$l_1$ and disjoing
from the three other lines.
In this case, we run an algorithm similar to that described in
\autoref{s.counting}, adding to~$S$ up to three sections
intersecting~$l_1$ in order
to increase the rank from $\rank P=18$ to the maximum~$20$.
By \autoref{prop.pseudo}\iref{i.pseudo.4}, each section
intersects exactly one line of each other parabolic
fiber; given the rich automorphism group, this observation leaves
relatively few possibilities for pairs and triples of sections.
Then, as in \autoref{obs.L}, we enumerate the geometric realizations of each
configuration of maximal rank and compute the number of lines, arriving at
the inequality $\ls|\Fn(S)|\le33$.
\endproof

\lemma\label{lem.tD4}
If $\pseudo\subset S$ is a pseudo-pencil with a fiber of type~$\tD_4$,
then $\ls|\pseudo|\le19$.
\endlemma

\proof
By \autoref{prop.pseudo}, one has $\deg\pseudo=\kappa(\tD_4)=6$ and all
fibers of~$\pseudo$ are of types \smash{$\tD_4$}, \smash{$\tA_5$},
or $\bA_p$, $1\le p\le4$. Trying all combinations one by one and arguing as
explained prior to \autoref{lem.tA3},
we arrive at the inequality stated. (In fact, the only type with
$\ls|\pseudo|=19$ is \smash{$2\tD_4\oplus\tA_5\oplus\bA_2\oplus\bA_1$}.)
\endproof

\subsection{The bound}\label{s.triang.free}
The following theorem is the principal result of this section.

\theorem\label{th.Segre}
If a geometric configuration $S$ is triangle free, then
$\ls|\Fn(S)| \leq 52$.
\endtheorem

\proof
We consider separately several cases, each time picking an appropriate affine
Dynkin diagram $\aD\subset\Fn(S)$ and using~\eqref{eq.pseudo.count} to
estimate the number of lines,
which is
$\ls|\pseudo(\aD)|+\ls|\{\text{lines intersecting a vertex of~$\aD$}\}|$.

First, assume that the maximal valency of a line in~$S$ is at most~$3$.
If $\Fn(S)$ is elliptic, then $\ls|\Fn(S)|\le19$. Otherwise, $\Fn(S)$
contains an affine Dynkin diagram; pick one \smash{$\aD\subset\Fn(S)$}
of the minimal
Milnor number~$\mu$. Using the classification of affine Dynkin diagrams,
we conclude that the number of lines that are not in~$\aD$ and adjacent to a
vertex of $\aD$ is at most $2n_1+n_2\le\mu+3$, where $n_i$ is the number of
vertices of~\smash{$\aD$} of valency~$i$. Since $2\le\mu\le18$,
by~\eqref{eq.pseudo.count} and \autoref{cor.pseudo},
\[*
\ls|\Fn(S)|\le\mu+3+\ls|\pseudo(D)|\le\mu+\frac{18}\mu+21\le40.
\]

Now, assume that $S$ has a line of valency at least~$4$ and is quadrangle
free. Let $l_0$ be a line of maximal valency, and pick four lines
$l_1,\ldots,l_4$ adjacent to~$l_0$ so that $\val l_1\ge\ldots\ge\val l_4$.
Then, $\aD:=\{l_0,\ldots,l_4\}$ is a subgraph of type~$\tD_4$ and,
by~\eqref{eq.pseudo.count} and \autoref{lem.tD4},
\[*
\ls|\Fn(S)|\le\val l_0+\val l_1+\val l_2+\val l_3+\val l_4+11.
\]
The sum of the valencies in the latter expression is estimated using
\autoref{lem.quad.3} (and the assumption $\val l_3,\val l_4\le\val l_2$),
and we obtain $\ls|\Fn(S)|\le38$.

Finally, assume that $\Fn(S)$ has a quadrangle, \ie, a $4$-cycle
$l_1,l_2,l_3,l_4$, which can be regarded as a subgraph~\smash{$\aD$} of
type~\smash{$\tA_3$}. Assume that $\ls|\Fn(S)|\ge46$ and
apply \eqref{eq.tA3}:
each of the first two terms is bounded by~$18$ by \autoref{lem.0-q},
and $\ls|\pseudo|\le20$ by \autoref{lem.tA3}; hence,
$\ls|\Fn(S)|\le52$.
\endproof

\remark\label{rem.triang.free}
The idea that triangle free configurations of lines in quartics should be
treated separately is also due to B.~Segre, and his geometric proof~\cite{Segre}
of the bound $\ls|\Fn(S)|\le64$ for such configurations can easily be
modified to get $\ls|\Fn(S)|\le60$. Our bound $\ls|\Fn(S)|\le52$ given by
\autoref{th.Segre} can be improved to $\ls|\Fn(S)|\le50$: in
\autoref{lem.tA3}, the few types with $\ls|\pseudo|=19$ or~$20$ can be ruled
out similar to \smash{$5\tA_3\oplus\bA_1$}. Probably, this better bound is
still not sharp:
currently, the best known example of triangle free configurations has
$37$ lines.
\endremark

\section{Proofs}\label{S.pencils}

In this section, we prove the principal results of the paper,
\viz. \autoref{th.unique}, \autoref{th.56}, and
\autoref{ad.counts}.

\subsection{Large configurations}\label{s.large.conf}
All proofs
are based on the following
statement,
which bounds the number of lines in a
geometric
configuration containing a
plane.
With further applications in mind, we state it in a slightly stronger form.

\proposition\label{prop.56}
If a geometric configuration~$S$ contains a plane, then either
\roster*
\item
$S$ is isomorphic to \config{64},
\config{60}, \config{60.2}, \config{56}, \config{56.real}, \config{q56},
\config{54}, \config{q54},
\config{52.1}, \config{52.2}, \config{52.4}, \config{52.5-3}, \config{52.5},
or~\config{52.0},
or
\item
one has $\ls|\Fn(S)|\le52$ and $\max\{\val l\,|\,l\in\Fn(S)\}\le17$, or
\item
one has $\ls|\Fn(S)|<52$.
\endroster
\endproposition

\proof
Assume that $\ls|\Fn(S)|\ge52$.
If $S$
has a pencil of type $(6,\q)\impr$,
\autoref{th.6-q} implies that
$S\cong\config{64}$, \config{60}, \config{54}, \config{52.1}, \config{52.2},
or \config{52.0}.
Hence, from now on we can also assume that $S$ does not
have such a pencil. In particular,
\[*
v:=\max\{\val l\,|\,l\in\Fn(S)\}\le18;
\]
if $v\le15$, then $\ls|\Fn(S)|=52$
by \autoref{cor-Segre}.

Pick a maximal pencil~$\pencil$ such that $\ls|\pencil|=v$.
By~\eqref{tablichka}, this pencil~$\pencil$ has
a $3$-fiber $\{m_1,m_2,m_3\}$, which we order so that
$\val m_1\le\val m_2\le\val m_3$.
We have
\[*
\val m_1 + \val m_2 + \val m_3 = \ls| \Fn(S) | + 8 - v \ge 42;
\]
hence $\val m_3\ge14$. Then $\mult m_3\ge2$ by~\eqref{tablichka} again,
and one can find another plane
$\{s_0,s_1,s_2,m_3\}$ containing~$m_3$.
The lines $s_0,s_1,s_2$ are sections of~$\pencil$, and they satisfy
the inequality
\[*
\val s_0 + \val s_1 + \val s_2 = \ls|\Fn(S)|+8-\val m_3.
\]
Assuming that $\val s_0\ge\val s_1\ge\val s_2$, we obtain
\[
3\val s_0\ge\ls|\Fn(S)|+8-\val m_3.
\label{eq.valency}
\]

Let $v=18$. We need to show that $\val s_0\ge15$; then,
Theorems~\ref{th.6-0.primitive}, \ref{th.4-6}, and~\ref{th.5-3},
would imply that $S\cong\config{56}$,
\config{52.4}, \config{52.5-3}, or~\config{52.5}.
If $\val m_3\le17$, the desired inequality $\val s_0\ge15$ follows
from~\eqref{eq.valency}. If $\val m_3=18$ \emph{and} $\val s_0\le14$, we
repeat the same argument, taking $m_3$ and~$s_0$ for~$l$ and~$m_3$,
respectively, and obtaining a section $s_0'$ of the new pencil $\pencil(m_3)$ of
valency $\val s_0'\ge16$.

If $v=16$ and $\ls|\Fn(S)|>52$, the same argument as above produces a
pencil~$\pencil'$
(not necessarily the original one)
and section $s_0'$ of~$\pencil'$ such that
$\ls|\pencil'|=\val s_0'=16$; hence, Lemmas~\ref{lem.3-7}, \ref{lem.4-4},
and~\ref{lem.5-1} imply that $S\cong\config{56.real}$, \config{q56}, or
\config{q54}.

Finally, let $v=17$. If $\ls|\Fn(S)|\ge54$, we use the same argument to get a
pencil~$\pencil'$ and section~$s_0'$ of~$\pencil'$ such that $\ls|\pencil'|=17$
and $\val s_0'\ge16$; hence, by Lemmas~\ref{lem.4-5} and~\ref{lem.5-2}, we have
$S\cong\config{60.2}$. If $\ls|\Fn(S)|=53$, the argument may fail as one may
have $\val s_0\le15$ and $\val m_3=16$.
But in the latter case, starting with $\pencil':=\pencil(m_3)$, we obtain a
section $s_0'$ of $\pencil'$ such that $\val s_0'\ge16$;
this is a contradiction to Lemmas~\ref{lem.3-7}, \ref{lem.4-4},
and~\ref{lem.5-1} (if $\val s_0'=16$) or~\ref{lem.4-5} and~\ref{lem.5-2} (if
$\val s_0'=17$;
in this latter case,
when applying the lemmas, we regard~$m_3$ as a section of $\pencil(s_0')$).
\endproof

\subsection{Real configurations}\label{s.config.real}
In the next statement, we consider a configuration~$S$
equipped with a ``real structure'', \ie,
involutive automorphism $S\to S$, $a\mapsto\bar a$.
For such a configuration, the \emph{real part} is the subconfiguration
$S_\R:=\{a\in S\,|\,\bar a=a\}$.
We let $\Fn_\R(S):=\Fn(S_\R)$ and call the lines contained in $\Fn_\R(S)$
\emph{real}.

\proposition\label{prop.real.plane}
Let $S$ be a geometric configuration equipped with an involutive automorphism
$a\mapsto\bar a$, and assume that $\ls|\Fn_\R(S)|>48$. Then any plane
$\Ga\subset\Fn(S)$ is \emph{totally real}, \ie, $\Ga\subset\Fn_\R(S)$.
\endproposition

\proof
Consider a plane $\alpha = \{a_1, a_2, a_3, a_4\}$.
Let $r$ be the number of real lines in $\alpha$,
and let $r_i$ be the number of real lines in $\pencil(a_i)$, $i = 1,\ldots,4$.
The following formula is a
straightforward modification of the conclusion of \autoref{cor-Segre}:
\[*
\ls| \Fn_\R(X) | = r _1 + r_2 + r_3 + r_4 - 2r. \label{Segre-real-equality}
\]
If $a_i$ is real, then $r_i \le \ls|\pencil(a_i)|\le 20$ by~\eqref{tablichka}.
Otherwise,
$r_i \le \ls|\pencil(a_i) \cap \pencil(\bar a_i)|$, which
does not exceed $2$ or~$10$ if $a_i\cdot\bar a_i=1$ or~$0$,
respectively, see \autoref{lem.2-pencils}.

Consider the conjugate plane~$\bar\Ga$. If $\alpha \cap \bar\alpha = \varnothing$,
then
$r = 0$ and $\ls | \Fn_\R(X)| \le 40$.
If $\ls| \alpha \cap \bar\alpha | = 1$ (\ie, $r = 1$),
then $\ls| \Fn_\R(X)| \le 48$.
If $\ls| \alpha \cap \bar\alpha | > 1$,
then $\alpha = \bar\alpha$ by \autoref{two-planes}
and $r_i \le 2$ for each non-real line $a_i$;
hence, $\ls| \Fn_\R(X)| \le 16r + 8$ and, since $r \ne 3$,
we conclude that $r = 4$, \ie, $\Ga\subset\Fn_\R(S)$.
\endproof

The following corollary is a real counterpart of \autoref{th.Segre}.

\corollary\label{Segre_real}
Let $X \subset \Cp3$ be a nonsingular real quartic.
If $\ls| \Fn_\R(X) | > 52$,
then $\Fano_\R(X)$ contains a plane\rom;
moreover, any plane in $\Fano(X)$ is contained in $\Fano_\R(X)$.
\endcorollary

\proof
Clearly, $\Fano_\R(X)$ is the real part of the Fano configuration $\Fano(X)$
with respect to the involution $a\mapsto-\conj_*a$ induced by the real
structure. The configuration $\Fano(X)$ is geometric
(see \autoref{th.complex})
and it contains a plane (see \autoref{th.Segre});
there remains to apply \autoref{prop.real.plane}.
\endproof

\subsection{Proof of \autoref{th.unique}}\label{proof.64}
According to \autoref{th.complex},
the Fano configuration $\Fano(X)$ is
geometric and, since we assume $\ls|\Fn(X)|>52$,
\autoref{th.Segre} implies that this configuration
contains a plain. Then, by
\autoref{prop.56}, $\Fano(X)$ is isomorphic to
one of the configurations listed in \autoref{tab.large},
and the statement of the theorem
follows from \autoref{lem.unique} and \autoref{app.complex}.
(The quartic corresponding to~\config{64} is
identified as Schur's quartic since both contain $64$
lines.)
\qed

\subsection{Proof of \autoref{th.56}}\label{proof_th.56}
The real Fano configuration $\Fano_\R(X)$ is
geometric (see \autoref{th.complex})
and, assuming that $\ls|\Fn_\R(X)|>52$,
this configuration contains a plain due to \autoref{Segre_real}.
Then,
the statement of the corollary follows from \autoref{prop.56} and
\autoref{cor.real}.
\qed

\subsection{Proof of \autoref{ad.counts}}\label{proof.counts}
The statement is an immediate consequence
of \autoref{lem.counts} and \autoref{th.complex} (for
lines in complex quartics) or \autoref{cor.real} (for real lines in real
quartics).
\qed

\section{The known examples}\label{S.examples}

\subsection{Schur's quartic}\label{s.Schur}
The following example
is more than $130$ years old: it goes back to F.~Schur~\cite{Schur:quartics}
(see also \cite{barth:1983,boissiere:2007}).
According to our \autoref{th.unique}, this is the \emph{only} nonsingular
quartic containing $64$ lines, and its configuration of lines is \config{64}.

Consider the quartic~$\quartic{64}$ given by the equation
\[
\Gf(z_0,z_1)=\Gf(z_2,z_3),\qquad
 \Gf(u,v):=u(u^3-v^3).
\label{eq.Schur}
\]
Let $k_0:=0$, $k_1:=1$, and $k_{2,3}:=(-1\pm i\sqrt3)/2$ be the four roots of
$\Gf(u/v,1)$. Then, $\quartic{64}$ contains the sixteen lines
\[
z_1=k_rz_0,\quad z_3=k_sz_2,\quad r,s=0,\ldots,3.
\label{eq.16}
\]
Besides,
$\quartic{64}$ contains the line
\[*
l_0:=\{z_0=z_2,\ z_1=z_3\}.
\]
Finally, observe that $\Gf$ is the ``most symmetric'' polynomial of degree
four: its zero locus $\{k_0,k_1,k_2,k_3\}\subset\Cp1$ has $j$-invariant~$0$,
\ie, $\Gf$ is
invariant under a subgroup $G\cong\AG4\subset\PGL(2,\C)$.
This subgroup lifts a a subgroup $\bG\subset\GL(2,\C)$ preserving $\Gf$
literally, not just up to a factor; it is generated by
\[*
\frac1{\sqrt3}\bmatrix1&-1\\-2&-1\endbmatrix,
 \bmatrix1&0\\0&\epsilon\endbmatrix\in\GL(2,\C),\quad
 \epsilon^3=1,\epsilon\ne1,
\]
and the kernel of the projection $\bG\onto G$ is the
central subgroup $H\cong\CG4$ generated by $i\id$.
Letting~$\bG$ act separately on $(z_0,z_1)$ and $(z_2,z_3)$, we obtain a
subgroup $\Aut_0\quartic{64}:=\bG\odot\bG\subset\Aut \quartic{64}$, where the central
product is the quotient of $\bG\times\bG$ by the diagonal
$H\subset H\times H$.
The stabilizer of~$l_0$ is the diagonal $\bG/H\subset\Aut_0\quartic{64}$; hence,
its orbit consists of $48$ distinct lines, and $X$ contains $16+48=64$ lines.

A computation of the intersection matrix reveals that the sixteen
lines~\eqref{eq.16}
are distinguished: each is contained in six planes~$\Ga$ such that
$\quartic{64}\cap\Ga$ splits into four lines, whereas any other line
is contained in
four such planes. Hence, any \hbox{(anti-)}\penalty0automorphism of~$\quartic{64}$
preserves
the pair of lines $m_{ij}:=\{z_i=z_j=0\}$,
$(i,j)=(0,1)$ or $(2,3)$.
It follows that $\Aut \quartic{64}$ is an extension of
the group $\Aut_0\quartic{64}$ preserving each of $m_{01}$, $m_{23}$
by the involution $z_0\leftrightarrow z_2$, $z_1\leftrightarrow z_3$
interchanging $m_{01}\leftrightarrow m_{23}$.
This group has order~$1152$.
As a consequence, we have the following statement.

\proposition\label{prop.X64}
Up to automorphism, the quartic $\quartic{64}$ has four real structures,
\viz. those sending $[z_0:z_1:z_2:z_3]$ to
\[*
[\bar z_0:\bar z_1:\bar z_2:\bar z_3],\quad
[\bar z_0:\bar z_1:i\bar z_2:i\bar z_3],\quad
[\bar z_2:\bar z_3:\bar z_0:\bar z_1],\quad
[\bar z_2:\bar z_3:-\bar z_0:-\bar z_1].
\]
The numbers of real lines are $8$, $4$, $28$, and~$4$, respectively.
\endproposition

\proof
Denote by $\bar\,$ the standard complex conjugation, and extend its action to
matrices. Then, any real structure
on~$\quartic{64}$ is
$\Gs_g\:z\mapsto g\bar z$, where $g\in\Aut \quartic{64}$ is such that $g\bar g=\id$.
Two real structures $\Gs_g$, $\Gs_{g'}$ are isomorphic if and only if
one has $g'=h\1g\bar h$ for some $h\in\Aut \quartic{64}$.

The set of lines real with respect to a real structure
$\Gs_g$ is found as follows.
A line $l\subset \quartic{64}$ as in~\eqref{eq.16} is
uniquely
determined by its
``endpoints'' $l\cap m_{01}$, $l\cap m_{23}$, and the set of all eight
endpoints is preserved by any (anti-)automorphism of~$\quartic{64}$.
Hence, such a line is real if and only if $\Gs_g$ preserves its pair of
endpoints;
there are four such lines for any~$g$.
The other lines constitute the orbit $\Aut_0\quartic{64}/G$ of~$l_0$,
where $G=\bG/H$ is the diagonal. Since $\bar l_0=l_0$, a line $h l_0$ is
$\Gs_g$-real
if and only if $h\1g\bar h\in G$.
Now, both statements are easily proved using
\GAP~\cite{GAP4}.
\endproof

\subsection{A real quartic with $56$ real lines}\label{s.56}
To our knowledge, this example is new.

Below, we make use of bihomogeneous polynomials, \ie, algebraic curves in the
product $\Cp1\times\Cp1$. For the sake of simplicity, we use the affine
coordinates $u:=z_0/z_1$, $v:=z_2/z_3$ in the two copies of~$\Cp1$.

Fix
$\e:=\pm\sqrt2$ and consider the polynomials
\[*
\poly_1(u,v):=-3v+v^3+2\e u,\qquad
\poly_2(u,v):=2\e u^3-v+3u^2v
\]
of bidegree $(1,3)$ and $(3,1)$, respectively.
The quartic $Y:=\quartic{56.real}$ in question is given by the polynomial
\[
z_1z_3^3\poly_1\Bigl(\frac{z_0}{z_1}, \frac{z_2}{z_3}\Bigr)
 -z_1^3z_3\poly_2\Bigl(\frac{z_0}{z_1}, \frac{z_2}{z_3}\Bigr),
\label{eq.56}
\]
or, explicitly,
\[*
3\e z_0^2z_1z_2+3\e z_1z_2z_3^2-\e z_1^3z_2-\e z_1z_2^3+4z_0^3z_3-4z_0z_3^3.
\]
Below, we show that $Y$ contains $56$ real lines; by \autoref{th.unique}, this
configuration of lines is \config{56.real}, and $Y$ is the only real quartic
with this property.

The quartic~$Y$ contains the two lines
\[
m_1:=\{z_0=z_1=0\},\quad
m_2:=\{z_2=z_3=0\}.
\label{eq.m}
\]
The curves $\{\poly_k=0\}\subset\Cp1\times\Cp1$, $k=1,2$, intersect at $10$
real points, see \autoref{tab.points}.
\table
\caption{The solutions to $\poly_1(u,v)=\poly_2(u,v)=0$}\label{tab.points}
\centerline{\vbox{\openup2pt%
\halign{\strut\ $#$\hss&&\qquad$#$\hss\ \cr
\noalign{\hrule\vspace{2pt}}
P_1(-1+\e, -1+\e)&A_1(1/\e, -2)&B_1(\infty,\infty)\cr
P_2(1+\e, -1-\e)&A_2(1/2, \e)&B_2(0,0)\cr
P_3(1-\e, 1-\e)&C_1(-1/\e, 2)\cr
P_4(-1-\e, 1+\e)&C_2(-1/2, -\e)\cr
\noalign{\vspace{2pt}\hrule}\crcr}}}
\endtable
Each such point~$L(u,v)$ gives rise to the line
\[
l:=\{z_0=uz_1,\ z_2=vz_3\}
\label{eq.l}
\]
through $[u:1:0:0]\in m_2$ and $[0:0:v:1]\in m_1$; it is
contained in~$Y$.

The intersection of~$Y$ with each of the six planes
shown in \autoref{tab.planes} splits into four
lines; twelve of the resulting $24$ lines (some of which coincide)
are among~\eqref{eq.m} and~\eqref{eq.l}, see \autoref{tab.planes},
and the twelve others are new and distinct.

\table
\caption{The six special planes}\label{tab.planes}
\hbox to\hsize{\hss\vbox{\openup2pt%
\halign{\strut\ $\deMaple#$:\hss&&\qquad\deMaple#\hss\cr
\noalign{\hrule\vspace{2pt}}
\omit\strut\ Plane\hss&New lines&Old lines\ \cr
\noalign{\hrule\vspace{2pt}}
z[1]=0&$(z[0]-z[3])*(z[0]+z[3])$&$m_1$, $b_1$\cr
z[2]=0&$(z[0]-z[3])*(z[0]+z[3])$&$m_2$, $b_2$\cr
z[1]=\e*z[0]&$(z[0]+z[2]-z[3])*(z[0]-z[2]+z[3])$&$m_1$, $a_1$\cr
z[1]=-\e*z[0]&$(z[0]+z[2]+z[3])*(z[0]-z[2]-z[3])$&$m_1$, $c_1$\cr
z[2]=\e*z[3]&$(z[0]+z[1]+z[3])*(z[0]+z[1]-z[3])$&$m_2$, $a_2$\cr
z[2]=-\e*z[3]&$(z[0]-z[1]-z[3])*(z[0]-z[1]+z[3])$&$m_2$, $c_2$\cr
\noalign{\vspace{2pt}\hrule}\crcr}}\hss}
\endtable


Finally, the ten skew lines~\eqref{eq.l} constitute sixteen quadruples
$(l_1,l_2,l_3,l_4)$, each lying in a quadric, see
\autoref{tab.quadrics}. The equation of this quadric~$Q$ is
\[*
z_1z_3\,\chi\Bigl(\frac{z_0}{z_1}, \frac{z_2}{z_3}\Bigr)=0,
\]
where $\chi(u,v)$ is the polynomial given in \autoref{tab.quadrics}
(see also \autoref{rem.quadrics} below).
The intersection $Y\cap Q$ is a bidegree~$(4,4)$ curve in~$Q$. Since it
contains four skew generatrices of~$Q$, it must split into $l_1,\ldots,l_4$
and four generatrices of the other family.
Two of them are $m_1$, $m_2$, and the two others are new.
It is straightforward
that the $16\times2=32$ lines thus obtained are all
real (see also \autoref{rem.real} below),
pairwise distinct (as the sixteen quadrics
are distinct),
and distinct from~\eqref{eq.m}, \eqref{eq.l}, and
the lines in \autoref{tab.planes}
(as they are disjoint from $m_1\cup m_2$).

\table
\caption{The sixteen special quadrics}\label{tab.quadrics}
\hbox to\hsize{\hss\vbox{\openup2pt%
\halign{\strut\ \lowercase{$#$}:\hss&&\qquad$#$\ \hss\cr
\noalign{\hrule\vspace{2pt}}
(A_1,B_1,C_1,B_2)&u+2\e v\cr
(A_2,B_2,C_2,B_1)&u-2\e v\cr
(P_1,P_3,B_1,B_2)&u-v\cr
(P_2,P_4,B_1,B_2)&u+v\cr
(P_3,P_4,B_1,A_1)&1+\e u+v\cr
(P_1,P_2,B_1,C_1)&1-\e u-v\cr
(P_2,P_3,B_2,A_2)&\e u-v+uv\cr
(P_1,P_4,B_2,C_2)&\e u-v-uv\cr
(P_1,P_4,A_1,C_1)&\e-2\e u-v+uv\cr
(P_2,P_3,A_1,C_1)&\e+2\e u+v+uv\cr
(P_1,P_2,A_2,C_2)&1-2\e u+v-\e uv\cr
(P_3,P_4,A_2,C_2)&1+2\e u-v-\e uv\cr
(P_1,P_3,A_1,A_2)&-3\e+4+(2\e-2)u-(2\e-2)v+\e uv\cr
(P_1,P_3,C_1,C_2)&-3\e+4-(2\e-2)u+(2\e-2)v+\e uv\cr
(P_2,P_4,A_1,C_2)&3\e+4+(2\e+2)u+(2\e+2)v+\e uv\cr
(P_2,P_4,A_2,C_1)&3\e+4-(2\e+2)u-(2\e+2)v+\e uv\cr
\noalign{\vspace{2pt}\hrule}\crcr}}\hss}
\endtable

Summarizing, we obtain $2+10+12+32=56$ real lines in~$Y$.

\remark\label{rem.quadrics}
Let $u_1,\ldots,u_4\in m_1$ and $v_1,\ldots,v_4\in m_2$ be two quadruples,
where, as above, we let $u:=z_0/z_1$ and $v:=z_2/z_3$.
Then the lines $l_i:=(u_iv_i)$, $i=1,\ldots,4$, \cf.~\eqref{eq.l},
lie in a quadric if and only if the quadruples $(u_i)$ and $(v_i)$ are
isomorphic, \ie, their cross-ratios are equal.
When this is the case,
the quadruples are related by a fractional linear transformation,
$v_i=f(u_i)$ for $i=1,\ldots,4$, and
the equation of the quadric is obtained from
$z_2/z_3=f(z_1/z_0)$
by clearing the denominators.
\endremark

\subsection{Further properties of $\quartic{56.real}$}\label{s.56.more}
Let $Y:=\quartic{56.real}$ be the quartic constructed in the previous section.
The following statements are straightforward.
\roster
\item
The lines~$m_1$ and~$m_2$ are disjoint.
\item
The
lines~\eqref{eq.l} are pairwise disjoint; each of them
intersects~$m_1$ and~$m_2$.
\endroster
Let~$\Ga$ be a plane as in
\autoref{tab.planes}. Then
$Y\cap\Ga$ splits into $m_i$, $l$, and a pair $r_1$, $r_2$,
where $i=1$ or~$2$ and $l$ is one of the lines~\eqref{eq.l},
see \autoref{tab.planes}.
\roster[3]
\item
The lines~$r_1$ and~$r_2$ intersect~$m_i$, $l$, and each other;
they are disjoint from $m_{3-i}$ and any line~$l'\ne l$ as in~\eqref{eq.l}.
\endroster
This observation confirms the fact that all twelve lines thus obtained are
pairwise distinct and distinct from~\eqref{eq.m} and~\eqref{eq.l}.
Note that, according to \autoref{tab.planes}, the plane~$\Ga$ is completely
determined by the line $l\subset\Ga$ as in~\eqref{eq.l};
hence, we can use the notation $\Ga(l)$ and $r_{1,2}(l)$.

Finally, pick a quadruple $(l_1,l_2,l_3,l_4)$ as in \autoref{tab.quadrics},
let~$Q$ be the corresponding quadric,
and let $n_1$, $n_2$ be the two extra lines (other
than $m_1$, $m_2$) in $Y\cap Q$.
The remaining observations follow from the properties of the generatrices
of~$Q$; in particular,
the intersection $Y\cap Q$ may contain at most four generatrices of
each family and, if a line intersects three generatrices of the same family,
it lies in~$Q$.
\roster[4]
\item
The lines~$n_1$ and~$n_2$ are disjoint from~$m_1\cup m_2$; they intersect
each of $l_1,l_2,l_3,l_4$ and are disjoint from all other lines as
in~\eqref{eq.l}.
\item\label{intr.5}
If a line $l$ as in~\eqref{eq.l} is distinct from all~$l_i$, $i=1,\ldots,4$, the
lines~$n_1$, $n_2$ and $r_{1,2}(l)$
can be indexed so that $\fb(n_i\cap r_j)=\Gd_{ij}$ is the Kronecker symbol.
\endroster
In more details, the intersection matrix can be computed using explicit
equations of all lines. We leave this exercise to the reader.

\remark\label{rem.real}
Statement~\iref{intr.5} proves also that $n_1$ and~$n_2$ are real:
if they were complex conjugate, they would have to intersect the same
real line~$r_1$ or~$r_2$.
\endremark

We conclude with a description of the automorphism group $\Aut \quartic{56.real}$.

\proposition\label{prop.Aut56}
The group $\Aut \quartic{56.real}\subset\PGL(4,\C)$
is generated by
\roster*
\item
the reflections $z_i\mapsto\rho_iz_i$ with
$\rho_i=\pm1$ and $\rho_0\rho_3=\rho_1\rho_2$,
\item
the transposition $z_1\leftrightarrow z_2$,
\item
the order~$4$ map $z_0\mapsto z_3$, $z_3\mapsto-z_0$, and
\item
the involution $z_0\mapsto(z_0+z_3)/\e$, $z_3\mapsto(z_0-z_3)/\e$.
\endroster
This group has order~$32$\rom;
it acts faithfully on the set of lines contained in~$\quartic{56.real}$.
\endproposition

\proof
Computing the intersection matrix, one can see that there are exactly four
pairs $(l_1,l_2)$ of skew lines
such that $l_1$ and~$l_2$ intersect ten other common lines.
In turn, these pairs split
into four quadrangles: one is $(m_1,m_2)$, $(b_1,b_2)$, and the other is
formed by the four remaining lines in the planes $\{z_1=0\}$ and $\{z_2=0\}$,
see \autoref{tab.planes}.
The last involution in the statement interchanges the two quadrangles. The
other transformations preserve the quadrangle $(m_1,m_2)$, $(b_1,b_2)$ and,
hence, the coordinate tetrahedron; they can easily be listed.
The last two statements are straightforward.
\endproof

\remark\label{rem.Y.aut}
All automorphisms of~$\quartic{56.real}$ are real with respect to the standard
complex conjugation $c\:[z_i]\mapsto[\bar z_i]$. Hence, the last statement of
\autoref{prop.Aut56} implies that $c$ is the only real structure on~$\quartic{56.real}$
with respect to which all $56$ lines are real.
(In fact, up to automorphism $\quartic{56.real}$ has six real structures: they are
enumerated by the conjugacy classes of the involutions in $\Aut \quartic{56.real}$.)
\endremark

\remark\label{rem.Y.rational}
By
rescaling $u\mapsto\e u$, one can make
the quartic~\quartic{56.real} defined
over~$\Q$; however, some of the lines are still defined over the quadratic
number field $\Q(\e)$ only. To see this, one can observe that the
cross-rations of some of the quadruples of points in $m_1$ cut by the lines
as in~\eqref{eq.l} are irrational, see \autoref{tab.points}.
\endremark

\subsection{A few other quartics}\label{s.60}
In
this concluding section, we describe
briefly a few other quartics with large configurations of lines,
for which we do not know explicit equations.
The existence (and uniqueness, when it holds) is given by the existence of
the corresponding $\L$-configurations, see \autoref{tab.list},
and the results of \autoref{arithmetical_reduction}.
Other properties, \eg, groups of projective automorphisms, classes of real
structures, \etc., can easily be computed using the corresponding properties
of configurations
and Nikulin's theory of lattice extensions; however, we
omit these details.

\subsubsection{The quartics menioned in \autoref{th.unique}}\label{ss.th.unique}
By \autoref{lem.unique} and \autoref{th.complex}, for each
of the four configurations
$S=\config{60}$, \config{q56}, \config{54}, \config{q54},
there
exists
a unique, up to projective equivalence, quartic~$X$ such that
$\Fano(X)\cong S$;
this
quartic can be chosen real, see
\autoref{symmetric-reflexive}.
We denote these quartics by
\quartic{60}, \quartic{q56}, \quartic{54}, \quartic{q54}, respectively.
Besides, for $S=\config{60.2}$ or \config{56}, there is a unique pair of
nonequivalent
complex conjugate quartics $X$, $\bar X$ such that
$\Fano(X)\cong\Fano(\bar X)\cong S$; these pairs are denoted by
\quartic{60.2}, \bquartic{60.2} and
\quartic{56}, \bquartic{56}, respectively.
Together with \quartic{64} (see~\autoref{s.Schur}) and \quartic{56.real}
(see \autoref{s.56} and \autoref{s.56.more}), these
surfaces make a complete list of
quartics containing more than $52$ lines.

\table
\caption{Configurations with many $\L$-realizations}\label{tab.Y}
\def\*{\raise3pt\hbox{$\scriptstyle\sqrt{}$}}
\hbox to\hsize{\hss\vbox{\halign{\strut\quad#\quad\hss&\hss#\hss\quad&
\hss#\hss\quad&\hss#\hss\quad&\hss#\hss\quad&&\hss$#$\hss\quad\cr
\noalign{\hrule}%
 \vphantom{\Big(}\hss$S$\hss&$\ls|\Fn|$&t.r.&ref&sym&\ls|\OG_h(S)|&
 \discr S&T:=S^\perp&(r,c)
\cr
\noalign{\vspace{-1pt}\hrule\vspace{2pt}}
\config{60.2}&60&&&&240&
\<\frac65\>\oplus\<\frac{10}{11}\>&[4,1,14]&(0,1)
\cr
\config{56}&56&&&&128&%
\<\frac{15}8\>\oplus\<\frac{15}8\>&[8,0,8]&(0,1)
\cr
\config{52.5}&52&\*&\*&\*&8&\<\frac12\>\oplus\<\frac32\>\oplus\<\frac4{19}\>
        &[2,0,38]&(1,1)\cr
  &&&&&&&[8,2,10]
\cr
\config{52.real}&52&\*&\*&\*&8&\<\frac6{79}\>
        &[2,1,40]&(1,2)\cr
  &&&&&&&[4,1,20]\cr
  &&&&&&&[8,1,10]
\cr
\config{q52.2}&52&&&&64&\<\frac14\>\oplus\<\frac54\>\oplus\<\frac25\>
        &[8,4,12]&(0,1)
\cr
\config{51}&51&&\*&\*&12&
        \<\frac43\>\oplus\<\frac2{29}\>
        &[6,3,16]&(1,1)\cr
  &&&&&&&[4,1,22]
\cr
\config{50.2}&50&&$\times2$&$\times2$&12&
        \<\frac74\>\oplus\<\frac58\>\oplus\<\frac43\>
        &[4,0,24]&(2,0)
\cr
\config{50}&50&&&&16&
        \<\frac74\>\oplus\<\frac58\>\oplus\<\frac43\>
        &[4,0,24]&(0,1)
\cr
\config{48}&48&\*&\*&\*&8&\<\frac12\>\oplus\<\frac5{16}\>\oplus\<\frac23\>
        &[2,0,48]&(1,0)
\cr
\config{48.2}&48&\*&\*&\*&8&\<\frac25\>\oplus\<\frac4{19}\>
        &[2,1,48]&(2,1)\cr
  &&&&&&&[8,1,12]\cr
  &&&\*&\*&&&[10,5,12]
\cr
\noalign{\vspace{2pt}\hrule}\crcr}}\hss}
\endtable

\subsubsection{Large configurations of real lines}\label{ss.real}
Arguing as in the proof of \autoref{lem.unique}, it is not difficult to
classify the $\L$-realizations of the four other \config{Y}-type
configurations listed in \autoref{tab.list}; we summarize the results in
\autoref{tab.Y}. This table is organized similar to \autoref{tab.large},
with the
last column showing the numbers~$r$, $c$ of, respectively, real and pairs of
complex conjugate quartics with the given configuration of lines.
Note, though, that, with the only exception of \autoref{prop.unique.real}
below, we never assert the uniqueness of the real form: considering the large
automorphism groups, it is likely not unique, \cf. \autoref{prop.X64} and
\autoref{rem.Y.aut}.

If $S=\config{48.2}$, the natural homomorphism
$\rho\:\OG_h(S)\to\Aut\discr S$ maps $\OG_h(S)$ onto the index~$2$ subgroup
$\Aut\discr_5S$; in the other three cases, $\rho$ is surjective.
It follows that,
in all four cases, the weak isomorphism classes of $\L$-realizations are
classified by the transcendental lattices $T:=S^\perp$.
In three cases, there are several isomorphism classes;
however, only one of them is totally reflexive. In view of
\autoref{app.real},
this fact merits a separate statement.

\proposition\label{prop.unique.real}
For each \config{Y}-type configuration~$S$ listed in \autoref{tab.list},
there  is a unique, up to real projective equivalence, real quartic
$Y$ such that $\Fano_\R(Y)\cong S$.
The real part of this real quartic is a connected surface of genus~$10$.
\done
\endproposition

The configuration $S=\config{48.2}$ admits another reflexive
$\L$-realization, which is not totally reflexive;
thus, the corresponding quartic~$X$ can be chosen real, but some of the
lines contained in~$X$ are necessarily complex conjugate.
(Note that, unlike
the case of
Schur's quartic~\quartic{64}, see \autoref{prop.X64},
or the maximizing real quartic \quartic{56.real}, see
\autoref{rem.Y.aut},
this quartic~$X$ and the quartic~$Y$ given by \autoref{prop.unique.real}
are not just distinct real forms of the same surface: $X$ and $Y$ are not
projectively equivalent even over~$\C$.)

\remark\label{rem.min.field}
\autoref{tab.Y} suggests also that the quartics $\quartic{Y}_*$ realizing
each of the configurations $\config{Y}_*:=\config{52.5}$, \config{52.real},
\config{48}, \config{48.2} should be Galois conjugate over a certain
algebraic number field~$\Bbb K$ of degree~$3$, $5$, $1$, $4$,
respectively, so that this field~$\Bbb K$ is the minimal field of definition
of~$\quartic{Y}_*$. In particular, \quartic{52.5} and \quartic{52.real} are
probably not defined over~$\Q$, \cf. open problems at the end of
\autoref{s.results}.
\endremark

\subsubsection{Configurations with many realizations}\label{ss.many}
For completeness, we describe also the few configurations from
\autoref{tab.list} that admit more than one geometric realization.
For all configurations of maximal rank, the computation runs exactly as in
the proof of \autoref{lem.unique} and can easily be automated.
Omitting the straightforward details, we summarize the results in
\autoref{tab.Y}.
(The meaning of the columns is explained in \autoref{ss.real}.)
In the line containing~\config{50.2}, the symbol ``$\times2$'' means
that there are two distinct
geometric $\L$-realizations, which are both reflexive.

Any other configuration in \autoref{tab.list} admits a unique geometric
$\L$-realization, and this realization is reflexive.
In particular, for a configuration~$S$ as in
\autoref{tab.list}, a
geometric $\L$-realization is reflexive if and only if
it is symmetric. Currently, we do not know whether this is a common property
of configurations: in some similar $K3$-related problems, it may not hold
(\cf. the existence of a connected real equisingular family of simple
plane sextics containing no real curves discovered in~\cite{degt:geography}.)

\remark\label{rem.abstract.K3}
The isomorphism type of a singular $K3$-surface (\ie, one of Picard
rank~$20$) is determined by its transcendental lattice. Analyzing
Tables~\ref{tab.large} and~\ref{tab.Y}, one can observe that the quartics
\quartic{60} and \quartic{q56} are isomorphic as abstract $K3$-surfaces;
a similar statement holds for
the seven quartics realizing the configurations
\config{54}, \config{50.2}, and \config{50}.
On the other hand, each of the configurations \config{52.5},
\config{52.real}, \config{51}, \config{48.2} is realized by several distinct
$K3$-surfaces.
\endremark

\subsubsection{Families with parameters}\label{ss.Z}
Finally, worth mentioning are the
configurations $S=\config{52.0}$, \config{50.0}, \config{49}
in \autoref{tab.list}.
Recall
that the dimension of the equilinear moduli space
$\Omega'(S)/\!\PGL(4, \C)$,
\cf. the proof of \autoref{th.complex},
equals $20-\rank S$; hence, we obtain $1$-parameter families of
distinct quartics sharing the same combinatorial type of configurations of
lines.

The connectedness of each family follows from \autoref{th.complex} and
a
computation based on
the results of \cite{Miranda.Morrison:book}, covering
indefinite transcendental lattices.
We have
\roster*
\item
if $S=\config{52.0}$, then
$\CS=\<\frac12\>\oplus\<\frac12\>\oplus\<\frac58\>\oplus\<\frac43\>$
and $\OG^+(T)\onto\Aut\discr T$;
\item
if $S=\config{50.0}$, then
$\CS=\<\frac74\>\oplus\<\frac25\>\oplus\<\frac25\>$
and $\OG_h(S)\onto\Aut\CS$;
\item
if $S=\config{49}$, then
$\CS=\CV_2\oplus\<\frac54\>\oplus\<\frac67\>$
and
$\Im[\OG_h(S)\to\Aut\CS]=\Aut\CS_2$.
\endroster
The uniqueness of $T:=S^\perp$ in its genus and the assertion on $\OG^+(T)$
for $S=\config{52.0}$ follow from~\cite{Miranda.Morrison:book}.
Thus, in each case, there is a unique geometric $\L$-realization.
If $S=\config{52.0}$, this realization is totally reflexive,
\ie,
there is a
$1$-parameter family (not necessarily
connected) of real quartics~$Z$ such
that $\Fano_\R(Z)\cong\config{52.0}$.
For the other two configurations, for each involution
$a\in\Aut\discr S$, exactly one of $\pm a$
admits an involutive lift to $\OG_h(S)$. Hence,
these configurations are reflexive
(not totally)
and the corresponding equilinear families also
contain real quartics.

The existence of the family corresponding to \config{52.0}, with
exactly $52$ lines
in each quartic, as well as the non-uniqueness
of $\L$-realizations discussed in \autoref{ss.real} and \autoref{ss.many},
can be regarded as yet another justification for the
assumption $\ls|\Fn(X)|>52$
in \autoref{th.unique}: quartics with fewer lines are probably
more difficult to control.


{
\let\.\DOTaccent
\def\cprime{$'$}
\bibliographystyle{amsplain}
\bibliography{degt,sinan}
}

\end{document}